\RequirePackage{lineno}
\documentclass[11pt]{article}
\usepackage{float}
\usepackage{hyperref}
\usepackage[T1]{fontenc}
\usepackage[latin1]{inputenc}
\usepackage{amsmath,amsthm,amssymb}
\usepackage{mathabx}
\usepackage{xcolor}
\usepackage{authblk}
\usepackage[margin=1in]{geometry}
\usepackage{subcaption} 
\usepackage{mathrsfs}

\usepackage{graphicx}
\usepackage{amsbsy}
\usepackage{yfonts}
\usepackage[all]{xy}
\usepackage{pgfplots}

\providecommand{\U}[1]{\protect\rule{.1in}{.1in}}

\newtheorem{prop}{Proposition}[section]
\newtheorem{cor}[prop]{Corollary}
\newtheorem{defi}[prop]{Definition}

\newtheorem{lem}[prop]{Lemma}

\newtheorem{theo}[prop]{Theorem}

\newcommand{\vertiii}[1]{{\left\vert\kern-0.25ex\left\vert\kern-0.25ex\left\vert #1
    \right\vert\kern-0.25ex\right\vert\kern-0.25ex\right\vert}}

\newcommand{\AC}{\mathbb{A}}
\newcommand{\CC}{\mathbb{C}}

\newcommand{\EE}{\mathbb{E}}

\newcommand{\KK}{\mathbb{K}}
\newcommand{\LL}{\mathbb{L}}
\newcommand{\MM}{\mathbb{M}}
\newcommand{\NN}{\mathbb{N}}
\newcommand{\PP}{\mathbb{P}}
\newcommand{\QQ}{\mathbb{Q}}
\newcommand{\RR}{\mathbb{R}}

\newcommand{\XX}{\mathbb{X}}

\newcommand{\YY}{\mathbb{Y}}
\newcommand{\ZZ}{\mathbb{Z}}
\newcommand{\GG}{\mathbb{G}}

\newcommand{\Ba}{ {\cal B }}

\newcommand{\Da}{ {\cal D }}
\newcommand{\La}{ {\cal L }}

\newcommand{\Ea}{ {\cal E }}

\newcommand{\Va}{ {\cal V }}
\newcommand{\Ua}{ {\cal U }}
\newcommand{\Fa}{ {\cal F }}
\newcommand{\Ga}{ {\cal G }}

\newcommand{\Ia}{ {\cal I }}
\newcommand{\Xa}{ {\cal X }}
\newcommand{\Ma}{ {\cal M }}
\newcommand{\Ta}{ {\cal T}}
\newcommand{\Ha}{ {\cal H }}
\newcommand{\Ja}{ {\cal J }}
\newcommand{\Pa}{ {\cal P }}
\newcommand{\Za}{ {\cal Z }}

\newcommand{\point}{\mbox{\LARGE .}}

\newcommand{\cqfd}{\hfill\blbx \\}
\def\blbx{\hbox{\vrule height 5pt width 5pt depth 0pt}\medskip}
\def \KK{\mathbb{K}}
\def \PP{\mathbb{P}}
\def \RR{\mathbb{R}}
\def \SS{\mathbb{S}}
\def \EE{\mathbb{E}}
\def \QQ{\mathbb{Q}}

\def \CC{\mathbb{C}}
\def \LL{\mathbb{L}}
\def \ZZ{\mathbb{Z}}

\def \BB{\mathbb{B}}
\def \UU{\mathbb{U}}

\newcommand{\goodchi}{\protect\raisebox{2pt}{$\chi$}}

\numberwithin{equation}{section}

\makeatletter
\DeclareRobustCommand\frownotimes{\mathbin{\mathpalette\frown@otimes\relax}}
\newcommand{\frown@otimes}[2]{
  \vbox{
    \ialign{##\cr
      \hidewidth$\m@th#1{}_\frown$\kern-\scriptspace\hidewidth\cr
      \noalign{\nointerlineskip\kern-.1pt}
      $\m@th#1\otimes$\cr
    }
  }
}
\makeatother

\begin{document}

\title{A duality formula and a particle Gibbs sampler for continuous time Feynman-Kac measures on path spaces}

\author[$1$]{Marc Arnaudon}
\author[$2$]{Pierre Del Moral\thanks{P. Del Moral was supported in part by the Chair Stress Test, RISK Management and Financial Steering, led by the French Ecole Polytechnique and its Foundation and sponsored by BNP Paribas.}}
\affil[$1$]{{\small Univ. Bordeaux, CNRS, Bordeaux INP, IMB, UMR 5251, F.33400, Talence, France}}
\affil[$2$]{{\small INRIA, Bordeaux Research Center, Talence  \& CMAP, Polytechnique Palaiseau, France}}
\date{}

\maketitle

\begin{abstract}
Continuous time Feynman-Kac measures on path spaces are central in applied probability, partial differential equation theory, as well as in quantum physics.
This article presents a new duality formula between normalized Feynman-Kac distribution and their mean field particle interpretations.
Among others, this formula allows us to design a reversible particle Gibbs-Glauber sampler for continuous time
Feynman-Kac  integration on path spaces. We also provide new   Dyson-Phillips semigroup  expansions, as well as novel uniform propagation of chaos estimates for continuous time genealogical tree based particle models with respect to the time horizon and the size of the systems. Our approach is self contained and it is based on a novel stochastic perturbation analysis and backward  semigroup techniques.
These techniques  allow to obtain sharp quantitative estimates of
the convergence rate to equilibrium of particle Gibbs-Glauber samplers. To the best of our knowledge these
results are the first of this kind for continuous time Feynman-Kac measures.\\

\emph{Keywords} : Feynman-Kac formulae, interacting particle systems, genealogical trees, ancestral lines, Gibb-Glauber dynamics, propagation of chaos properties,  contraction inequalities, Dyson-Phillips expansions.\\

\emph{Mathematics Subject Classification} : 60K35, 60H35, 37L05, 47D08.

\end{abstract}

\section{Introduction}

Feynman-Kac measures on path spaces are central in applied probability as well as in biology and quantum physics. They also arise 
in a variety of application domains such as in estimation and control theory, as well as a rare event analysis. For a detailed review 
on Feynman-Kac  measures and their application domains we refer to the books~\cite{d-2004,d-2013,dm-2000,dm-penev}, 
see also the more recent articles~\cite{cox,kyprianou} on branching processes and neutron transport equations and the references therein.

Their mean field type particle interpretations are defined as a system of particles 
 jumping a given rate uniformly onto the population. From the pure numerical viewpoint, this interacting jump transition can be interpreted as an acceptance-rejection
 scheme with a recycling of rejected particles by duplicating the selected ones. Feynman-Kac interacting particle models encapsulate a variety of algorithms such as
 the diffusion Monte Carlo used to solve Schr\" odinger ground states, see for instance the series of articles~\cite{caffarel,cances-tony,dm-sch,mathias,mathias-2,tony,tony-rs} as well as section 23.5 and chapter 27 in~\cite{dm-penev} and the references therein. 
 
 Their discrete time versions are encapsulated a variety of well known algorithms such as particle filters~\cite{d-1996} (a.k.a. sequential Monte Carlo methods
 in Bayesian literature~\cite{cappe-moulines,d-2004,d-2013,dm-2000,ddg}),  the go-with the winner~\cite{aldous},
 as well as the self-avoidind random walk pruned-enrichment algorithm by Rosenbluth and Rosenbluth~\cite{rosen}, and many others. This list is not exhaustive (see also the references therein).
 The research monographs~\cite{d-2004,d-2013} provide a detailed discussion on these subjects with precise reference pointers.
 
The seminal article~\cite{adh-2010} by Andrieu, Doucet and Holenstein  introduced a new way to combine Markov chain Monte Carlo methods  with
discrete generation particle methods. A variant of
the method, where ancestors are resampled in a forward pass, was developed by Lindsten,
 Sch\"on and Jordan in~\cite{lindsten}, and Lindsten and Sch\"on~\cite{lindsten2}. In all of these studies, the validity of
the particle conditional sampler is assessed by interpreting the model as a traditional Markov chain Monte Carlo sampler
on an extended state space. The central idea is first to 
design a detailed  encoding of the ancestors at each level in terms of random maps on integers, and then
to extend the "target" measure on a sophisticated state space incapsulating these iterated random sequences. 
 
 In a more recent article~\cite{DMKP:16}, the authors provide an alternative and we believe more 
 natural interpretation of these particle Markov chain Monte Carlo methods in terms of a duality formula extending the well known unbiasedness properties
of Feynman-Kac particle measures on many-body particle measures. This article also provides sharp quantitative estimates of the convergence rate to equilibrium of the models
with respect to the time horizon and the size of the systems. The analysis of these models, including backward particle Markov chain Monte Carlo samplers  has been further developed in~\cite{dmj-18-1,dmj-18-2}.

The main objective of the present article is to extend these methodologies to continuous time Feynman-Kac measures on path spaces.

The first difficulty comes from the fact that the discrete time analysis~\cite{dmj-18-1,dmj-18-2,DMKP:16} only applies to simple
genetic type particle models, or equivalently to branching models with pure multinomial selection schemes. Thus, these results do not apply to discrete time approximation of continuous time models based on geometric type jumps, and any density type argument cannot be applied.

  In contrast with their discrete time version, continuous time
Feynman-Kac particle models are not described by conditionally independent local transitions, but in terms of 
 interacting jump processes. 

The analysis of continuous time genetic type particle models is not so developed as their discrete time versions. For instance, uniform convergence estimates are available for continuous time   Feynman-Kac models  with stable processes~\cite{dm-2000,dm-stab,dm-sch,mathias}. Nevertheless, to the best of our knowledge,
 sharp estimates for path space models and genealogical tree based particle samplers in continuous time have never been discussed in the literature.
 These questions are central in the study the convergence to equilibrium of particle Gibbs-Glauber sampler on path spaces. 

In the present article we provide a duality formula for continuous time Feynman-Kac measures on path-spaces (cf. theorem~\ref{theo-duality}). This formula on generalogical tree based particle models that can be seen
as an extension of well known unbiasedness properties of Feynman-Kac models to their many body version (defined in section~\ref{sec-many-body-FK}). The second main result of the article is to design and to analyze the stability properties of a particle Gibbs-Glauber sampler of  path space (cf. theorem~\ref{theo-PG-intro}). Our approach combines a perturbation analysis of nonlinear stochastic semigroups with propagation of chaos techniques (cf. section~\ref{sec-perturbation}). Incidentally these techniques also provide with little efforts new uniform propagation of chaos estimates w.r.t. the time horizon (cf. corollary~\ref{cor-pchaos}).

\subsection{Statement of the main results}\label{sec-1-ref}
Let $(X_t,V_t)$ be a continuous time Markov process and a bounded non negative function on some metric space $(S,d_S)$.
We denote by $D_t(S)$ the set  
of c\`adl\`ag paths  from $[0,t]$ to $S$.
Unless otherwise is stated, as a rule in the further development of the article $$\widehat{X}_t:=(X_s)_{s\leq t}\in \widehat{S}:=\cup_{t\geq 0}~D_t(S)$$ stands for the historical process of some process $X_t$. In this notation, we extend $V_t$ to $D_t(S)$
 by setting $\widehat{V}_t(\widehat{X}_t)=V_t(X_t)$. 
 
The Feynman-Kac probability measures  $\QQ_t$ on the path spaces $D_t(S)$  associated with  $(X_t,V_t)$ are defined for any bounded measurable function
 $F$ on $D_t(S)$ by the path integration formula
 \begin{equation}\label{def-FK-intro}
 \QQ_t(F):=\int~F(\omega)~\QQ_t(d\omega):= \ZZ_t^{-1}~\EE\left(F(\widehat{X}_t)~ \exp{\left[-\int_0^t\widehat{V}_s(\widehat{X}_s)ds\right]}\right)
 \end{equation}
where $\ZZ_t$ stands for the normalizing constant. 

We also let $\eta_t$ be the terminal time marginal of the measures $\QQ_t$. In this notation, for any bounded measurable function
 $f$ on $S$ we have
 \begin{equation}\label{def-eta-intro}
 \eta_t(f):=\int~f(x)~\eta_t(dx):= \ZZ_t^{-1}~\EE\left(f(X_t)~ \exp{\left[-\int_0^tV_s(X_s)ds\right]}\right)
 \end{equation}
In addition, the normalizing constant $\ZZ_t$ is given by the easily checked free energy formula
 $$
\ZZ_t=\EE\left(\exp{\left[-\int_0^tV_s(X_s)ds\right]}\right)=\exp{\left[-\int_0^t\eta_s(V_s)ds\right]}
$$
For a more thorough discussion on these Feynman-Kac models we refer the reader to section~\ref{ref-def-D}
and section~\ref{sec-historical-FK}. A selected list of application areas are also discussed in section~\ref{sec-appli-FK}.
 
 The path integration formulae (\ref{def-FK-intro}) can rarely be solved analytically and their numerical solving often require  extensive calculations.
 One strategy is to interpret these probability  measures on path space in terms the occupation measures of
the ancestral lines of a genetic type interacting jump process~\cite{dm-2000-moran,dm-2000,dm-sch,dm-penev}.  

These particle interpretations are defined as follows:

 Consider a system of $N$ particles $\xi_t=(\xi^i_t)_{1\leq i\leq N}\in S^N$ evolving independently as $X_t$ with jump rate $V_t$;
at each jump time the particle jumps onto a particle uniformly chosen in the pool. 

This jump can be interpreted 
as the death of the particle and the instantaneous birth of an offspring of a particle uniformly chosen in the pool.
After this birth each offspring evolves as independent copies of the process $X_t$. 

When the $i$-th  particle duplicates it becomes the parent of two offsprings. Running back in time we can trace the whole genealogy of each particle. 

Let $A_{s,t}(i)\in [N]$ be the index of the ancestor of of the $i$-th particle $\xi^i_t$ at level $s\leq t$. In this notation, the $i$-th ancestral line $\XX^i_t:=(\xi^{A_{s,t}(i)}_s)_{0\leq s\leq t}\in D_t(S)$ of the $i$-th particle $\xi^i_t$ at time $t$ is represented in a synthetic way by the following backward ancestral path 
$$
\xi^{A_{0,t}(i)}_0\leftarrow\ldots\leftarrow\xi^{A_{s,t}(i)}_s\leftarrow\ldots\leftarrow\xi^i_t
$$
In this interpretation, the genealogical tree  of the $N$ particles $\xi_t=(\xi^i_t)_{1\leq i\leq N}$ is defined by  the $N$ ancestral lines $\XX_t=(\XX^i_t)_{1\leq i\leq N}\in D_t(S)^N$.

We underline that the $N$ ancestral lines $\XX_t=(\XX^i_t)_{1\leq i\leq N}\in D_t(S)^N$ of length $t$ of the $N$ individuals  $\xi_t=(\xi^i_t)_{1\leq i\leq N}$ is also defined forward in time by a 
system of $N$ path-valued particles evolving independently as the historical process $\widehat{X}_t$, with jump rate $\widehat{V}_t$ on $\widehat{S}$.
At each jump time an ancestral line jumps onto another uniformly chosen ancestral line in the pool. Between jumps
the ancestral lines evolves as independent copies of the historical process $\widehat{X}_t$.

Let $\Ia$ be an uniform random variable on the index set $[N]:=\{1,\ldots,N\}$, independent of particle model discussed above.

Given $\XX_t$ the distribution of a randomly chosen ancestral line $\XX^{\Ia}_t$ is given by the occupation measure  of the genealogical tree. In addition, in some weak sense we have the convergence
 \begin{equation}\label{def-m-X}
m(\XX_t):=\frac{1}{N}\sum_{i\in [N]}\delta_{\XX^i_t}~\longrightarrow_{N\rightarrow\infty}~\QQ_t\quad \mbox{\rm and}\quad m(\xi_t):=\frac{1}{N}\sum_{i\in [N]}\delta_{\xi^i_t}~\longrightarrow_{N\rightarrow\infty}~\eta_t
\end{equation}
For a more detailed discussion and precise estimates, we refer to~\cite{dm-2000,dm-2007,dm-sch,mathias}, as well as to theorem~\ref{cor-intro} and section~\ref{sec-non-asympt} in the present article.

  The historical process $\widehat{\xi}_t=\left(\xi_s\right)_{s\leq t}\in D_t(S)^N$ encapsulates the ancestral lines at any time horizon, including ancestors at any time level $s$ with no currently living descendants at a time horizon $t\geq s$. 
  The process $\widehat{\xi}_t$ has the same jump rate as $\xi_t$. When a jump occurs, we keep and extend all
 historical  trajectories from the randomly selected particles at that time. 
 
 For a more detailed description of the processes $(\xi_t,\widehat{\xi}_t,\XX_t)$, we refer to section~\ref{hist-gen-sec} dedicated to historical and genealogical tree evolutions, see for instance figure~\ref{figure-xi} and the complete genealogical tree presented in figure~\ref{fig-coalescents}.

Given $\Ia$, the dual process $\widehat{\goodchi}_t\in D_t(S)^N$ is defined as $\widehat{\xi}_t$
but the $\Ia$-th ancestral line is a frozen path with the same law as the historical process $\widehat{X}_t$. More precisely, the jump 
of the $\Ia$-th particle onto the $j$-th one  is replaced
by  the jump of
the $j$-th particle onto the $\Ia$-th frozen one. In addition, at these jumps times we flip the $j$-th and $\Ia$-th row historical path up to the present time horizon.
 A schematic picture of this jump and historical permutation at time $t$ is given below
   \begin{figure}[H]
\centering
\begin{subfigure}[c]{0.15\textwidth}
\centering    
\resizebox{\linewidth}{!}{   
\xymatrix@C=4em@R=1.5em{
\bullet&\ar@<-.25ex>@{-}[l]_{\Ia}\bullet\\
\bullet&\ar@<-.25ex>@{-}[l]_{j}\bullet 
}}
 \end{subfigure} $=\widehat{\goodchi}_{t-}$~~~~$\Longrightarrow$ ~~~$   \widehat{\goodchi}_t=$
\begin{subfigure}[c]{0.15\textwidth}
\centering
\resizebox{\linewidth}{!}{   
\xymatrix@C=4em@R=1.5em{
\bullet\hskip-.15cm&\hskip-.14cm\ar@<-.1ex>@{-}[l]_{j}\circ\ar@<0.1ex>@{.>}[d] \bullet\\
\bullet\hskip-.15cm&\hskip-.14cm\ar@<-.1ex>@{-}[l]_{\Ia}~\bullet
}
}
 \end{subfigure}
\end{figure}
The duality terminology comes from the fact that $\widehat{\goodchi}_t$ coincides with the conditional stochastic process $\widehat{\xi}_t$ given a frozen ancestral line, under some many-body Feynman-Kac measure (cf. theorem~\ref{theo-duality} and the Bayes' formula (\ref{intro-equation-duality-ref})). We also underline that   $ \widehat{\goodchi}_{t}:=( \widehat{\goodchi}_{t}(s))_{s\leq t}\in D_t(S)^N$ is not necessarily the historical path of an auxiliary Markov process.  

For a more detailed description of $\widehat{\goodchi}_t$ in terms of generators we refer to section~\ref{sec-particle-frozen-line} dedicated to particle models with a frozen ancestral line, see for instance (\ref{ref-wG-m-k}), the evolution diagrams given in figure~\ref{figure-xi-3-goodchi} as well as  the graphical description of the process $\widehat{\goodchi}_t$ provided in figure~\ref{figure-good-chi}.

The genealogical type tree $\YY_t=(\YY^i_t)_{1\leq i\leq N}\in D_t(S)^N$ of the dual process discussed above is also 
defined as $\XX_t$. The main difference is that the $\Ia$-th ancestral line is frozen and $\YY^{\Ia}_t\stackrel{law}{=}\widehat{X}_t$. 

 The remaining $(N-1)$ path-valued  particles
$\YY^-_t:=(\YY^j_t)_{j\not=\Ia}$ are defined as above with a rescaled jump rate $(1-1/N)\widehat{V}_t$, and including an auxiliary
jump rate $2\widehat{V}_t/N$ at which the path-particle jump onto the first frozen ancestral line.

For a more detailed description of these dual ancestral lines, we also refer to section~\ref{sec-particle-frozen-line}, see for instance the graphical description of the dual genealogical tree  provided in figure~\ref{fig-coalescents-goodchi} as well as the generator of the process in path space defined in (\ref{ref-GG-k-replacement}).\\

A realization of the genealogical tree associated with $N=3$ particles with $2$ interacting jumps and the first frozen ancestral line  is illustrated below in figure~\ref{fig-intro-YY}
.
 \begin{figure}[H]
\begin{center}
\vskip-.3cm\hskip.2cm
\xymatrix@C=4em@R=1.5em{
\circ&\ar@<-.25ex>@{-}[l]\circ&\ar@<-.25ex>@{-}[l]\circ\ar@<-.25ex>@{->}[d]&\\
\circ&\ar@<-.25ex>@{-}[l]\circ\ar@<-.25ex>@{->}[d]&\circ\ar@<-.2ex>@{-}[ld]\ar@{-}[rd]^{\YY^3_t}&\circ\ar@<-.25ex>@{-}[l]_{\YY^2_t}\\
&\circ\ar@<.2ex>@{.}[ld]&&\circ\\
\circ&&\circ\ar@{.}[r]^{\YY^1_t}\ar@<.25ex>@{.}[lu]&\circ   \\
&\ar@{-}[l]&\ar@{-}[l]_{\mbox{\footnotesize time axis $[0,t]$}}\ar@{->}[r]     &                   }
\end{center}
\caption{A genealogical tree associated with $N=3$ particles with $2$ interacting jumps. The couple of arrows stands for the interacting jumps, the dotted line represents the frozen ancestral line with $\Ia=1$. The vertical axis stands for the state space $S$.}\label{fig-intro-YY}
\end{figure}

The first main result of the article is the following duality formula.

\begin{theo}[Duality formula]\label{theo-duality}
For any time horizon $t\geq 0$,  any $N\geq 2$ and any bounded measurable function $F$ on $D_t(S)^{(2N)}:=(D_t(S)^N\times D_t(S)^{N})$ we have the almost sure formula
\begin{equation}\label{intro-equation-duality}
\EE\left( F(\widehat{\xi}_t,\XX_t)~\exp{\left[-\int_0^t~m(\XX_s)(\widehat{V}_s)ds\right]}~\right)
=
\EE\left( F(\widehat{\goodchi}_t,\YY_t)~ \exp{\left[-\int_0^t\widehat{V}_s(\YY^{\Ia}_s)ds\right]}~|~\Ia\right)
\end{equation}
\end{theo}
The proof of the above theorem is provided in section~\ref{duality-sec}, see for instance theorem~\ref{theo-duality-in} and corollary~\ref{cor-duality}. We obtain the normalizing constant $\ZZ_t$ by choosing the unit function $F=1$ in (\ref{intro-equation-duality}). Observe that for any $F$ on $D_t(S)$ we have
\begin{equation}\label{intro-equation-duality-unbiased}
(\ref{intro-equation-duality})\Longrightarrow
\EE\left( F(\XX^{\Ia}_t)~\exp{\left[-\int_0^t~m(\XX_s)(\widehat{V}_s)ds\right]}~\right)=\ZZ_t\times\QQ_t(F)
\end{equation}
The last assertion yields the rather well known 
unbiasedness property of the occupation measure of the ancestral lines; see for instance~\cite{dm-2000} and references therein.
Formula (\ref{intro-equation-duality}) can also be interpreted as a Bayes formula for many body probability measures on the product space $D_t(S)^{(2N)}$. Therefore, using the conditional distributions we can  design a Gibbs sampler. Nevertheless to target $\QQ_t$, we do not need to store the complete ancestral tree.  
A more judicious target measure is the probability measure $\Pi_t$  on $D_t(S)^{N+1}$ defined by
$$
\Pi_t(F)=\frac{1}{\ZZ_t}~
\EE\left( F(\XX_t,\XX^{\Ia}_t)~\exp{\left[-\int_0^t~m(\XX_s)(\widehat{V}_s)ds\right]}~\right)
$$
Let $\pi_t$ be the marginal of $\Pi_t$ w.r.t. the variable $\XX_t$. Theorem~\ref{theo-duality}  yields the Bayes  formula
\begin{equation}\label{intro-equation-duality-ref}
  \Pi_t(d(z_1,z_2))=\pi_t(dz_1)~ \AC_t(z_1,dz_2)=\QQ_t(dz_2)~ \BB_t(z_2,dz_1)
\end{equation}
with the Markov  transitions $ \AC_t$ and $\BB_t$ defined by
 $$
 \AC_t(z_1,dz_2):=\PP(\XX^{\Ia}_t\in dz_2~|~\XX_t=z_1)=:m(z_1)(dz_2)
  \quad\mbox{\rm and}\quad
 \BB_t(z_2,dz_1):=\PP(\YY_t\in dz_1~|~\YY_t^{\Ia}=z_2)  
 $$ 
The detailed proof of the above assertion is provided in section~\ref{pg-section}, on page~\pageref{intro-equation-duality-ref-proof}. Equivalently, under the probability measure $\Pi_t$, any randomly chosen ancestral line is distributed with $\QQ_t$. In addition, under $\Pi_t$ the conditional distribution of the genealogical tree $\XX_t$ given a selected ancestral line coincides with the one of the genealogical tree $\YY_t$ given the frozen selected ancestral line.

Next we provide an analytic description of the Bayes formula (\ref{intro-equation-duality-ref}). Integrating (\ref{intro-equation-duality-ref}) w.r.t. the first coordinate  for any $z_1\in D_t(S)^N$ we have
$$
 \AC_t(z_1,dz_2)\ll (\pi_t\AC_t)(dz_2):=\int \pi_t(dy_1)~\AA_t(y_1,dz_2)=\QQ_t(dz_2)
$$
Thus, we can define the dual operator $\AC^{\star}_{t,\pi_t}$ from $\LL_1(\pi_t)$ into $\LL_1(\QQ_t)$ given for any
$F\in \LL_1(\pi_t)$ by the Radon-Nikodym formula
$$
\AC^{\star}_{t,\pi_t}(F):=\frac{d\pi_t^F\AC_t}{d\QQ_t}\quad \mbox{\rm with}\quad \pi_t^F(dy):=\pi_t(dy)~F(y)
$$
For a more thorough discussion on dual Markov transitions, we refer the reader to~\cite{dm-ledoux-miclo,revuz}.
In addition, for any conjugate integers $1/p+1/q=1$ with $p,q\geq 1$ 
we have
$$
(F,G)\in\left(\LL_p(\pi_t)\times\LL_q(\QQ_t)\right)\qquad \QQ_t\left(\AC^{\star}_{t,\pi_t}(F)~G\right)=\pi_t\left(F~\AC_t(G)\right)
$$
We define in the same way the dual operator $\BB_{t,\QQ_t}^{\star}$ from  $\LL_1(\QQ_t)$ into $\LL_1(\pi_t)$ by
$$
\BB^{\star}_{t,\QQ_t}(G):=\frac{d\QQ_t^G\,\BB_t}{d\pi_t}\quad \mbox{\rm with}\quad \QQ_t^G(dy):=\QQ_t(dy)~F(y)
\Longrightarrow
\pi_t\left(F~\BB^{\star}_{t,\QQ_t}(G)\right)=\QQ_t\left(\BB_t(F)~G\right)
$$

Let $P_{s,t}$ be the Markov semigroup of the reference process $X_t$, that is
$$
P_{s,t}(x,dy):=\PP\left(X_t\in dy~|~X_s=x\right)
$$

We consider the following regularity  condition
\begin{equation}\label{H0}
(H_0)\quad~\exists h>0\quad\mbox{\rm s.t.}\quad \forall t\geq 0\quad
\forall x\in S\quad
\rho(h)~\mu_{t,h}(dy)~\leq P_{t,t+h}(x,dy)~\leq \rho(h)^{-1}~\mu_{t,h}(dy)
\end{equation}
 for some constant $\rho(h)>0$  and some collection of probability measures $\mu_{t,h}$ on $S$  indexed by $h>0$ and $t\geq 0$ whose values do not depend on the parameters $(x,y)$.  A discussion on the above condition is provided at the end of this section.
 
Let $\mbox{\rm osc}(F)$ and
$\Vert \mu_1-\mu_2 \Vert_{\tiny\rm tv}$  be the oscillation of a function $F$
on $D_t(S)$  and the total variation distance between probability measures $\mu_1$ and $\mu_2$ (cf. section~\ref{ref-basic-notation} for a 
more precise definition).

 In this notation, the second main result of the article can be stated basically as follows.
\begin{theo}[Particle Gibbs-Glauber dynamics]\label{theo-PG-intro}
For any time horizon $t\geq 0$ the measure $\QQ_t$ is reversible w.r.t.  the Markov transition $\KK_t:=\BB_t\AC_t$ on $D_t(S)$ defined for any bounded measurable function $F$
on $D_t(S)$ and any path $x\in D_t(S)$ by the formula
$$
\KK_t(F)(x):=\EE\left(m(\YY_t)(F)~|~\YY_t^{\Ia}=x\right)
$$
In addition, when $(H_0)$ is satisfied, for $n\geq 1$, any probability measure $\mu$ on  $D_t(S)$  and any bounded function $F$ on $D_t(S)$ we have the contraction inequality
\begin{equation}\label{estimates-intro}
 \Vert\mu \KK_t^n-\QQ_t\Vert_{\tiny\rm tv}\leq(c~(t\vee 1)/N)^n~\left\Vert \mu -\QQ_t \right\Vert_{\tiny\rm tv}
\end{equation}
with the $n$-th iterate $ \KK_t^n$ of the Markov transition $\KK_t$ defined sequentially  by
the integral formula 
$$
 \KK_t^n(x_1,dx_3):=\int~\KK_t^{n-1}(x_1,dx_2)~\KK_t(x_2,dx_3)
$$

In the above display,  $c$ stands for  some explicit finite constant whose value doesn't depend on the parameters $(F,t,n,N)$. 
\end{theo}

The proof of the above theorem is provided in section~\ref{pg-section}.

For any given time horizon $t\geq 0$, the integral operator
$\KK_t$ is the probability transition of a  discrete generation 
Markov chain $\XX^{(n)}_t$  taking values in the path space $D_t(S)$ and indexed by the integer parameter $n\in\NN$; that is we have
$$
\PP\left(\XX^{(n+1)}_t\in dz~|~\XX^{(n)}_t=x\right)=\KK_t(x,dz):=\int~\BB_t(x,dy)~\AC_t(y,dz)
$$
Initially, we can sample a genealogical tree $\XX_t=(\XX^i_t)_{i\in [N]}$, then we pick randomly an ancestral line $\XX_t^{\Ia}$ and set $\XX^{(0)}_t=\XX_t^{\Ia}$.

For any given $x\in D_t(S)$ and $y\in D_t(S)^N$, we summarize  the overlapping transition of the particle Gibbs sampler graphically as follows: 
\begin{linenomath*}
 $$
\left\{
\begin{array}{rcl}
\YY^{(n)}_t&=&y\\
\XX^{(n)}_t&=&x
\end{array}
\right\}\stackrel{\BB_t}{\longrightarrow}
\left\{
\begin{array}{rcl}
\YY^{(n+1)}_t&=&\overline{y}\sim \left(\YY_t~|~\YY^{\Ia}_t=x\right)\\
\XX^{(n)}_t&=&x
\end{array}
\right\} \stackrel{\AC_t}{\longrightarrow} \left\{
\begin{array}{rcl}
\YY^{(n+1)}_t&=&\overline{y}\\
\XX^{(n+1)}_t&=&\overline{x}\sim m\left(\overline{y}\right)
\end{array}
\right\}.
$$
\end{linenomath*}

A realization of the overlapping transition $\XX^{(n)}_t\leadsto \XX^{(n+1)}_t$ 
 for a genealogical tree with $N=3$ ancestral lines  is illustrated by the following schematic  diagram:
 \begin{figure}[H]
\begin{center}
\vskip-.3cm\hskip.2cm
\xymatrix@C=4em@R=1.5em{
\circ&&\circ\ar@<.2ex>@{-}[ld]\ar@<-.5ex>@{->}[ld]\ar@{-}[rd]&\circ\ar@<-.25ex>@{-}[l]\ar@<-1ex>@{->}[l]_{\XX^{(n+1)}_t}\\
\circ&\circ\ar@<.2ex>@{.}[ld]\ar@<-.6ex>@{->}[ld] &&\circ\\
\circ&&\circ\ar@{.}[r]^{\XX^{(n)}_t}\ar@<.25ex>@{.}[lu]&\circ \\
&\ar@{-}[l]&\ar@{-}[l]_{\mbox{\footnotesize time axis $[0,t]$}}\ar@{->}[r]     &                                    }
\end{center}
\caption{A realization of the transition $\XX^{(n)}_t\leadsto \XX^{(n+1)}_t$ of a particle Gibbs sampler on an genealogical tree
 with $N=3$ ancestral lines. The dotted and plain lines account together for  the three paths in the genealogical tree $\YY^{(n+1)}_t$; the frozen ancestral line is represented by the dotted line  $\XX^{(n)}_t$ selected at rank $(n)$. The sequence of arrows stands for the selected ancestral line $\XX^{(n+1)}_t$ at rank $(n+1)$. The vertical axis stands for the state space $S$.}
 \end{figure}
As shown in (\ref{intro-equation-duality-unbiased}), under the many-body measure $\Pi_t$ a randomly selected ancestral line $\XX^{\Ia}_t$ is unbiased but this random path is biased w.r.t. the distribution of the interacting particle system. Propagation of chaos type estimates allow to quantify this bias, see~\cite{dm-2000,dm-sch}, as well as chapter 15 in the monograph~\cite{d-2013} for discrete time generation particle systems.  

In this connexion, we underline that  the Particle Gibbs-Glauber dynamics presented in theorem~\ref{theo-PG-intro} allows to improve the precision of 
the conventional particle interpretations of Feynman-Kac measures by sampling sequentially a series of particle Gibbs samplers on path spaces with frozen trajectories. Each iteration of the sampler 
reduces the distance between the distribution of the ancestral path and the desired target measure. 
For instance, whenever $(H_0)$ is met and $\XX^{(0)}_t=\XX_t^{\Ia}$, the propagation of chaos property  presented in (\ref{ref-bias-path-space}) and (\ref{estimates-intro}) yield the following theorem

\begin{theo}\label{cor-intro}
For any $n\geq 1$ we have the inequalities
$$
\left\Vert  \mbox{\rm Law}\left(\XX^{\Ia}_t\right)-\QQ_t\right\Vert_{\tiny\rm tv}\leq c~t/N
\quad\mbox{\rm and}\quad
\left\Vert \mbox{\rm Law}\left(\XX^{(n)}_t\right)-\QQ_t \right\Vert_{\tiny\rm tv}\leq (c~t/N)^n\times\left\Vert  \mbox{\rm Law}\left(\XX^{\Ia}_t\right)-\QQ_t\right\Vert_{\tiny\rm tv}
$$
\end{theo}

 We end this section with some comments on our regularity assumptions. 
 
 Condition (\ref{H0}) is satisfied for jump-type
elliptic diffusions on compact manifolds $S$ with a bounded jump rate, see for instance the pioneering work of Aronson~\cite{aronson}, Nash~\cite{nash} and Varopoulos~\cite{varopoulos} on Gaussian estimates for heat kernels on manifolds. The estimates (\ref{estimates-intro}) are also met under weaker regularity assumptions such as 
conditions $(H_1)$ and $(H_2)$ stated in (\ref{ref-H1}) and (\ref{ref-H2}). Nevertheless these conditions depends on the stability properties of the semigroup of the unknown Feynman-Kac measures $\eta_t$;  thus this type of condition is difficult to check in practice. 

Also recall that Feynman-Kac semigroups for time homogeneous models equipped with a 
reversible reference process $X_t$ can also be turned into conventional Markov semigroups of $h$-processes, see for instance~\cite{dm-2000,dm-sch,mathias} as well as chapter 27 in~\cite{dm-penev} and references therein. In this context, the long time behavior of Feynman-Kac semigroups can be discussed in terms of the spectral properties of the $h$-process. For quadratic potential functions and Ornstein-Uhlenbeck reference processes the Feynman-Kac model discussed above reduces to the harmonic oscillator. In some situations, the potential function can be chosen so that the exponential weight in  (\ref{def-FK-intro}) is an exponential change of measure, see for instance section 4.2 in~\cite{dm-2000}, as well as chapter 18 in~\cite{dm-penev}. In this context the corresponding Feynman-Kac semigroup also coincides with the semigroup of a conventional Markov process. 
In more general cases, the spectral properties of the $h$-processes are unknown.
 
\subsection{Illustrations and comments}
This section gives some comments on the impact of the above results on some application domain areas. We also provide a detailed discussion on some numerical aspects of the particle Gibbs-Glauber dynamics introduced above as well as some comparisons with existing literature on interacting particle systems. 

\subsubsection{Some application domains}\label{sec-appli-FK}

 As mentioned in the introduction, 
the Feynman-Kac measures (\ref{def-FK-intro}) and their mean field particle interpretations appear in wide variety of applications including in biology, physics,  as well as in signal processing and mathematical finance.  \\
\indent  Continuous time models arise when the process $X_t$
 is derived from physical or natural evolution principles, such as continuous time signals in target tracking filtering problems~\cite{singer},
stochastic population dynamics describing species competition and populations growths~\cite{lande},
Langevin gradient-type diffusions including their overdampted versions describing the evolution of a particle in a fluid~\cite{langevin}, as well as Brownian fluctuations of atomic structures in molecular chemistry~\cite{kubo}, 
and many others. \\
\indent  The potential function $V_t$ depends on the problem at hand. In nonlinear filtering, it represents the log-likelihood of the robust optimal filter. In population dynamics, $V_t$ can be interpreted as a killing rate of a branching process. In statistical physics and quantum mechanics, it represents the ground state energy (a.k.a. local energy) of a physical system, including molecular and atomic systems. It is clearly out of the scope of the present article to enter into the details of all of these models. For a more thorough discussion on these application domain areas, we refer to the books~\cite{d-2004,d-2013,dm-2000,dm-penev} and the reference therein. \\
\indent In most cases we are mainly interested in computing the final-time marginal of the Feynman-Kac measures (\ref{def-FK-intro}).
For instance, in nonlinear filtering these measures represent the robust optimal filter, while the path space measures represents the full conditional distributions of the random trajectories of the signal w.r.t. the observation process. Thus, they also solve the smoothing problem by estimating the signal states at any given time using observations from larger time intervals.

Apart from few notable exceptions such as for linear-Gaussian models 
in Kalman-Bucy filtering theory and for the harmonic oscillator in  the spectral theory of Schr\"odinger operators, the flow of final-time marginal measures has no finite recursion and cannot be solved analytically.

 In signal processing literature, the interacting particle system $\xi_t$ discussed above is also known as a particle filter on path space.
 
  In Quantum Monte Carlo literature, the particle system $\xi_t$ discussed above is also known as the population Monte Carlo algorithm and the particles $\xi^i_t$ are often referred as walkers or replica. 
  
  This class of processes can be interpreted as Moran type interacting particle systems~\cite{moran-1,moran-2}.
They can also be seen as Nanbu type interpretation of a particular spatially homogeneous generalized Boltzmann equation~\cite{dm-2000-moran,mel}.

\subsubsection{Practical and numerical aspects}
 In some particular instances, the random paths of the process $X_t$ can be sampled exactly on any time discretization
mesh. This class of models includes linear-Gaussian and geometric-type Brownian models, as well as  some piecewise deterministic processes and some classes of one-dimensional jump-diffusion processes~\cite{beskos,beskos-2,beskos-3,beskos-4,casella}. Discretization-free simulation procedures for general diffusion processes based on sequential importance sampling techniques have also been developed in~\cite{fearnhead}.
Using these discretization-free simulation procedures the interacting jump particle systems discussed in this article, including the particle 
Gibbs-Glauber dynamics can be sampled perfectly using conventional Poisson thinning techniques
(a.k.a. Gillespie's algorithm~\cite{gillespie}). The resulting particle sampler provides an estimate of the marginal of the Feynman-Kac measures (\ref{def-FK-intro}) on the random paths w.r.t. any time discretization mesh. 

 More generally, the simulation of the random trajectories of $X_t$ requires to discretize the time parameter. For a more thorough discussion on the time discretization of stochastic processes we refer to the seminal book by Kloeden and Platen~\cite{kloeden}.

  This additional level of approximation
may also corrupt some  statistical properties of the continuous time process. 
For instance, the reversible properties of overdampted Langevin diffusions are lost for any Euler-Maryuama discretization of the underlying diffusion. 
In this context, a Metropolis-Hastings type adjustment (a.k.a. MALA) is required to recover the reversibility property w.r.t. some prescribed target invariant measure~\cite{roberts}. From the physical viewpoint, the random paths simulated by MALA algorithms are based {\em on auxiliary non physical rejection-type transitions} so that they loose their initial physical interpretation. Therefore, in physics and statistics, the unajusted Langevin algorithm (a.k.a. ULA) is often preferred  to describe the "true" random trajectories of the system. Under appropriate global Lipschitz conditions on the gradient of the confinement potential function several bias-type estimates can be found in~\cite{dalayan,durmus}.\\
\indent  In the same vein, the sampling of the particle Gibbs-Glauber dynamics described in theorem~\ref{theo-PG-intro} requires some  Euler-type discretization as soon as the underlying process $X_t$ cannot be directly sampled. 
In this situation, one natural strategy is to consider the discrete time version of the Feynman-Kac measures $\eta_t$ defined as in (\ref{def-FK-intro}) by replacing $X_t$ by some discrete time approximation (see for instance chapter 5 in~\cite{d-2013} and the references therein). In this context, several 
discrete time approximations of the particle Gibbs-Glauber dynamics discussed above can be designed using the discrete time particle Gibbs samplers discussed in~\cite{adh-2010,DMKP:16}. In contrast with MALA algorithms the reversible-type properties of the resulting 
Gibbs samplers in discrete time are preserved w.r.t. to the discrete-time version of the target Feynman-Kac measures. In addition, these discrete time approximations are not based on any type of auxiliary Metropolis-Hasting rejection so that they preserve their physical interpretations. \\
\indent Several bias-type estimates between continuous and discrete time Feynman-Kac measures can be found in~\cite{d-2013,dm-jacob-lee-murrau-peters,dm-jacod-02}. Most of these estimates are concerned with the time discretization of the terminal-time marginal of the Feynman-Kac measures (\ref{def-FK-intro}), including uniform estimates w.r.t. the time horizon.
The extension of these results to 
path space models remains an important open research question.

\subsubsection{Comparisons with diffusion type particle models}
 The interacting particle systems discussed in the present article differ from nonlinear and interacting diffusion processes arising in fluid mechanics and granular flows~\cite{bene-1,bene-2,mckean-1,mckean-2,tamura,tamura-2}. In this context, the interaction mechanism is encapsulated in the drift of diffusion-type particles.
One common feature of these interacting processes is the nonlinearity of the distribution flow associated with these stochastic processes.\\
\indent  One natural idea is to interpret the mean field particle systems associated with these processes as a stochastic perturbation of a nonlinear process. This interpretation allows to enter the stability properties of the  nonlinear process into the convergence analysis of these particle algorithms. This  technique has been developed in~\cite{dm-guionnet,dm-guionnet-2,dm-2000} for discrete time Feynman-Kac models and further extended in~\cite{mathias}  to continuous time models.  The extension to stochastic diffusion flows and McKean-Vlasov type nonlinear diffusions are developed in~\cite{adm-19,dm-ssd-19}.

Theorem~\ref{theo-key-decom} in the present article also provides a novel backward stochastic perturbation formula which simplifies the stability analysis of these models and provides sharp propagation of chaos estimates.\\
\indent We underline that the stochastic perturbation techniques discussed above and in the present article differs from the  log-Sobolev functional techniques~\cite{malrieu,malrieu-2}, entropy dissipation approaches~\cite{cmcv-2006,cordero},
as well as gradient flows in Wasserstein metric spaces, optimal transportation inequalities~\cite{gbg-13,cmcv-2006,cattiaux,otto,otto-2} and the more recent variational approach~\cite{adm-18} currently used in the analysis of
gradient type flow interacting diffusions. \\
\indent In this connection, we mention that the backward perturbation 
analysis developed in the present article relies on weak Taylor expansions of the evolution semigroup of Feynman-Kac measures. We  project to extend these expansions to nonlinear diffusions in a forthcoming article. \\
\indent The duality formula and the particle Gibbs-Glauber dynamics introduced in this article
open up a whole new avenue of research questions.\\
\indent Recall that the Feynman-Kac measures (\ref{def-FK-intro}) can be interpreted as the distribution of the random paths of a non absorbed 
particle evolving as $X_t$ and killed at rate $V_t$. This class of models are often referred as particle models in absorbing medium with soft obstacles~\cite{dd-2004,dm-2000,dm-sch}. A natural research project is to extend this framework to absorbing medium with hard obstacles~\cite{dm-villemonais,delyon-cerou-18,denis-14}.\\
\indent Another important question is to extend the Taylor expansions of the Gibbs sampler 
developed in ~\cite{DMKP:16} to continuous time models. One possible route is to combine the weak Taylor expansions developed in~\cite{dm-patras-ruben} for particle approximating measures with the backward analysis developed in the present article. \\
\indent We mention that  the perturbation analysis developed in~\cite{DMKP:16} allows to destimate the $\LL_p$-decays rates to equilibrium in terms of the norm of integral operators. In this connection, one important question is to quantify with more precision the exponential convergence rates to equilibrium of the Particle Gibbs-Glauber dynamics stated in theorem~\ref{theo-PG-intro}.

\section{Some preliminary results}
\subsection{Basic notation}\label{ref-basic-notation}
Let $\Ba(E)$ be the Banach space of bounded functions $f$ on some measurable space $(E,\Ea)$ equipped with the uniform norm $\Vert f\Vert:=\sup_{x\in E}\vert f(x)\vert$.
Also let $ \mbox{\rm Osc}(E)\subset \Ba(E)$ be the subset of functions $f$ with at most. unit oscillations; that is  s.t. $\mbox{\rm osc}(f):=\sup_{x,y}\vert f(x)-f(y)\vert\leq 1$.  

We also let $\Ma(E)$ be the set of finite signed measures on $E$, $\Ma_+(E)\subset \Ma(E)$ the subset of positive measures and $\Pa(E)\subset \Ma_+(E)$ the subset of probability measures. 
Given a measure $\mu$ on $E$ we write $\mu(f)$ the Lebesgue integral given by
$$
\mu(f)=\int~\mu(dx)~f(x)
$$
The total variation norm on the set $\Ma(E)$ is defined by 
\begin{eqnarray}
\Vert \mu\Vert_{\tiny\rm tv}&:=&\sup{\left\{ \vert \mu(f)\vert~:~f\in \mbox{\rm Osc}(E)\right\}}=\frac{1}{2}~\sup{\left\{ \vert \mu(f)\vert~:~f\in \Ba(E),~\Vert f\Vert\leq 1\right\}}
\end{eqnarray}
We denote by $ [s,t]_n$ the collection Weyl chambers defined for any $n\geq 1$ by
\begin{equation}
[s,t]_n:=\left\{(r_1,\ldots,r_n)\in [s,t]^n~:~s\leq r_1\leq\ldots\leq r_n\leq t\right\}\label{weyl-chamber-def}
\end{equation}
We denote by $dr=dr_1\times\ldots\times dr_n$ the Lebesgue measure on $[s,t]_n$.

For a given $N\geq 1$, we let $\langle N\rangle$ be the semigroup of mappings from $[N]:=\{1,\ldots,N\}$ into itself equipped with the composition of mappings and the neutral element $e(i)=i$.  Also let  $\left[N\right]^2_0\subset [N]^2$ be the subset of indices
       $$
\left[N\right]^2_0:=\{(i,j)\in [N]^2~:~i\in [N]\quad \&\quad j\in [N]-\{i\}\}
     $$

\subsection{Integral operators}
We introduce some integral operator notation needed from the onset.
For any bounded positive integral operator  $Q(x,dy)$ and any
$(\mu,f,x)\in (\Ma(E)\times\Ba(E)\times E)$ we define by $\mu Q\in \Ma(E)$ and $Q(f)\in \Ba(E)$ by the formulae
$$
(\mu Q)(dy):=\int\mu(dx)Q(x,dy)\quad\mbox{\rm and}\quad Q(f)(x):=\int~Q(x,dy)~f(y)
$$ 
 By Fubini theorem we have $\mu Q f:=\mu(Q(f))=(\mu Q)(f)$. Given a pair of operators $Q_1$ and
 $Q_2$ we denote by $Q_1Q_2:=Q_1\circ Q_2$ the composition of the operators defined for functions $f$ on $S$ by
 \begin{equation}\label{ref-composition}
(Q_1\circ Q_2)(f):=( Q_1Q_2)(f)=Q_1(Q_2(f))
 \end{equation}
 We also write $Q^n$ the $n$ iterate of $Q$ defined by
 the recursion $Q^n(f)=Q(Q^{n-1}(f))=Q^{n-1}(Q(f))$.
 
When $Q(1)>0$ we let $\overline{Q}$
be the Markov operator
$$
\overline{Q}~:~f\in \Ba(E)\mapsto \overline{Q}(f):=Q(f)/Q(1)\in \Ba(E)
$$
We also let $\phi$ be the mapping from $\Pa(E)$ into itself defined by
 \begin{equation}\label{ref-phi-Q}
\phi(\eta)=\eta Q^{\eta}\quad \mbox{\rm with}\quad Q^{\eta}:=\frac{Q}{\eta Q(1)}\Longrightarrow
 \eta Q^\eta(1)=1\quad\mbox{\rm and}\quad  \phi(\delta_x)(f)=\overline{Q}(f)(x)
\end{equation}
 Notice that
 $$
 Q^{\eta}(1)=~\mu Q^{\eta}(1)~Q^{\mu}(1)\Longrightarrow  (\mu Q^{\eta}(1))^{-1}= \eta Q^{\mu}(1)
 $$
 
The Dobrushin ergodic coefficient   $\beta_{\mbox{\tiny\rm dob}}(M)$ of a Markov transition
 $M(x,dy)$ from $E$ into itself is defined by
 $$
 \beta_{\mbox{\tiny\rm dob}}(M):=\sup_{x,y\in E}{\Vert M(x,\point)-M(y,\point)\Vert_{\tiny\rm tv}}=\sup{\left\{ \mbox{\rm osc}(M(f))~:~f\in \mbox{\rm Osc}(E)\right\}}
 $$
 For any $\mu_1,\mu_2\in \Pa(E)$ and any pair of Markov transition
 $M_1,M_2$ from $E$ into itself  we have
 \begin{equation}\label{ref-contraction}
  \beta_{\mbox{\tiny\rm dob}}(M_1M_2)\leq  \beta_{\mbox{\tiny\rm dob}}(M_1)~\beta_{\mbox{\tiny\rm dob}}(M_2)\quad \mbox{\rm and}\quad\Vert\mu_1 M-\mu_2 M\Vert_{\tiny\rm tv}\leq \beta_{\mbox{\tiny\rm dob}}(M)~\Vert\mu_1-\mu_2 \Vert_{\tiny\rm tv}
 \end{equation}

 \subsection{Taylor expansions}

  The Feynman-Kac semigroups discussed in section~\ref{ref-def-D} have the same form as the map
  $\phi$ discussed in (\ref{ref-phi-Q}) (see for instance (\ref{forward-FK-Gamma-red})).
  The stochastic perturbation analysis developed in section~\ref{sec-perturbation} is mainly based on a second order  Taylor expansion of these maps (see for instance the proof of proposition~\ref{prop-forward-backward-ref}, as well as the perturbation semigroup equation presented in theorem~\ref{theo-epsilon-delta} and the Aleeksev-Gr\"obner interpolation formula stated in theorem~\ref{theo-key-decom}). 
    
 To describe in some details these Taylor expansions, consider the collection of first order integral operators $ \partial_\eta\phi$ indexed by $\eta\in \Pa(E)$ and defined by
 
\begin{equation}\label{ref-osc-partial-phi-0}
 \partial_\eta\phi~:~f\in \Ba(E)\mapsto  \partial_\eta\phi(f)=Q^{\eta}\left[f-\phi(\eta)(f)\right]\in \Ba(E)\Longrightarrow
 \eta \partial_\eta\phi=0=\partial_\eta\phi(1)
\end{equation}
Rewritten in integral form, we have
$$
\partial_\eta\phi(f)(x)=\int~\partial_\eta\phi(x,dy)~f(y)\quad \mbox{\rm with}\quad
\partial_\eta\phi(x,dy):=Q^{\eta}(x,dy)-Q^{\eta}(1)(x)~\phi(\eta)(dy)
$$
 For any
 $\eta,\nu\in\Pa(E)$ we have the first order Taylor expansion
  \begin{eqnarray*}
\phi(\nu)-
\phi(\eta)&=&
\eta Q^{\nu}(1)~ \times~(\nu-\eta) \partial_\eta\phi \quad \mbox{\rm with}\quad
(\nu-\eta) \partial_\eta\phi(dy):=\int~(\nu-\eta)(dx) \partial_\eta\phi(x,dy)
\end{eqnarray*}

Also observe that
\begin{equation}\label{ref-osc-partial-phi}
\begin{array}{l}
\displaystyle\partial_\eta\phi(f)(x)=Q^{\eta}(1)(x)~\int~\eta(dy)~Q^{\eta}(1)(y)~\left(\overline{Q}(f)(x)-\overline{Q}(f)(y)\right)\\
\\
\Longrightarrow
\Vert \partial_\eta\phi(f)\Vert\leq \Vert Q^{\eta}(1)\Vert~\mbox{\rm osc}(\overline{Q}(f))\quad \mbox{\rm and}\quad
\Vert \phi(\nu)-
\phi(\eta)\Vert_{\tiny\rm tv}\leq   \beta_{\nu,\eta}(\phi)~\Vert \nu-
\eta\Vert_{\tiny\rm tv}
\end{array}
\end{equation}
with
$$
  \beta_{\nu,\eta}(\phi):=\left[\Vert Q^{\nu}(1)\Vert\wedge \Vert Q^{\eta}(1)\Vert\right]~  \beta_{\mbox{\tiny\rm dob}}(\overline{Q})
$$
More generally, using the identity
\begin{equation}\label{key-fraction}
\frac{1}{x}=\sum_{0\leq k<n}~(1-x)^k+\frac{(1-x)^{n}}{x}
\end{equation}
which is valid for any $x>0$ and $n\geq 1$, we check
the Taylor with remainder expansion
  \begin{eqnarray}\label{ref-mean-field-cv}
\phi(\nu)&=&\phi(\eta)+\sum_{1\leq k\leq n}~\frac{1}{k!}~(\nu-\eta)^{\otimes k}~ \partial_\eta^{k}\phi+\frac{1}{(n+1)!}~(\nu-\eta)^{\otimes (n+1)} ~\overline{\partial}^{n+1}_{\nu,\eta}\phi
 \end{eqnarray}
In the above display, $ \partial_\eta^{k}\phi$ stand for the collection of integral operators
 $$
 \partial_\eta^{k}\phi(f):=(-1)^{k-1}~k!~\left[Q^{\eta}(1)^{\otimes (k-1)}\otimes\partial_\eta\phi(f)\right]
\quad\mbox{\rm and}\quad
\overline{\partial}^{n+1}_{\nu,\eta}\phi:=\eta Q^{\nu}(1)~ \partial_\eta^{n+1}\phi
 $$
 
 For any $\mu,\eta\in \Pa(E)$ we have the decomposition
 \begin{eqnarray*}
 \partial_{\eta}\phi(f)=Q^{\eta}[f-\phi(\eta)f]
 &=&\mu Q^{\eta}(1)~ \left( \partial_{\mu}\phi(f)+Q^{\mu}(1)~[\phi(\mu)-\phi(\eta)](f)\right)
 \end{eqnarray*}
 
 \subsection{Carr\'e du champ operators}
 
 The generator $L$ of a stochastic process $X_t$ on some measurable space $(E,\Ea)$ provides a simple and natural way to define the evolution of 
 the random path of the process. 
 
 These operators are also the stepping stone of stochastic calculus. They play the same role as the vector fields associated with a dynamical system and they allow to derive integration by parts formula; see for instance section~\ref{sec-mean-field} in the context of mean-field particle processes.  The carr\'e du champ operator $\Gamma_{L}$ associated with the generator $L$  characterizes the predictable quadratic part of these integration by parts formula (see for instance (\ref{ref-ibp-mean-field}) and (\ref{ref-ibp-mean-field-2})).
 
 In this section we review some basic inequalities that allows to quantify the first order fluctuation term as well as the second order bias term in the  Aleeksev-Gr\"obner interpolation formula stated in theorem~\ref{theo-key-decom}. 
  
 The carr\'e du champ operator $\Gamma_L$ acts on an algebra of functions $\Da(L)\subset \Ba(E)$ and it is defined by the quadratic form
$$
(f,g)\in \Da(L)^2\mapsto
\Gamma_{L}(f,g)=L(fg)-fL(g)-gL(f)\in \Ba(E)
$$ 
When $f=g$ sometimes we write $\Gamma_{L}(f)$ instead of $\Gamma_{L}(f,f)$. We also have the Cauchy-Schwartz inequality
 \begin{equation}\label{ref-Gamma-CS-f-g}
 \vert\Gamma_{L}(f,g)\vert\leq\sqrt{\Gamma_{L}(f,f)\Gamma_{L}(g,g)} \quad\mbox{\rm and}\quad \Gamma_{L}(cf)=c^2~\Gamma_{L}(f)
 \end{equation}
The above inequality yields the estimate
 \begin{equation}\label{ref-Gamma-f-g}
\Gamma_{L}(f+g) =\Gamma_{L}(f)+\Gamma_{L}(g)+2\Gamma_{L}(f,g)\leq \left[\sqrt{\Gamma_{L}(f)}+\sqrt{\Gamma_{L}(g)}\right]^2
 \end{equation}
Let $L^d$ be some bounded jump-type generator of the
following form
$$
L^d(f)(u)=\lambda(u)~\int~(f(v)-f(u))~J(u,dv)
$$
for some bounded rate function $\lambda$ and some Markov transition $J$ on $E$. In this case, we have
$$
\Gamma_{L^d}(f,g)(u)=\int~L^d(u,dv)
 \left(\delta_{v}-\delta_{u}\right)^{\otimes 2}(f\otimes g)
$$
The convergence analysis of the particle measures (\ref{def-m-X}) developed in  section~\ref{sec-particle-flows} is also based
on the $n$-th order operators given by the formula
\begin{equation}\label{n-Gamma}
\Gamma^{(n)}_{L^d}(f_{1},\ldots, f_{n})(u):=\int~
L^d(u,dv)
 \left(\delta_{v}-\delta_{u}\right)^{\otimes n}(f_{1}\otimes\ldots \otimes f_{n})
\end{equation}

Applied to mappings of the form (\ref{ref-phi-Q}), for any $\mu,\eta\in \Pa(E)$ we  have the carr\'e du champ formula
\begin{equation}\label{ref-pr-final}
\begin{array}{l}
 (\eta Q^{\mu}(1))^{2}~\Gamma_L\left(Q^{\eta}(1), \partial_{\eta}\phi(f)\right) \\
 \\
 = \Gamma_L\left(Q^{\mu}(1), \partial_{\mu}\phi(f)\right)+~[\phi(\mu)-\phi(\eta)](f)~\Gamma_L\left(Q^{\mu}(1)\right)
 \end{array}
\end{equation}
for any $f\in \Da(L)$ as soon as $Q^{\eta}(1), \partial_{\eta}\phi(f)\in \Da(L)$.
 
\subsection{Historical processes}

The article discusses  several classes of particle models evolving in path spaces, such as the complete ancestral tree models and the genealogical tree based evolutions discussed in section~\ref{sec-1-ref}.

In this short section, we review some basic facts on the concatenation of paths and the description of historical semigroups.

  For any $s\leq t$, we denote by $D_{s,t}(S)$   the set 
of c\`adl\`ag paths  from $[s,t]$ to the metric space $(S,d_S)$. Fo any given $x=(x_r)_{s\leq r\leq t}\in D_{s,t}(S)$
we let $x_{-}\in D_{s,t}(S)$ the stopped path given by
$$
x_{-}(u):=\left\{
\begin{array}{lcl}
x_u&\mbox{\rm if}& u\in [s,t[\\
x_{t-}&\mbox{\rm if}& u=t
\end{array}
\right.
$$
For any $r\leq s\leq t$ and any $x=(x_{u})_{r\leq u\leq s}\in D_{r,s}(S)$ and $y=(y_{u})_{s\leq u\leq t}\in D_{s,t}(S)$ the concatenate path
$(x\vee y)\in D_{r,t}(S)$ is defined by
$$
(x\vee y)(u):=\left\{
\begin{array}{rcl}
x_u&\mbox{\rm if}& u\in [r,s[\\
y_u&\mbox{\rm if}& u\in [s,t]
\end{array}
\right.\Longrightarrow x\vee y=x_-\vee y
$$

  The set $D_{r,s}(S)$ can be embedded into $D_{r,t}(S)$ by considering the stopped process extension 
  $$
  x=(x_u)_{r\leq u\leq s} \in D_{r,s}(S)~\mapsto   x_{\wedge s}:=(x_{u\wedge s})_{r\leq u\leq t} \in D_{r,s}(S)\quad
  \mbox{\rm with}\quad u\wedge s:=\min{\{u,s\}}
  $$
When $s=t$, the set $D_{s,s}(S)$ reduces to $S$ and for any $x=(x_u)_{r\leq u\leq s}\in  D_{r,s}(S)$ and $y\in S$ the concatenate path
$(x\vee y)\in D_{r,s}(S)$ is given by to the c\`adl\`ag path
\begin{equation}\label{ref-concatenate}
(x\vee y)(u):=\left\{
\begin{array}{rcl}
x_u&\mbox{\rm if}& u\in [r,s[\\
y&\mbox{\rm if}& u=s
\end{array}
\right.\quad \Longrightarrow \quad x=x\vee x_s\quad\mbox{\rm and}\quad
x_-=x\vee x_{s-}
\end{equation}

  Let   $X_{s,t}(y)$ be the stochastic semigroup of the process $X_t$ starting at $X_s=y$ at time $s\leq t$.    The stochastic semigroup $\widehat{X}_{s,t}~:~D_s(S)\mapsto D_t(S)$ of the historical process $ \widehat{X}_t:=(X_s)_{s\leq t}$ is defined for any $x\in D_s(S)$ and $s\leq t$ by the stop-and-go formula
   \begin{equation}\label{def-Xsg-h}
   \widehat{X}_{s,t}(x)=x\vee\widetilde{X}_{s,t}(x)=x_-\vee\widetilde{X}_{s,t}(x)\in D_{t}(S)
      \end{equation}  
   with the mapping
   $$
   \widetilde{X}_{s,t}~: x=(x_u)_{0\leq u\leq s}\in D_s(S)\mapsto
   \widetilde{X}_{s,t}(x)=(X_{s,u}(x_s))_{s\leq u\leq t}\in D_{s,t}(S)
  $$

\subsection{Coalescent operators}

The complete ancestral as well as the genealogical tree  evolutions $(\widehat{\xi}_t,\XX_t)$ discussed in section~\ref{sec-1-ref} belong to the class of interacting 
jump processes in the space of c\`adl\`ag paths. 

When a jump occurs in a genealogical tree evolution, a path-particle is killed and instantaneously  other path-particle duplicates.
When a jump occurs in a complete ancestral tree evolution, a path-particle restarts its evolution from a new selected state. In this situation, the jump of the historical particle is characterized by the concatenation of the stopped process of the historical particle with a new terminal state, from which the particle restarts its free evolution. The dual process $\widehat{\goodchi}_t$ of the historical process incorporates
an additional permutation of the historical paths. 

To describe with some precision these jumps, we need to introduce some new objects.

    We denote by $\mathfrak{T}$ and $\mathfrak{C}$ the set of transpositions $\sigma_{\iota}$ and coalescent maps $c_{\iota}$       indexed by $\iota=(i,j)\in \left[N\right]^2_0$ and
       given for any $k\in [N]$ by
      \begin{equation}\label{def-sigma-c}
       \sigma_{i,j}(k):=\left\{
 \begin{array}{ccl}
 k &\mbox{\rm if}& k\not\in \{i,j\}\\
 j &\mbox{\rm if}& k=i\\
  i &\mbox{\rm if}& k=j
 \end{array}
 \right.\quad \mbox{\rm and}\quad
 c_{i,j}(k):=\left\{
 \begin{array}{ccl}
 k &\mbox{\rm if}& k\not=i\\
 j &\mbox{\rm if}& k=i
 \end{array}
 \right.
      \end{equation}

\begin{defi}
For any $a,b\in \langle N\rangle$  and $x=(x_r)_{s\leq r\leq t}\in D_{s,t}(S)^N$
we set
\begin{equation}\label{def-varsigma-ab}
  x^a:=(x^{a(i)})_{i\in [N]}\in D_{s,t}(S)^N\qquad  \CC_{a,b}(x):=x^a\vee x^b_t\quad \mbox{and}\quad
  \CC_{a}:=\CC_{a,a}
\end{equation}
\end{defi}
Observe that for any $a_1,b_1$ and any $a_2,b_2\in \langle N\rangle$ we have 
     $$
\CC_{a_1,b_1}\circ \CC_{a_2,b_2}=\CC_{a_2\circ a_1,b_2\circ b_1}\quad \mbox{\rm and}\quad
\CC_{a_1}(x)=x^{a_1}
$$
We also mention that for any $(i,j)\in [N]_0^2$ we have
\begin{equation}\label{ref-c-sigma-prop}
\CC_{c_{i,j}}\circ\CC_{\sigma_{j,i}}=\CC_{\sigma_{j,i}\circ c_{i,j}}=\CC_{c_{j,i}}\quad \mbox{\rm and therefore}\quad
\CC_{\sigma_{i,j},c_{j,i}}=\CC_{e,c_{i,j}}\circ \CC_{\sigma_{i,j}}
\end{equation}

Observe that for any index $k\in [N]$ we have the partition
$$
[N]^2_0:=\{(i,j)\in [N]_0^2~:~i\in [N]-\{k\}\}\cup  \{(k,j)\in [N]_0^2~:~j\in [N]-\{k\}\}
$$

\begin{defi}
For any $k\in [N]$ and any $(i,j)\in [N]_0^2$ we let
\begin{equation}\label{def-CC-k-1}
\left(\CC^{k}_{e,c_{i,j}},  \CC^{k}_{c_{i,j}}\right):= \left\{
  \begin{array}{lcl}
 \left( \CC_{e,c_{i,j}}, \CC_{c_{i,j}}\right)&\mbox{\rm if}& i\in [N]-\{k\}\\
  &&\\
  \left( \CC_{\sigma_{k,j},c_{j,k}},  \CC_{c_{j,k}}\right)&\mbox{\rm if}& i=k
  \end{array}
  \right.
  \end{equation}
  \end{defi}
  
 A schematic picture of the jumps $x\leadsto \CC_{e,c_{i,k}}(x)=x\vee x^{c_{i,k}}_t$ and $x\leadsto 
 \CC_{\sigma_{i,k},c_{i,k}}(x)=  x^{\sigma_{i,k}}\vee x^{c_{i,k}}_t$  is given below
   \begin{figure}[H]
\centering
\begin{subfigure}[c]{0.15\textwidth}
\centering    
\resizebox{\linewidth}{!}{   
\xymatrix@C=4em@R=1.5em{
\bullet&\ar@<-.25ex>@{-}[l]_{k}\bullet\\
\bullet&\ar@<-.25ex>@{-}[l]_{i}\bullet 
}}
 \end{subfigure} $=x$~~~~$\Longrightarrow$ ~~$ \CC_{e,c_{i,k}}(x)=$
 \begin{subfigure}[c]{0.15\textwidth}
\centering
\resizebox{\linewidth}{!}{
\xymatrix@C=4em@R=1.5em{
\bullet\hskip-.15cm&\hskip-.14cm\ar@<-.1ex>@{-}[l]_{k}~\bullet\\
\bullet\hskip-.15cm&\hskip-.14cm\ar@<-.1ex>@{-}[l]_{i}\circ\bullet\ar@<-0.1ex>@{.>}[u]}
}
 \end{subfigure} ~~$\mbox{\rm and}$~~$   \CC_{\sigma_{i,k},c_{i,k}}(x)=$
\begin{subfigure}[c]{0.15\textwidth}
\centering
\resizebox{\linewidth}{!}{   
\xymatrix@C=4em@R=1.5em{
\bullet\hskip-.15cm&\hskip-.14cm\ar@<-.1ex>@{-}[l]_{i}\circ\ar@<0.1ex>@{.>}[d] \bullet\\
\bullet\hskip-.15cm&\hskip-.14cm\ar@<-.1ex>@{-}[l]_{k}~\bullet
}
}
 \end{subfigure}
\end{figure}
We underline  that the jump in $ \CC_{e,c_{i,k}}(x)$ occurs on the $k$-th row trajectory while the jump in $ \CC_{\sigma_{i,k},c_{i,k}}(x)$ occurs on the $i$-th row trajectory. 

 A schematic picture of the path-valued jump $x\leadsto \CC_{c_{i,k}}(x)=x^{c_{i,k}}$  is given below
   \begin{figure}[H]
\centering
\begin{subfigure}[c]{0.15\textwidth}
\centering    
\resizebox{\linewidth}{!}{   
\xymatrix@C=4em@R=1.5em{
\bullet&\ar@<-.25ex>@{-}[l]_{k}\bullet\\
\bullet&\ar@<-.25ex>@{-}[l]_{i}\bullet 
}}
 \end{subfigure} $=x$~~~~$\Longrightarrow$ ~~$ \CC_{c_{i,k}}(x)=$
 \begin{subfigure}[c]{0.15\textwidth}
\centering
\resizebox{\linewidth}{!}{
\xymatrix@C=4em@R=1.5em{
\bullet&\ar@<-.1ex>@{-}[l]_{k}\bullet\\
\bullet&\ar@<-.1ex>@{-}[l]_{k}\bullet}
}
 \end{subfigure}
\end{figure}

Between jumps all the particle models  discussed in this article evolve independently as independent copies of the reference process $X_t$ introduced in section~\ref{sec-1-ref}.
Let  $X_{s,t}^i$ with $i\in [N]$ be $N$ independent copies of the stochastic semigroup
$X_{s,t}$. We extend the stochastic semigroups introduced $(\ref{def-Xsg-h})$ to product spaces and we let $\Xa_{s,t}$, $\widehat{\Xa}_{s,t}$ and $\widetilde{\Xa}_{s,t}$ the stochastic semigroups 
with the $i$-th coordinates mappings defined for any for any $z\in S^N$ and any $x\in D_{s}(S)^N$ by the formulae
\begin{equation}\label{def-Xa-wXa}
\begin{array}{l}
\Xa^i_{s,t}(z):=X^i_{s,t}(z^i)\qquad
\widehat{\Xa}^i_{s,t}(x):=\widehat{X}^i_{s,t}(x^i)\quad\mbox{\rm and}\quad
\widetilde{\Xa}^i_{s,t}(x):=\widetilde{X}^i_{s,t}(x^i)\\
\\
\Longrightarrow\quad     \widehat{\Xa}_{s,t}(x)=x\vee\widetilde{\Xa}_{s,t}(x)=x_-\vee\widetilde{\Xa}_{s,t}(x)\in D_{t}(S)^N
\end{array}
\end{equation}

We also denote $\widehat{\xi}_{s,t}$ and $\widehat{\goodchi}_{s,t}$     the stochastic evolution semigroup 
 of the processes $\widehat{\xi}_t$ and  $\widehat{\goodchi}_t$; that is, for any $s\leq t$ we have the stochastic evolution semigroup formulae
 $$
 \widehat{\xi}_{s,t}(\widehat{\xi}_s)=\widehat{\xi}_t \quad \mbox{\rm and}\quad
  \widehat{\goodchi}_{s,t}(\widehat{\goodchi}_s)=\widehat{\goodchi}_t
 \quad \mbox{\rm with the initial condition}\quad
 \widehat{\xi}_{0}= \xi_{0}= \widehat{\goodchi}_{0}
 $$
 We define in the same vein the stochastic evolution semigroups $\XX_{s,t}$ and $\YY_{s,t}$ 
 of the genealogical tree processes $\XX_t$ and  $\YY_t$.

 At jump times, these path-space interacting jumps stochastic flows are described  by the composition of the coalescent operators presented in (\ref{def-varsigma-ab}) and (\ref{def-CC-k-1}) with the stochastic evolution semigroups (\ref{def-Xa-wXa}).

\begin{defi}
For any $a\in \mathfrak{C}$  and any $s\leq t$ and  $k\in [N]$ 
we define
\begin{equation}\label{def-wTa}
\begin{array}{rclcrcl}
\widehat{\Xa}^a_{s,t}&:=&\CC_a\circ 
\widehat{\Xa}_{s,t}&\hskip1cm& \widehat{\Xa}^{k,a}_{s,t}&:=&\CC^k_{a}\circ 
\widehat{\Xa}_{s,t}\\
 &&&&&&\\
 \widehat{\Ta}_{s,t}^{a}&:=&\CC_{e,a}\circ\widehat{\Xa}_{s,t}&&
 \widehat{\Ta}_{s,t}^{\,k,a}&:=&\CC_{e,a}^k\circ\widehat{\Xa}_{s,t}
\end{array}
\end{equation}
\end{defi}
Observe that
$$
\widehat{\Xa}^a_{s,t}(x)=\CC_a(
\widehat{\Xa}_{s,t}(x))=\left(\widehat{\Xa}^{a(i)}_{s,t}(x)\right)_{i\in [N]}=\left(\widehat{X}^{a(i)}_{s,t}\left(x^{a(i)}\right)\right)_{i\in [N]}
$$
This yields for any $r\leq s\leq t$ and any $a,b\in\langle N\rangle$ and  $x\in D_{r}(S)^N$ the composition formula
$$
\left(\widehat{\Xa}^a_{s,t}\circ\widehat{\Xa}^b_{r,s}\right)(x)=\widehat{\Xa}^a_{s,t}\left(\widehat{\Xa}^b_{r,s}(x)\right)=\left(\widehat{X}^{a(i)}_{s,t}\left(\widehat{X}^{(b\circ a)(i)}_{r,s}\left(x^{(b\circ a)(i)}\right)\right)\right)_{i\in [N]}
$$
More generally, for any sequence of mappings  $a_n\in\langle N\rangle$ any $x\in D_{t_0}(S)^N$ and any non decreasing time steps
$t_n$ we have the formula
$$
\begin{array}{l}
\displaystyle\left(\widehat{\Xa}_{t_{n-1},t_n}^{a_n}\circ\ldots \circ \widehat{\Xa}_{t_{0},t_1}^{a_1}\right)(x)\\
\\
\displaystyle=\left(\left(\widehat{X}^{a_{n}(i)}_{t_{n-1},t_n}\circ \widehat{X}^{a_{n-1,n}(i)}_{t_{n-2},t_{n-1}}\circ \ldots \circ\widehat{X}^{a_{1,n}(i)}_{t_{0},t_1}\right)\left(x^{a_{1,n}(i)}\right)\right)_{i\in [N]}\quad\mbox{\rm with}\quad a_{p,n}:=a_{p}\circ a_{p+1}\circ\ldots\circ a_n
\end{array}
$$
Finally observe that for any $a_1,\ldots,a_n\in\mathfrak{C}$ and any $k\in [N]$ and $x\in D_{t_0}(S)^N$ we have
\begin{equation}\label{ref-ancestral-k-proof}
\begin{array}{l}
\displaystyle   \CC^{k}_{c_{i,j}}=1_{i\not=k}~\CC_{c_{i,j}}+1_{i=k}~ \CC_{c_{j,k}}\\
\\
\displaystyle\Longrightarrow\quad\left(\widehat{\Xa}_{t_{n-1},t_n}^{k,a_n}\circ\ldots \circ \widehat{\Xa}_{t_{0},t_1}^{k,a_1}\right)(x)^k=\widehat{X}_{t_0,t_n}^k(x^k)
\end{array}
\end{equation}

The composition of the stochastic flows $\widehat{\Ta}^a_{s,t}$ is slightly more involved.
For instance, for any $s\leq r\leq t$ and any $a,b\in \mathfrak{C}$ and $x\in D_s(S)^N$ we have the concatenate path formulae
$$
\widehat{\Ta}^a_{s,t}(x)= \widehat{\Xa}_{s,t}(x)\vee
\Xa_{s,t}^a(x_s)\in D_t(S)^N
$$
In addition, we have the composition rule
$$
\left( \widehat{\Ta}^a_{r,t}\circ \widehat{\Ta}^b_{s,r}\right)(x)
=\widehat{\Xa}_{s,r}(x)\vee \left(\widetilde{\Xa}_{r,t}\circ \widehat{\Ta}^b_{s,r}\right)(x)\vee \left(\Xa^a_{r,t}\circ \Xa_{s,r}^{b}\right)(x_s)
$$
The proof of the first assertion is immediate. To check the second, observe that
\begin{eqnarray*}
 \widehat{\Ta}^a_{r,t}\left(\widehat{\Ta}^b_{s,r}(x)\right)
 &=&\widehat{\Xa}_{r,t}\left(\widehat{\Xa}_{s,r}(x)\vee
\Xa_{s,r}^b(x_s)\right)\vee
\Xa_{r,t}^a(\Xa^b_{s,r}(x_s))\\
&=&\widehat{\Xa}_{s,r}(x)\vee \widetilde{\Xa}_{r,t}\left(\widehat{\Ta}^b_{s,r}(x)\right)\vee
\Xa_{r,t}^a(\Xa^b_{s,r}(x_s))\quad(\mbox{\rm by (\ref{def-Xa-wXa})})
\end{eqnarray*}

This ends the proof of the	 lemma.
\cqfd

\newpage
A graphical description of these compositions is given below in figure~\ref{fig-coalescents-0}.   \begin{figure}[H]
\centering
$ 
x=$
\begin{subfigure}[c]{0.12\textwidth}
\centering    
\resizebox{\linewidth}{!}{   
\xymatrix@C=4em@R=1.5em{
&\ar@<-.25ex>@{.}[l]_{k}\bullet \\
&\ar@<-.25ex>@{.}[l]_{j_1}~\bullet\\
&\ar@<-.25ex>@{.}[l]_{j_2}~\bullet\\
0 & s\ar@<-.25ex>@{-}[l]
}}
\end{subfigure} 
 $ \quad\leadsto\quad$
\begin{subfigure}[c]{0.5\textwidth}
\centering    
\resizebox{\linewidth}{!}{   
\xymatrix@C=4em@R=1.5em{
\ar@<-.25ex>@{.}[r]^k&\bullet&\ar@<-.25ex>@{-}[l]_{k}\circ\ar@<-.25ex>@{.>}[d]&\bullet\ar@<-.2ex>@{-}[ld]_k&\ar@<-.25ex>@{.}[l]_{k}\\
&&\bullet\ar@<.2ex>@{-}[ld]_{j_1}&&\\
\ar@<-.25ex>@{.}[r]^{j_1}&\bullet&&\bullet\ar@{.}[ru]^{j_1}\ar@<.25ex>@{-}[lu]_{j_1}& \\
\ar@<-.25ex>@{.}[r]^{j_2}&\bullet&\ar@<-.25ex>@{-}[l]_{j_2}\bullet&\ar@<-.25ex>@{-}[l]_{j_2}\circ\ar@<-.25ex>@{.>}[u]&\ar@<.25ex>@{.}[lu]^{j_2} &      \\
0 & s\ar@<-.25ex>@{-}[l] & r\ar@<-.25ex>@{-}[l]& t  \ar@<-.25ex>@{-}[l]&  \ar@<-.25ex>@{.}[l] &
}}
 \end{subfigure} 
 \caption{Genealogical tree associated with the composition $\widehat{\Ta}^{c_{j_2,j_1}}_{r,t}\circ \widehat{\Ta}^{c_{k,j_1}}_{s,r}$}\label{fig-coalescents-0}
\end{figure}
The dotted lines  $(\stackrel{j}{\ldots})$  represent the paths associated with $x$, the 
plain lines $(\stackrel{j}{-\!\!\!-\!\!\!-})$ represent the paths associated with the stochastic semigroup $X^j$. 
A synthetic description of the composition associated with the graph in figure~\ref{fig-coalescents-0} is given below in figure~\ref{fig-states-Ta}.
 \begin{figure}[H]
\centering
 $\widehat{\Ta}^{\,c_{k,j_1}}_{s,r}(x)=$
\begin{subfigure}[c]{0.2\textwidth}
\centering   
\resizebox{\linewidth}{!}{   
\xymatrix@C=4em@R=1.5em{
&\ar@<-.25ex>@{--}[l]_{k}\bullet&\ar@<-.25ex>@{-}[l]_{k}\circ\ar@<+.3ex>@{.>}[d]\bullet \\
&\ar@<-.25ex>@{--}[l]_{j_1}\bullet&\ar@<-.25ex>@{-}[l]_{j_1}~\bullet\\
&\ar@<-.25ex>@{--}[l]_{j_2}\bullet&\ar@<-.25ex>@{-}[l]_{j_2}~\bullet\\
0&s\ar@<-.25ex>@{-}[l] & \ar@<-.25ex>@{-}[l]r
}}
\end{subfigure} 
\centering \quad $\mbox{\rm and}\quad  \left(\widehat{\Ta}^{\,c_{j_2,j_1}}_{r,t}\circ \widehat{\Ta}^{\,c_{k,j_1}}_{s,r}\right)(x)=$
\begin{subfigure}[c]{0.3\textwidth}
\centering    
\resizebox{\linewidth}{!}{   
\xymatrix@C=4em@R=1.5em{
&\ar@<-.25ex>@{--}[l]_{k}\bullet&\ar@<-.25ex>@{-}[l]_{k}\circ\ar@<+.3ex>@{.>}[d]\bullet &\ar@<-.25ex>@{-}[l]_k\bullet&\\
&\ar@<-.25ex>@{--}[l]_{j_1}\bullet&\ar@<-.25ex>@{-}[l]_{j_1}~\bullet&\ar@<-.25ex>@{-}[l]_{j_1}~\bullet\\
&\ar@<-.25ex>@{--}[l]_{j_2}\bullet&\ar@<-.25ex>@{-}[l]_{j_2}\bullet &\ar@<-.25ex>@{-}[l]_{j_2}\circ\bullet\ar@<-.4ex>@{.>}[u]   \\
0& s\ar@<-.25ex>@{-}[l]& r\ar@<-.25ex>@{-}[l]& t\ar@<-.25ex>@{-}[l]& }}
\end{subfigure} 
 \caption{ }\label{fig-states-Ta}
\end{figure}

\subsection{Empirical measures}
 
We fix some integer $N\geq 2$ and for any $1\leq i<j\leq N$ and  $x=(x^i)_{1\leq i\leq N}$ on some $N$-fold cartesian product $E^N$  of some measurable space $(E,\Ea)$ we set
  \begin{eqnarray*}
x^{-i}&=&\left(x^1,\ldots,x^{i-1},x^{i+1},\ldots,x^N\right)\in E^{N}\\
x^{-\{i,j\}}&=&\left(x^1,\ldots,x^{i-1},x^{i+1},\ldots,x^{j-1},x^{j+1},\ldots,x^N\right)\in E^{N-2}
 \end{eqnarray*}
For any $1\leq i\leq N$ and  $x=(x^i)_{1\leq i\leq N}\in E^N$ we consider the functions
  \begin{eqnarray}
\varphi_{x^{-i}}~:~u\in E&\mapsto& \varphi_{x^{-i}}(u)=\left(x^1,\ldots,x^{i-1},u,x^{i+1},\ldots,x^N\right)\in E^N\nonumber\\
m~:~x\in S^N&\mapsto& m(x)=\frac{1}{N}\sum_{i\in [N]}\delta_{x^i}\in\Pa(E)\label{ref-varphi}
 \end{eqnarray}

 In this notation, the generator $\La_t$ of the stochastic flow $x\in S^N\mapsto \Xa_{s,t}(x)\in S^N$ introduced in (\ref{def-Xa-wXa}) is given for any
$F\in \Da(\La)$ and $x\in D_t(S)^N$ by the formula
$$
\La_t(F)(x)=\sum_{i\in [N]}L_{t}(F_{x^{-i}})(x^i)
$$
In the above display,  $\Da(\La)\subset \Ba(S^N)$ stands for the set of 
functions $F\in \Ba(S^N)$ s.t. for any $x\in S^N$ we have
$$
F_{x^{-i}}:=F\circ \varphi_{x^{-i}}\in \Da(L)
$$ 
 
Let 
$X=(X^i)_{1\leq i\leq N}$ be  $N$ independent random samples from some distribution $\eta\in\Pa(S)$. 
Using  (\ref{ref-mean-field-cv}) we have the first order expansion
  \begin{equation}\label{second-order-iid}
  \begin{array}{l}
\phi(m(X))-\phi(\eta)\\
\\
=(m(X)-\eta) \partial_\eta\phi
-\eta(Q^{m(X)}(1))\times (m(X)-\eta)(Q^{\eta}(1))\times(m(X)-\eta)\partial_\eta\phi
\end{array}
 \end{equation}
 Several estimates can be derived from the above decomposition. For instance using Cauchy-Schwartz inequality we have
 the bias estimate
 $$
 \log{(Q(1)(x)/Q(1)(y))}\leq q\Longrightarrow
N~ \vert\EE\left[\phi(m(X))(f)\right]-\phi(\eta)(f)\vert\leq e^q~ \mbox{\rm osc}(\overline{Q}(f))
 $$
 The stochastic perturbation analysis discussed above will be used repeatedly in section~\ref{sec-particle-flows} dedicated to the convergence analysis of the empirical measures associated with the genetic particle models discussed in section~\ref{sec-1-ref}.
 For instance, theorem~\ref{theo-key-decom} is the extended version of the second order expansion  (\ref{second-order-iid}) to interacting particle occupation measures.  
 
 \section{ A brief review on Feynman-Kac measures}
 \subsection{Exponential maps}
 For any $s\leq t$ and $N\geq 1$ we let $\Za_{s,t}$ be  the exponential map defined by
\begin{equation}\label{def-Za-xi-t}
\Za_{s,t}~:~x=(x_r)_{0\leq r\leq t}\in D_t(S)^N\mapsto
\Za_{s,t}(x):=\exp{\left(-\int_{s}^{t} m(x_r)(V_r)~dr\right)}\in [0,1]
\end{equation}
To simplify notation, when $s=0$ sometimes we write $\Za_{t}$ instead of $\Za_{0,t}$.
 When $N=1$, the map $\Za_{s,t}$ reduces to the exponential map ${Z}_{s,t}$  defined by
 \begin{equation}\label{def-Zt}
z=(z_r)_{0\leq r\leq t}\in D_t(S)\mapsto
Z_{s,t}(z):=\exp{\left(-\int_{s}^{t} V_r(z_r)~dr\right)}\in [0,1]\end{equation}
Observe that for any $s\leq t$ we have
$$
Z_{t}:=Z_{0,t}\Longrightarrow
Z_t(\widehat{X}_t)=Z_s(\widehat{X}_{s})~Z_{s,t}(\widehat{X}_t)\quad \mbox{\rm with the historical process $\widehat{X}_t=(X_s)_{0\leq s\leq t}$}
$$
To clarify the presentation when there are no possible confusions we write $Z_t(X)$ and $Z_{s,t}(X)$ instead of $Z_t(\widehat{X}_t)$ and $Z_{s,t}(\widehat{X}_t)$.

\subsection{Evolution semigroups}\label{ref-def-D}

Consider the flow of Feynman-Kac measures $(\gamma,\eta):t\in \RR_+:=[0,\infty[\mapsto(\gamma_t,\eta_t)\in (\Ma_+(S)\times\Pa(S))$ defined for any  $f\in \Ba(S)$ by the formulae
 \begin{equation}\label{ref-FK}
\eta_t(f)=\gamma_t(f)/\gamma_t(1)\quad\mbox{\rm with}\quad
\gamma_t(f):=\EE\left(f(X_t)~Z_t(X)\right)
\end{equation}
Observe that
$$
\log{\EE\left(Z_t(X)\right)}=- \int_0^t\eta_s(V_s)ds
$$
This shows that
\begin{equation}\label{def-over-V}
 \overline{Z}_t(X):={Z_t(X)}/\,{\EE(Z_t(X))}=\exp{\left[-\int_0^t \overline{V}_s(X_s)ds\right]}\quad
 \mbox{\rm with}\quad  \overline{V}_t:= V_t-\eta_t(V)
 \end{equation}
We also consider the Feynman-Kac semigroup
\begin{equation}\label{sg-Q}
Q_{s,t}(f)(x)=\EE\left(f(X_t)~Z_{s,t}(X)~|~X_s=x\right)\quad\mbox{\rm with}\quad Z_{s,t}(X):= Z_{t}(X)/ Z_{s}(X)
\end{equation}
When $V=0$ the semigroup $Q_{s,t}$ resumes to the Markov semigroup $P_{s,t}$ of the reference process $X_t$.

\begin{defi}
 The mathematical model $(\gamma_t,\eta_t,Q_{s,t})$ defined above and the measure $\QQ_t$ defined in (\ref{def-FK-intro}) is called the Feynman-Kac model associated with
the reference process and the potential function $(X_t,V_t)$.
\end{defi}

We further assume that the (infinitesimal) generators $L_t$ of $X_t$ are well defined on some common sub-algebra $\Da(S)\subset \Ba(S)$, and 
for any $s<t$ we have $Q_{s,t}(\Ba(S))\subset \Da(S)$. 

We let $\Va_t(f)=V_tf$ the multiplication operator on $\Ba(S)$. We also let 
$L_t=L^c_t+L^d_t$ be the decomposition of the generator $L_t$ in terms of a diffusion-type operator $L^c_t$ and a bounded jump-type generator of the
following form
$$
L^d_t(f)(u)=\lambda_t(u)~\int~(f(v)-f(u))~J_t(u,dv)
$$
for some bounded rate function $\lambda_t$ and some Markov transition $J_t$ on $S$. 

In this notation, for any $f\in \Da(S)$ and $s\leq t$ we have
\begin{equation}\label{sg-DaS}
\partial_t\gamma_t(f)=\gamma_t(L_t^V(f))\quad \mbox{\rm with}\quad L_t^V=L_t-\Va_t\Longrightarrow \gamma_t=\gamma_s Q_{s,t}
\end{equation}

The semigroup associated with the normalized Feynman-Kac measures $\eta_t$ is given for any $s\leq t$ by the formula
\begin{equation}\label{forward-FK-Gamma-red}
\eta_t=\phi_{s,t}(\eta_s):=\frac{\eta_s Q_{s,t}}{\eta_s Q_{s,t}(1)}\Longrightarrow \partial_t\eta_t(f)=\Lambda_{t}(\eta_t)(f):=\eta_t(L_t^V(f))+\eta_t(V_t)~\eta_t(f)
\end{equation}
with the collection of functional linear operators
$$
\Lambda_{t}(\eta)~:~f\in \Da(S)\mapsto
\Lambda_{t}(\eta)(f):=\eta(L_t^V(f))+\eta(V_t)~\eta(f)\in \RR
$$

Finally we recall that $\eta_t=\mbox{\rm Law}(Y_t)$ can be interpreted as the law of a nonlinear Markov process $Y_t$ associated with the collection of generators $L_{t,\eta}$ defined for any  $(\eta,f,x)\in (\Pa(S),\Da(S)\times S)$  by
 \begin{equation}\label{ref-mckv}
L_{t,\eta}(f)(x):=L_t(f)(x)+V_t(x)~\int~(f(y)-f(x))~\eta(dy)\quad\Longrightarrow\quad \Lambda_{t}(\eta)=\eta L_{t,\eta}
\end{equation}
\subsection{Path space measures}
\subsubsection{Historical processes}\label{sec-historical-FK}
 Consider a Feynman-Kac model $(\gamma^{\prime}_t,\eta^{\prime}_t,Q^{\prime}_{s,t},\QQ^{\prime}_t)$ associated with
some auxiliary Markov process $X^{\prime}_t$ on some metric space $(S^{\prime},d_{S^{\prime}})$, and
some bounded non negative potential functions $V^{\prime}_t$  on $S^{\prime}$. 
For instance, for any function $f$ on $S^{\prime}$ and any function $F$ on $D_t(S^{\prime})$ we have
$$
\eta_t^{\prime}(f)~\propto~\EE\left(f(X^{\prime}_t)~ \exp{\left[-\int_0^tV^{\prime}_s(X^{\prime}_s)ds\right]}\right)\quad \mbox{\rm and}\quad
 \QQ^{\prime}_t(F)~\propto~\EE\left(F(\widehat{X}^{\prime}_t)~ \exp{\left[-\int_0^t\widehat{V}^{\prime}_s(\widehat{X}^{\prime}_s)ds\right]}\right)
$$
with the historical process
$$
\widehat{X}^{\prime}_t:=\left(X^{\prime}_s\right)_{s\leq t}\quad \mbox{\rm and the potential function}\quad
\widehat{V}^{\prime}_t\left(\widehat{X}^{\prime}_t\right):=V^{\prime}_t(X^{\prime}_t).
$$
Assume that the process $X_t$ discussed in (\ref{ref-FK}) is given by the historical process 
\begin{equation}\label{FK-historical}
X_t:=\widehat{X}^{\prime}_t\in D_t(S^{\prime})\quad\mbox{\rm and set}\quad V_t(X_t):=V^{\prime}_t(X^{\prime}_t)
\end{equation}
In this situation, we have $\QQ_t^{\prime}=\eta_t$, with the measure $\eta_t$ defined in (\ref{ref-FK}); that is we have
$$
\QQ^{\prime}_t(F)~\propto~\EE\left(F(X_t)~ \exp{\left[-\int_0^tV_s(X_s)ds\right]}\right)~\propto~\eta_t(F)
$$

In the reverse angle, consider a collection of Feynman-Kac measures $(\widehat{\gamma}_t,\widehat{\eta}_t,\widehat{Q}_{s,t})$  defined as $(\gamma_t,\eta_t,Q_{s,t})$ by replacing in (\ref{ref-FK}) and (\ref{sg-Q}) the pair $(X_t,V_t)$ by 
 the historical process
\begin{equation}\label{FK-historical-w}
\widehat{X}_t:=\left(X_s\right)_{s\leq t}\in \widehat{S}:=\cup_{t\geq 0}D_t(S)\quad\mbox{\rm and the potential function}\quad \widehat{V}_t(\widehat{X}_t):=V_t(X_t)
\end{equation}
For instance, replacing $(X_t,V_t)$ by $(\widehat{X}_t,\widehat{V}_t)$ in (\ref{ref-FK}) and (\ref{sg-Q}) for any $F$ on $D_t(S)$ we have the formula
$$
\widehat{\eta}_t(F)~\propto~\EE\left(F(\widehat{X}_t)~ \exp{\left[-\int_0^t\widehat{V}_s(\widehat{X}_s)ds\right]}\right)
$$
This shows that $\widehat{\eta}_t=\QQ_t$, with the measure $\QQ_t$ introduced in (\ref{def-FK-intro}). We underline that $\QQ_t$ differs from the law of the the historical process $\widehat{Y}_t:=(Y_s)_{s\leq t}$ of the nonlinear process $Y_t$ discussed in the end of section~\ref{ref-def-D}.

To avoid unnecessary technical discussions we distinguish these models with the following terminology.

\begin{defi}
The Feynman-Kac model  associated with the historical process
$\widehat{X}_t$ and the potential function $\widehat{V}_t$ is called the path-space version of the Feynman-Kac model  associated with the process
${X}_t$ and the potential function ${V}_t$. In the reverse angle the Feynman-Kac model associated with the pair $(X_t,V_t)$ is called the marginal model of the Feynman-Kac model associated with the pair $(\widehat{X}_t,\widehat{V}_t)$.
\end{defi}
Let $L_t$ be the generator of $X_t$
defined on some common sub-algebra $\Da(L)\subset \Ba(S)$. 
In this situation, the generators $\widehat{L}_t$ of the historical process $\widehat{X}_t$ can be defined on some
common domain $\Da(\widehat{L})\subset \Ba(\widehat{S})$  in two different ways:

 The more conventional approach 
is to consider cylindrical functions
\begin{equation}\label{cylindrical-ref}
f(\widehat{X}_t):=\varphi(X_{s_1},\ldots,X_{s_n},X_t)
\end{equation}
 that only depend on a given finite collection of time horizons $s_i\leq s_{i+1}<t$, with $1\leq i<n$, and some bounded functions $\varphi$
 from $S^{n+1}$ into $\RR$. The regularity of the "test" function $\varphi$ depends on the process at hand.  For jump-type processes, no additional regularity is required. For diffusion-type processes
 the function is often required to be compactly supported and twice differentiable. 
 
 Another elegant and more powerful approach is to use the
 functional It\^o calculus theory developed by B. Dupire in an unpublished article~\cite{dupire} dating from 2009, and further developed in~\cite{cont,yuri}.
This path-dependent stochastic calculus allows to consider more general functions such as running integrals or running maximum of the process $X_t$.  It also allows to consider diffusion-type processes with a drift and a diffusion term that depends on the history of the process.

These path space models and their applications in the context of genealogical tree based particle models 
and historical processes were also presented in the early 2000s in previous work of one of the authors with Miclo~\cite{dm-2000,dm-2001-g,dm-2007}.

The path space $\widehat{S}$ is equipped with a time-space metric $d_{\widehat{S}}$ so that $(\widehat{S},d_{\widehat{S}})$ is a complete and separable metric space (cf. for instance proposition 1.1.13 and theorem 1.1.15 in~\cite{yuri-2}). The smoothness properties of continuous function $f$ on $S$
are defined in terms of  time and space functional derivatives. Thus, for diffusion-type historical processes $\widehat{X}_t$, the generator $\widehat{L}_t$ 
 is defined on functions $F\in\Da(\widehat{L})$ which a differentiable w.r.t. the time parameter and, as before
twice differentiable  with compactly supported derivatives (cf. for instance theorem 1.3.1 in~\cite{yuri-2}, as well as the pathwise change of variable formula stated in theorem 5.3.6 in~\cite{bally-Cont}). For instance, for the cylindrical functions discussed in (\ref{cylindrical-ref}) we have
$$
\widehat{L}_t(f)(\widehat{X}_t)=L_t\left(\varphi(X_{s_1},\ldots,X_{s_n},\point)\right)(X_t)
$$
as soon as the section $y\mapsto\varphi(x_1,\ldots,x_n,y)$ belongs to $\Da(L)$.
In the same vein, for path-integral functionals, we have
$$
f(\widehat{X}_t):=g(X_t)+\int_0^t~h(X_s)~ds\quad\Longrightarrow\quad \widehat{L}_t(f)(\widehat{X}_t)=L_t(g)(X_t)+h(X_t)
$$
as soon as  $g\in \Da(L)$ and $h\in\Ba(S)$. We also have the jump formula
$$
\Delta f(\widehat{X}_t):= f(\widehat{X}_t)- f(\widehat{X}_{t-})=
f(\widehat{X}_{t-}\vee X_t)- f(\widehat{X}_{t-})
$$
It is clearly not the scope of the article to describe in full details the functional It\^o calculus on path space. We rather refer the reader to the article~\cite{yuri}
and the Ph.D thesis of Saporito~\cite{yuri-2}. 
The discrete  time version of these historical processes calculus in the context of path-space  interacting particle systems can be found in~\cite{dm-2001-g}, see also~\cite{d-2004,d-2013} and references therein. 

In the further development of the article we shall use these ideas back and forth. We already mention that  the mean field particle interpretation of the Feynman-Kac measures associated with an historical process $\widehat{X}_t$ coincides with the genealogical tree-based particle evolution of the marginal model.

 Let  $\Da(\widehat{\La})$ be the set of 
functions $F$ on $\widehat{S}^N:=\cup_{t\geq 0}D_t(S)^N$ s.t. for any $x\in D_t(S)^N$ we have
$$
F_{x^{-i}}:=F\circ \varphi_{x^{-i}}\in \Da(\widehat{L})
$$ 
The generator $\widehat{\La}_t$ of the historical stochastic flow $\widehat{\Xa}_{s,t}(x)$ is given for any
$F\in \Da(\widehat{\La})$ and $x\in D_t(S)^N$ by the formula
$$
\widehat{\La}_t(F)(x)=\sum_{i\in [N]}\widehat{L}_{t}(F_{x^{-i}})(x^i)
$$
  We choose a generator  $\widehat{\La}_t^{(2)}$ on some domain $\Da(\widehat{\La}^{(2)})$  for a coupled stochastic flow $  \widehat{\Xa}_{s,t}^{(2)}$ of the following form
  $$
  \widehat{\Xa}_{s,t}^{(2)}(x,y):=(\widehat{\Xa}_{s,t}(x),\widehat{\Xa}_{s,t}(y))\quad \mbox{\rm with $(x,y)\in D_s(S)^{(2N)}:=\left( D_s(S)^{N}\times D_s(S)^{N}\right)$.}
  $$
For instance, we can define $ \widehat{\Xa}_{s,t}^{(2)}(x,y):=(\widehat{\Xa}_{s,t}(x),\widehat{\Xa}_{s,t}(y))$ firstly when $x_s=y_s$ by
$$
\widehat{\Xa}_{s,t}(x)=x\vee \widetilde{\Xa}_{s,t}(x)\qquad
\widehat{\Xa}_{s,t}(y)=y\vee \widetilde{\Xa}_{s,t}(y)
\quad \mbox{\rm with}\quad
\widetilde{\Xa}_{s,t}(x)=\widetilde{\Xa}_{s,t}(y)
$$
When $x_s\not=y_s$, the processes $\widehat{\Xa}_{s,t}(x)$ and $\widehat{\Xa}_{s,t}(y)$ are chosen independent.

\subsubsection{McKean measures}
\begin{defi}
    The McKean measure
   $\MM_t$ is the distribution on $D_t(S)$ of the historical process $\widehat{Y}_t:=(Y_s)_{s\leq t}$ of the nonlinear process  $Y_t$ discussed in (\ref{ref-mckv}). 
  \end{defi} 
   
   Observe that the $t$-marginal of $\MM_t$ coincides with $\eta_t$, that is we have 
   $$
   F(\widehat{Y}_t)=f(Y_t)\Longrightarrow \int F(y)~\MM_t(dy)=\EE(F(\widehat{Y}_t))=\EE(f(Y_t))=\int f(z)~ \eta_t(dz)
   $$
   The generator $\widehat{L}_{t,\eta_t}$ of $\widehat{Y}_t$ only depends on the marginal probability $\eta_t$ on $S$ and it is given for any function $F\in \Da(\widehat{L})$ and any $x=(x_s)_{s\leq t}\in D_{t}(S)$ by the formula
   $$
   \widehat{L}_{t,\eta_t}(F)(x)=\widehat{L}_t(F)(x)+\widehat{V}_t(x)~\int~(F(x\vee z)-F(x))~\eta_t(dz)
   $$
where $\widehat{L}_t$ stands for the generator of the historical process $   \widehat{X}_t:=(X_s)_{s\leq t}$ and $\widehat{V}_t$ is the potential function on path space given in (\ref{FK-historical-w}).

Following (\ref{FK-historical-w}), the Feynman-Kac measures $\QQ_t=\widehat{\eta}_t$ associated with the historical process and the potential function
$(\widehat{X}_t,\widehat{V}_t)$ on path space $\widehat{S}$ satisfy for any $F\in \Da(\widehat{L})$ the evolution equation
$$
\partial_t\QQ_t(F)=\QQ_t(\LL_{t,\QQ_t}(F))
$$
with the generator $\LL_{t,\QQ_t}$ defined for any $F\in \Da(\widehat{L})$ and $x\in D_t(S)$ by the formula
\begin{equation}\label{red-QQ-nonlinear}
\LL_{t,\QQ_t}(F)(x)=\widehat{L}_t(F)(x)+\widehat{V}_t(x)~\int~(F(y)-F(x))~\QQ_t(dy)
\end{equation}
Arguing as above, $\QQ_t=\mbox{\rm Law}(\overline{Y}_t)$ can be interpreted as the law of a nonlinear Markov process $\overline{Y}_t\in D_t(S)$ associated with the collection of generators $\LL_{t,\QQ_t}$. Here again, besides the fact that $\overline{Y}_t$ is a random path, we underline that $\overline{Y}_t$ is not the historical path of an auxiliary Markov process.

The historical process $\widehat{Y}_t$ is a jump type process taking values in the path space $\widehat{S}:=\cup_{t\geq 0}~D_t(S)$.  
Let $(T_n)_{n\geq 0}$ be a collection of jump times  occurring at a rate $\widehat{V}_{t}(\widehat{Y}_{t}):=V_{t}(Y_{t})$ that only depends on the terminal value of the historical process. We use the convention $T_0=0$ when $n=0$.
Also let $(\widehat{Y}^n_t)_{n\geq 0}$  be a collection of $D_t(S)$-valued independent random variables with common law $\QQ_t$. The terminal time variables $Y^n_t:=\widehat{Y}^n_t(t)$ are independent random variables with common law $\eta_t$.

Between the jump times $T_n$, the process $\widehat{Y}_t$ evolves as the historical process $\widehat{X}_t$ of  the reference process $X_t$. At each jump time $T_n$ the predictable path $\widehat{Y}_{T_n-}$ jumps onto   the  c\`adl\`ag path $  \widehat{Y}_{T_n}$ given by the concatenate path
$$
   \widehat{Y}_{T_n}:= \widehat{Y}_{T_n-}\vee Y^{n}_{T_n}
   $$
   In the above display, $ Y^{n}_{T_n}$ stands for a random sample from $\eta_{T_n}$
   independent of $ \widehat{Y}_{T_n-}$.
      Thus, for any $T_n\leq t<T_{n+1}$ we have
   $$
     \widehat{Y}_t=   \widehat{X}_{T_n,t}(   \widehat{Y}_{T_n})=
      \widehat{Y}_{T_n}\vee \widetilde{X}_{T_n,t} ( \widehat{Y}^n_{T_n})\quad\mbox{\rm and}\quad
       \widehat{Y}_{T_n}  = \widehat{Y}_{T_{n-1}}\vee \widetilde{X}_{T_{n-1},T_n}(\widehat{Y}^{n-1}_{T_{n-1}})
   $$

We summarize the above discussion with the following proposition.
\begin{prop}\label{prop-McK}
For any time horizon $t\geq 0$ we have
   $$
   \widehat{Y}_t=\sum_{n\geq 0}
   \left(\widetilde{X}_{T_0,T_1}(\widehat{Y}^0_{T_0})\vee \widetilde{X}_{T_1,T_2}(\widehat{Y}^1_{T_1})
   \vee\ldots\vee \widetilde{X}_{T_n,t}(\widehat{Y}^n_{T_n})\right)~1_{T_n\leq t<T_{n+1}}
   $$
The McKean measure $\MM_t$  is given for any measurable function $F$ on $D_t(S)$ by
   $$
   \begin{array}{l}
  \displaystyle\MM_t(F)
 =  \sum_{n\geq 0}\int_{[0,t]_n}~  \EE\left\{\left[
   \prod_{0\leq p< n}    Z_{r_p,r_{p+1}}\left(\widehat{X}_{r_p,r_{p+1}}(\widehat{Y}^{p}_{r_{p}})\right)~\widehat{V}_{r_{p+1}}\left(
   \widehat{X}_{r_{p},r_{p+1}}(\widehat{Y}^{p}_{r_{p}})\right)\right]~\right.\\
   \\
 \hskip 2cm \displaystyle  \left. \times~{Z}_{r_n,t}\left(\widehat{X}_{r_n,t}(\widehat{Y}^{n}_{r_{n}})\right)~~F\left(\widetilde{X}_{r_0,r_1}(\widehat{Y}^0_{r_0})\vee \ldots\vee\widetilde{X}_{r_n,t}(\widehat{Y}^n_{r_n})\right)\right\}~dr_1\ldots dr_n
   \end{array}
   $$

   \end{prop}
 
  When the function $F(\widehat{X}_t)=f(X_t)$ only depends on the terminal time, we recover the fact that
   \begin{equation}\label{MM-eta}
 \MM_t(F)=\gamma_t(f)~e^{\int_0^t\eta_s(V_s)~ds}=\EE\left(f(X_t)~e^{-\int_0^t(V_s(X_s)-\eta_s(V_s))ds}\right)=\EE(f(Y_t))=\eta_t(f)
   \end{equation}
     The detailed proof of the above assertion is provided in the appendix, on page~\pageref{MM-eta-proof}.

\subsection{Some regularity conditions}

This section discusses in some details the two main regularity conditions used in the further development of the article.

Firstly,  observe that the semigroup $P_{s,t}$ associated with the historical process $X_t=\left(X^{\prime}_s\right)_{s\leq t}$ discussed in (\ref{FK-historical}) never satisfies the regularity 
condition $(H_0)$ stated in (\ref{H0}). Nevertheless it may happen that the semigroup $P_{s,t}^{\prime}$ associated with $X^{\prime}_t$ satisfies condition $(H_0)$. In this situation, to avoid repetition or unnecessary long discussions we simply say that
$(H_0^{\prime})$ is met.

We also use  the following weaker conditions:
\begin{eqnarray}
(H_1)&&\exists \alpha<\infty\quad \exists \beta>0\nonumber\\
&& \mbox{\rm s.t.}\quad\forall s\leq t\quad\forall f\in \Ba(S)\quad
 \mbox{\rm osc}(\overline{Q}_{s,t}(f))\leq \alpha~e^{-\beta(t-s)}~\mbox{\rm osc}(f) \label{ref-H1}
\\
&&\nonumber\\
(H_2)&&\exists q<\infty\quad \mbox{\rm s.t.}\quad\forall s\leq t\quad\forall x,y\in S\qquad  \log{\left(Q_{s,t}(1)(x)/Q_{s,t}(1)(y)\right)}\leq q\label{ref-H2}
\end{eqnarray}

As before when the semigroup $\overline{Q}_{s,t}^{\prime}$ and $Q_{s,t}^{\prime}$ of associated with a marginal Feynman-Kac model satisfy condition $(H_i)$, to avoid repetition or unnecessary long discussions we simply say that
$(H_i^{\prime})$ is met.
We mention that
$$
(H_0)\Longrightarrow (H_1) \Longrightarrow  (H_2)
$$
The proof of the l.h.s. implication can be found in~\cite{dm-stab} (see for instance remark 2.2 and the contraction inequalities stated in proposition 2.3). In this context, the parameters $(\alpha,\beta)$ do not depends on the measure $\mu_{t,h}$ discussed in (\ref{H0}). 
To check the second implication observe that
\begin{equation}\label{ref-expo-Q}
\log{\left(Q_{s,t}(1)(x)/Q_{s,t}(1)(y)\right)}=\int_s^t~\left[\phi_{s,u}(\delta_y)(V_u)-\phi_{s,u}(\delta_x)(V_u)\right]~du
\end{equation}
This implies that
$$
(H_1)\Longrightarrow (H_2)\quad\mbox{\rm with}\quad q= \alpha \beta^{-1}~\mbox{\rm osc}(V)\quad\mbox{\rm
with}\quad \mbox{\rm osc}(V):=\sup_{t\geq 0}\mbox{\rm osc}(V_t)
$$
Using (\ref{ref-osc-partial-phi}) we also have
\begin{eqnarray}
(H_2)&\Longrightarrow&
\Vert \partial_{\eta}\phi_{s,t}(f) \Vert\leq e^q~\mbox{\rm osc}(f)\quad \left(\mbox{\rm since}~~ \mbox{\rm osc}(\overline{Q}_{s,t}(f))\leq~\mbox{\rm osc}(f) \right)\nonumber\\
(H_1)&\Longrightarrow&
\Vert \partial_{\eta}\phi_{s,t}(f) \Vert\leq r~e^{-\beta(t-s)}~\mbox{\rm osc}(f)\quad \mbox{\rm with}\quad
r= \alpha~ e^{q}\quad \mbox{\rm and}\quad q= \alpha \beta^{-1}~ \mbox{\rm osc}(V)\label{ref-r}
\end{eqnarray}
We return to the historical process  $X_t=\left(X^{\prime}_s\right)_{s\leq t}$ discussed in (\ref{FK-historical}). 
In this case, for any $x_s=(x^{\prime}(u))_{u\leq s}\in D_s(S^{\prime})$ we have
$$
Q_{s,t}(f)(x_s)=Q^{\prime}_{s,t}(f^{\prime})(x^{\prime}_s)
$$
in the above display, $f$ and $f^{\prime}$ stand for some bounded measurable functions on the path space $D_t(S^{\prime})$ and on $S^{\prime}$ such that 
$$
\forall y_t=(y^{\prime}_u)_{u\leq t}\in D_t(S^{\prime})\qquad f(y_t)=f^{\prime}(y^{\prime}_t)
$$
This implies that
\begin{equation}\label{ref-H2-historical}
(H_1^{\prime})\Longrightarrow
(H_2)\quad\mbox{\rm is met with}\quad q=\alpha\beta^{-1}~\mbox{\rm osc}(V)
\quad\mbox{\rm and}\quad
\Vert \partial_{\eta}\phi_{s,t}(f) \Vert\leq e^q~\mbox{\rm osc}(f)
\end{equation}

In the reverse angle, the Feynman-Kac semigroups on path space $\widehat{Q}_{s,t}$ discussed in (\ref{FK-historical-w})  is defined for any $s\leq t$, and any $x\in D_s(S)$ and $f\in \Ba(D_t(S))$ 
 by the formula
\begin{equation}\label{def-w-Q}
\widehat{Q}_{s,t}(f)(x):=
\EE\left(f(\widehat{X}_{t})~ Z_{s,t}\left(\widehat{X}_{t}\right)~|~\widehat{X}_{s}=x\right)
\end{equation}
with the exponential map $Z_{s,t}$ defined in (\ref{def-Zt}).
The above semigroup never satisfies the regularity 
condition $(H_1)$ stated in (\ref{H0}). Nevertheless arguing as in (\ref{ref-H2-historical}) permuting $(Q^{\prime}_{s,t},Q_{s,t})$ and $(Q_{s,t},\widehat{Q}_{s,t})$ 
we have 
\begin{equation}\label{ref-H-wH-2}
(H_1)\Longrightarrow (\widehat{H}_2)
\end{equation}
In the above display, $(\widehat{H}_2)$ stands for the condition defined as $(H_2)$ by replacing $Q_{s,t}$ by $\widehat{Q}_{s,t}$ and the state $S$ by the path space $D_s(S)$.

\subsection{Forward and backward equations}

\begin{prop}\label{prop-forward-backward-ref}
For any $s\leq t$ and $\eta\in \Pa(S)$ we have the Gelfand-Pettis forward and backward differential equations
\begin{equation}\label{forward-backward-ref}
 \partial_t\phi_{s,t}(\eta_s)=\Lambda_{t}(\phi_{s,t}(\eta_s))\quad \mbox{and}\quad
 \partial_s\phi_{s,t}(\eta)=-\Lambda_{s}(\eta)\partial_\eta\phi_{s,t}
\end{equation}
In addition, for any mapping $\phi$ of the form (\ref{ref-phi-Q}) we also have
\begin{equation}\label{forward-backward-ref-2}
\partial_t\phi\left(\phi_{s,t}(\eta)\right)=\Lambda_{t}(\eta)\partial_{\phi_{s,t}(\eta)}\phi
\quad \mbox{and}\quad\partial_s\phi\left(\phi_{s,t}(\eta)\right)=-\Lambda_{s}(\eta)\partial_{\phi_{s,t}(\eta)}\phi
\end{equation}
\end{prop}
 \proof
  The l.h.s. assertion in (\ref{forward-backward-ref}) is a direct consequence of (\ref{forward-FK-Gamma-red}).
  
 The r.h.s. differential formula in (\ref{forward-backward-ref}) as (\ref{forward-backward-ref-2})
 can be checked  using brute force and lengthy calculations from the backward evolution equations associated with Feynman-Kac semigroups.  
 
 A more judicious and more direct approach is 
 to apply a second oder perturbation approach.
For instance,  applying the second order Taylor with remainder expansion (\ref{ref-mean-field-cv}) to the mapping $\phi_{s,t}$ given in (\ref{forward-FK-Gamma-red}), for any $s+h\leq t$ we find that
  \begin{eqnarray*}
 \phi_{s,t}(\eta)&=&\phi_{s+h,t}\left(\eta+[\phi_{s,s+h}(\eta)-\eta]\right)\\
 &=&\phi_{s+h,t}(\eta)+[\phi_{s,s+h}(\eta)-\eta] \partial_\eta\phi_{s+h,t}
+\frac{1}{\phi_{s,s+h}(\eta) Q^{\eta}(1)}~\frac{1}{2}~[\phi_{s,s+h}(\eta)-\eta]^{\otimes 2} \partial_\eta^{2}\phi_{s+h,t}
 \end{eqnarray*}
On the other hand  we have
$$
\phi_{s,s+h}(\eta)=\eta+\Lambda_{s}(\eta)~h+\mbox{\rm O}(h^2)\quad\mbox{\rm and}\quad \phi_{s,s+h}(\eta) Q^{\eta}(1)=1+\mbox{\rm O}(h)
$$
This yields
 $$
 \phi_{s,t}(\eta)=\phi_{s+h,t}(\eta)+\Lambda_{s}(\eta)\partial_\eta\phi_{s+h,t}~~h+\mbox{\rm O}(h^2)
 $$
from which we check the backward evolution formula
$$
h^{-1}\left[ \phi_{s+h,t}\left(\eta\right)-\phi_{s,t}(\eta)\right]\longrightarrow_{h\rightarrow 0}\partial_s\phi_{s,t}(\eta)=-\Lambda_{s}(\eta)\partial_\eta\phi_{s,t}
 $$
For any mapping $\phi$ of the form (\ref{ref-phi-Q}), applying the third order Taylor with remainder expansion (\ref{ref-mean-field-cv}) we also have
 $$
   \begin{array}{l}
 \displaystyle \phi\left(\phi_{s+h,t}(\eta)\right)- \phi\left(\phi_{s,t}(\eta)\right)
=(\phi_{s+h,t}(\eta)-\phi_{s,t}(\eta)) \partial_{\phi_{s,t}(\eta)}\phi\\
  \\
  \displaystyle\hskip.1cm  +\frac{1}{2}~(\phi_{s+h,t}(\eta)-\phi_{s,t}(\eta))^{\otimes 2} \partial_{\phi_{s,t}(\eta)}^{2}\phi
+\frac{1}{ \phi_{s+h,t}(\eta)~Q^{\phi_{s,t}(\eta)}(1)}~\frac{1}{3}~(\phi_{s+h,t}(\eta)-\phi_{s,t}(\eta))^{\otimes 3} \partial_{\phi_{s,t}(\eta)}^{3}\phi
  \end{array}
 $$
 Arguing as above we check (\ref{forward-backward-ref-2}).
This ends the proof of the proposition.
\cqfd

\section{Interacting genetic type particle systems}
\subsection{Mean field particle processes}\label{sec-mean-field}

 \begin{defi}
The $N$-mean field particle interpretation of the nonlinear process discussed in (\ref{ref-mckv}) is defined by the Markov process
 $\xi_t=\left(\xi_t^i\right)_{1\leq i\leq N}\in S^N$ with generator $ \Ga_t$ given for any  $F\in  \Da(\La)$ and any $x=(x^i)_{1\leq i\leq N}\in S^N$
 by
 \begin{equation}
 \Ga_t(F)(x):=\sum_{i\in [N]}~L_{t,m(x)}(F_{x^{-i}})(x^i)\label{ref-mean-field-FK}
 \end{equation}
 with the collection of generators $L_{t,\mu}$ indexed by probability measures $\mu$ on $S$ defined in (\ref{ref-mckv}).
 \end{defi}

Let $\Fa:=(\Fa_t)_{t\geq 0}$, with
$\Fa_t=\sigma(\xi_s~:~s\leq t)$ be the filtration generated by the mean field particle model defined in (\ref{ref-mean-field-FK}).
Also let $\Da([0,T],S^N)$ be the set of functions $$F~:~(t,x)\in ([0,T]\times S^N)\mapsto F_t(x)\in \RR$$ with a bounded derivative w.r.t. the first argument  and s.t. $F_t\in \Da(\La)$, for any $t\in [0,T]$.
For any $F\in \Da([0,T],S^N)$, and any $T\geq 0$, applying the It\^o integration by part formula (see for instance section 15.5 in~\cite{dm-penev}) we have
 \begin{equation}\label{ref-ibp-mean-field}
dF_t(\xi_t)=\left[\partial_tF_t+\Ga_t(F_t)\right](\xi_t)~dt+d\Ma_t(F)
 \end{equation}
 In the above display
 $\Ma_t$ stands for a martingale random field on $ \Da([0,T],S^N)$ with angle bracket defined for any functions
 $F,G\in \Da([0,T],S^N)$ and any time horizon $t\in [0,T]$ by the formula
 \begin{equation}\label{ref-ibp-mean-field-2}
 \begin{array}{l}
\partial_t\langle \Ma(F),\Ma(G)\rangle_t=\Gamma_{\Ga_t}(F_t,G_t)(\xi_t)
\end{array}
 \end{equation}
 
Choosing functions of the form 
\begin{equation}\label{link-fF}
F_t(x)=m(x)(f_t)\quad\mbox{\rm and}\quad
G_t(x)=m(x)(g_t)\Longrightarrow \Gamma_{\Ga_t}(F_t,G_t)(\xi_t)=m(\xi_t)\Gamma_{L_{t,m(\xi_t)}}(f_t,g_t)
\end{equation}  
for some $f,g\in \Da([0,T],S)$,
we also check that the occupation measure $m(\xi_t)\in\Pa(S)$ satisfies the stochastic equation
\begin{equation}\label{def-Mt}
dm(\xi_t)(f_t)=\left[m(\xi_t)(\partial_tf_t)+\Lambda_{t}(m(\xi_t))(f_t)\right]~dt+ \frac{1}{\sqrt{N}}~dM_t(f)
\end{equation} 
with a martingale random field $M_t$ on $\Da([0,T],S)$  with angle brackets by the formula
$$
 \begin{array}{l}
\partial_t\langle M(f),M(g)\rangle_t\\
\\
\displaystyle=m(\xi_t)(\Gamma_{L_t}(f_t,g_t))+\int~m(\xi_t)(du)~m(\xi_t)(dv)~V_t(u)~(f_t(v)-f_t(u))(g_t(v)-g_t(u))~
\end{array}$$
With a slight abuse of notation we also write $M_t$ the extension of the random field $M_t$ to
$\Fa$-predictable functions $\Da([0,T],S)$.

In the further development of the article we write $(M^c_t,\Ma^c_t)$ and $(M^d_t,\Ma^d_t)$ the continuous and the discontinuous part of the martingales $(M_t,\Ma_t)$; as well as
$$
L_{t,\eta}=L_{t,\eta}^c+L_{t,\eta}^d\quad\mbox{\rm with}\quad L_{t,\eta}^c:=L_t^c
$$
The angle bracket of $\Ma^d_t$ is given for any  functions  $F,G\in \Da([0,T],S^N)$ and any time horizon $t\in [0,T]$ by the formula
$$
 \begin{array}{l}
\partial_t\langle \Ma^d(F),\Ma^d(G)\rangle_t\\
\\
\displaystyle=\sum_{i\in [N]}\int~\left[F_{t,\xi_t^{-i}}(v)-F_t(\xi_t)\right]\left[G_{t,\xi_t^{-i}}(v)-G_t(\xi_t)\right]~\left[V_t(\xi^i_t)~m(\xi_t)(dv)+
\lambda_t(\xi^i_t)~J_t(\xi^i_t,dv)\right]
\end{array}
$$

  \subsection{Historical and genealogical tree evolutions}\label{hist-gen-sec}
      
      The generator $ \Ga_t$  of the $N$-mean field particle process $\xi_t=\left(\xi_t^i\right)_{1\leq i\leq N}\in S^N$ discussed in (\ref{ref-mean-field-FK}) can be rewritten for any  $F\in  \Da(\La)$ and $x=(x^i)_{1\leq i\leq N}\in S^N$ in terms of the interacting jump operator
 \begin{equation}
 \Ga_t(F)(x)=\sum_{i\in [N]}~\left[  L_t(F_{x^{-i}})(x^i)+ V_t(x^i)~\int~(F_{x^{-i}}(y)-F(x))~m(x)(dy) \right]\label{ref-mean-field-FK-genetic}
 \end{equation}

 The process $\xi_t$ is called the $N$-mean field particle interpretation of the Feynman-Kac measures $\eta_t$
 defined in (\ref{ref-FK}).

 The process $\xi_t=(\xi_t^i)_{i\in [N]}\in S^N$ can be interpreted as a genetic type particle system with mutation associated with the generator $L_t$ and selection dictated by the potential function $V_t$. Between the jumps the particles $\xi^i_t$ evolve independently as independent copies of a process with generator $L_t$. At rate $V_t(\xi^i_t)$ the particle is killed and instantly a particle randomly selected in the pool duplicates.
 
   The historical process $\widehat{\xi}_t=\left(\xi_s\right)_{s\leq t}\in \widehat{S}^N$ coincides with the 
   $N$-mean field interpretation of the historical process $\widehat{Y}_t=(Y_s)_{s\leq t}\in \widehat{S}$ of the nonlinear process $Y_t$ defined in (\ref{ref-mckv}), with the path space $\widehat{S}$ defined in (\ref{FK-historical-w}).
   
The generator $\widehat{\Ga}_t$ of the process $\widehat{\xi}_t$ given for any  function $F\in \Da(\widehat{\La})$ and any  $x=(x_s)_{s\leq t}\in D_t(S)^N$
 by the  operator
\begin{equation}\label{ref-wG-m}
  \widehat{\Ga}_t(F)(x)
  =\sum_{i\in [N]}~\left[  \widehat{L}_t(F_{x^{-i}})(x^i)+\widehat{V}_t(x^i)~\int~(F_{x^{-i}}(x^i\vee z)-F(x))~m(x_t)(dz) \right]
 \end{equation}

 The  genealogical tree evolution associated with the genetic type particle system $\xi_t\in S^N$ coincides with
 the $N$-mean field interpretation  $\XX_t=(\XX^i_t)_{i\in [N]}\in D_t(S)^N$ of the $D_t(S)$-valued nonlinear process $\overline{Y}_t$ defined in (\ref{red-QQ-nonlinear}). The generator $\GG_t$ of  $\XX_t$ is given  for any function $F\in \Da(\widehat{\La})$  and  $x=(x_s)_{s\leq t}\in D_t(S)^N$
by the operator
 \begin{eqnarray}
   \GG_t(F)(x)&:=&\sum_{i\in [N]}~ \LL_{t,m(x)}(F_{x^{-i}})(x^i)\nonumber\\
   &=&\sum_{i\in [N]}~\left[  \widehat{L}_t(F_{x^{-i}})(x^i)+ \widehat{V}_t(x^i)~\int~(F_{x^{-i}}(y)-F(x))~m(x)(dy) \right]\label{ref-GG-GG-k}
\end{eqnarray}
with the collection of generators $\LL_{t,\mu}$ indexed by probability measures $\mu$ on $D_t(S)$ defined in (\ref{red-QQ-nonlinear}).

 Observe that $ \GG_t$ can be deduced from the operator $ \Ga_t$ by replacing in $(\ref{ref-mean-field-FK})$ the generator and the potential function
 $(L_t,V_t)$ and the state $S$ by $(\widehat{L}_t,\widehat{V}_t)$ and by the path space $\widehat{S}$.

  The process $\XX_t$ is called the $N$-mean field particle interpretation of the Feynman-Kac measures on path space $\widehat{\eta}_t=\QQ_t$
 defined in (\ref{red-QQ-nonlinear}).

In contrast with $\XX_t$, the historical process $\widehat{\xi}_t$ keeps track of all 
past-ancestral lines with no descendants. Thus $\widehat{\xi}_t$ can be interpreted as the complete ancestral tree associated with the genetic type particle system discussed above.

Rewritten in a slightly different form the jump generator  in (\ref{ref-wG-m}) is given by the formula
 $$
 \begin{array}{l}
  \displaystyle \widehat{V}_t(x^i)~\int~(F_{x^{-i}}(x^i\vee z)-F(x))~m(x_t)(dz)\\
 \\
 \displaystyle=\widehat{V}^-_t(x^i)~\frac{1}{N-1}\sum_{k\in [N]-\{i\}}
\left(F(x\vee x_t^{c_{i,k}})-F(x)\right)\quad\mbox{with}\quad
\widehat{V}^-_t:=\left(1-\frac{1}{N}\right)~\widehat{V}_t
 \end{array}$$
 This shows that  for any $i\in [N]$ and at rate $\widehat{V}_t(\widehat{\xi}^{\,i}_t)$ the $i$-th particle $\xi^i_t$ jumps onto a randomly chosen particle $\xi^k_t$, with $k\in [N]-\{i\}$. Also observe that the jump rate of the system $\widehat{\xi}_t$ is given by the stochastic intensity
 \begin{equation}\label{def-lambda-intensity}
\widehat{\lambda}_t(\widehat{\xi}_t)=N~m(\widehat{\xi}_t)(\widehat{V}_t^-)=(N-1)~m(\widehat{\xi}_t)(\widehat{V}_t)
 \end{equation}
 A graphical description of the historical process $\widehat{\xi}_t$ with three jump times $(T_1,T_2,T_3)$ and the jump coalescent maps  $(c_{k,j_1},c_{j_2,j_1},c_{k,j_2})$ is given below in figure~\ref{figure-xi}.
\begin{figure}[H]
\centering  $$\widehat{\xi}_t=\left(
\begin{subfigure}[c]{0.5\textwidth}
\centering    
\resizebox{\linewidth}{!}{   
\xymatrix@C=4em@R=1.5em{
\bullet&\ar@<-.25ex>@{-}[l]_{k}\circ\ar@<+.3ex>@{.>}[d]\bullet &\ar@<-.25ex>@{-}[l]_k\bullet&\ar@<-.25ex>@{-}[l]_k\circ\ar@<+.7ex>@{.>}[dd]\bullet&\ar@<-.25ex>@{-}[l]_{k}\bullet\\
\bullet&\ar@<-.25ex>@{-}[l]_{j_1}~\bullet&\ar@<-.25ex>@{-}[l]_{j_1}~\bullet&\ar@<-.25ex>@{-}[l]_{j_1}\bullet~&\ar@<-.25ex>@{-}[l]_{j_1}\bullet\\
\bullet&\ar@<-.25ex>@{-}[l]_{j_2}\bullet &\ar@<-.25ex>@{-}[l]_{j_2}\circ\bullet\ar@<-.4ex>@{.>}[u]&\ar@<-.25ex>@{-}[l]_{j_2}~~\bullet&\ar@<-.25ex>@{-}[l]_{j_2}\bullet\\
0 & T_1\ar@<-.25ex>@{-}[l]& T_2  \ar@<-.25ex>@{-}[l]&  \ar@<-.25ex>@{-}[l] T_3 &\ar@<-.25ex>@{-}[l]t 
}}
\end{subfigure} 
\right)$$
 \caption{Graphical description of the historical process $\widehat{\xi}_t$}\label{figure-xi}
\end{figure}
The plain lines $(\stackrel{j}{-\!\!\!-\!\!\!-})$ represent the path associated with N independent copies $X^j_{u,v}$ of the stochastic semigroup $X_{u,v}$ of the process $X_t$ with generator $L_t$.  The dotted vertical arrows represent the jumps from 
the $i_n$-th coordinate onto the $j_n$-th one. 

The genealogical tree $\XX_t$ associated with the historical process $\widehat{\xi}_t$ in figure~\ref{figure-xi}
is given below in figure~\ref{fig-coalescents}.

  \begin{figure}[H]
\centering
\begin{subfigure}[c]{0.5\textwidth}
\centering    
\resizebox{\linewidth}{!}{   
\xymatrix@C=4em@R=1.5em{
\bullet&\ar@<-.25ex>@{-}[l]_{k}\circ\ar@<-.25ex>@{.>}[d]&\bullet\ar@<-.2ex>@{-}[ld]_k&\ar@<-.25ex>@{.>}[ddd]\circ\ar@<-.25ex>@{-}[l]_{k}\\
&\bullet\ar@<.2ex>@2{-}[ld]_{j_1}&&~\bullet&\ar@<.2ex>@{.}[l]_{j_1}\\
\bullet&&\bullet\ar@{-}[ru]^{j_1}\ar@<.25ex>@2{-}[lu]_{j_1}& & \ar@<-.25ex>@2{-}[ld]_{k}   \\
\bullet&\ar@<-.25ex>@{-}[l]_{j_2}\bullet&\ar@<-.25ex>@{-}[l]_{j_2}\circ\ar@<-.25ex>@{.>}[u]&\ar@<.25ex>@2{-}[lu]^{j_2}\bullet  &    \ar@<-.25ex>@{-}[l]^{j_2} &     \\
0 & T_1\ar@<-.25ex>@{-}[l]& T_2  \ar@<-.25ex>@{-}[l]&  \ar@<-.25ex>@{-}[l] T_3 &\ar@<-.25ex>@{-}[l]t&  
}}
 \end{subfigure} 
 \caption{Description of the genealogical tree associated with the particle system $\xi_t$}\label{fig-coalescents}
\end{figure}
In the above display the $k$-th ancestral line $\XX^k_t$ is represented by the doubled lines $(\xymatrix@1{&\ar@2{-}[l]})$. A graphical description of the genealogical tree  $\XX_t$ associated with the particle evolution presented in figure~\ref{fig-coalescents} is given below, in figure~\ref{figure-XX-process}.
 \begin{figure}[H]
\centering  $$\XX_t=\left(
\begin{subfigure}[c]{0.3\textwidth}
\centering    
\resizebox{\linewidth}{!}{   
\xymatrix@C=4em@R=1.5em{
\bullet&\ar@<-.25ex>@{-}[l]_{j_1}\bullet &\ar@<-.25ex>@{-}[l]_{j_1}\bullet&\ar@<-.25ex>@{-}[l]_{j_2}\bullet&\ar@<-.25ex>@{-}[l]_{k}\bullet\\
\bullet&\ar@<-.25ex>@{-}[l]_{j_1}~\bullet&\ar@<-.25ex>@{-}[l]_{j_1}~\bullet&\ar@<-.25ex>@{-}[l]_{j_1}\bullet~&\ar@<-.25ex>@{-}[l]_{j_1}\bullet\\
\bullet&\ar@<-.25ex>@{-}[l]_{j_1}\bullet &\ar@<-.25ex>@{-}[l]_{j_1}&\ar@<-.25ex>@{-}[l]_{j_2}~~\bullet&\ar@<-.25ex>@{-}[l]_{j_2}\bullet\\
0 & T_1\ar@<-.25ex>@{-}[l]& T_2  \ar@<-.25ex>@{-}[l]&  \ar@<-.25ex>@{-}[l] T_3 &\ar@<-.25ex>@{-}[l]t 
}}
\end{subfigure} 
\right)$$
 \caption{Graphical description of the genealogical tree $\XX_t$}\label{figure-XX-process}
\end{figure}
\subsection{Particle models with a frozen ancestral line}\label{sec-particle-frozen-line}

\begin{defi}
    Let $\Ia$ be an uniform random sample on $ [N]$.
Given $\Ia$, let  $\widehat{\goodchi}_t\in D_t(S)^N$ 
 be the Markov process with initial condition $\widehat{\goodchi}_0=\xi_0$ and generator $ \widehat{\Ga}_{\Ia,t}$ defined for any function $F\in \Da(\widehat{\La})$  and $x=(x_s)_{s\leq t}\in D_{t}(S)^N$, by the formula
\begin{equation}\label{ref-wG-m-k}
\begin{array}{l}
   \displaystyle \widehat{\Ga}_{\Ia,t}(F)(x)=\sum_{i\in [N]-\{\Ia\}}~\left[  \widehat{L}_t(F_{x^{-i}})(x^i)+\widehat{V}_t(x^i)~\int~(F_{x^{-i}}(x^i\vee z)-F(x))~m(x_t)(dz) \right]\\
  \\
 \hskip3cm \displaystyle +~ \widehat{L}_t(F_{x^{-\Ia}})(x^{\Ia})+\frac{1}{N}~\sum_{i\in [N]-\{\Ia\}}~\widehat{V}_t(x^i)~(F(x^{\sigma_{i,\Ia}}\vee x^{c_{i,\Ia}}_t)-F(x))
\end{array}
 \end{equation}
with the concatenate operation $\vee$ defined in (\ref{ref-concatenate}). 
\end{defi}

Observe that  $ \widehat{\goodchi}_{t}:=( \widehat{\goodchi}_{t}(s))_{s\leq t}\in D_t(S)^N$ is not necessarily the historical path of an auxiliary Markov process.

To better connect the above generator with the operator $\widehat{\Ga}_{t}$ defined in (\ref{ref-wG-m}) observe that
$$
\begin{array}{l}
\widehat{\Ga}_{t}(F)(x)-\widehat{\Ga}_{k,t}(F)(x)
\\
\\
=\displaystyle\frac{1}{N}~\widehat{V}_t(x^k)~\sum_{i\in [N]-\{k\}}(F(x\vee x_t^{c_{k,i}})-F(x))-\frac{1}{N}~\sum_{i\in [N]-\{k\}}~\widehat{V}_t(x^i)~(F(x^{\sigma_{i,k}}\vee x^{c_{i,k}}_t)-F(x))
\end{array}
$$
This yields the following lemma.
\begin{lem}\label{lem-connect-generators}
The generators $\widehat{\Ga}_{t}$ and $\widehat{\Ga}_{k,t}$   defined in (\ref{ref-wG-m}) and  (\ref{ref-wG-m-k}) 
are connected for any function $F\in \Da(\widehat{\La})$  and $x=(x_s)_{s\leq t}\in D_{t}(S)^N$ by the formula
$$
\widehat{\Ga}_{t}(F)(x)-m(x)(\widehat{V}_t)~F(x)=\widehat{\Ga}_{k,t}(F)(x)-\widehat{V}_t(x^k)~F(x)+\frac{1}{N}~\sum_{i\in [N]-\{k\}}~\epsilon_{k,i}(F)(x)
$$
with the bounded integral operator $\epsilon_{i,k}$ defined by the anti-symmetric functions
$$
\epsilon_{k,i}(F)(x):=\widehat{V}_t(x^k)~F(x\vee x_t^{c_{k,i}})-\widehat{V}_t(x^i)~F(x^{\sigma_{i,k}}\vee x^{c_{i,k}}_t)=-\epsilon_{k,i}(F)(x^{\sigma_{i,k}})
$$
\end{lem}

Nevertheless, given $\Ia=k$, the end-state process $\zeta_t$ defined by 
  $$\zeta_t:=\widehat{\goodchi}_t(t)\in S^N
  $$ 
  is a Markov process. Given $\Ia=k$, its generator $ \Ga_{k,t}$ is
 given for any  $F\in  \Da(\La)$ and $x=(x^i)_{1\leq i\leq N}\in S^N$
 by the formula
 \begin{equation}\label{ref-zeta-GG}
\begin{array}{l}
   \displaystyle \Ga_{k,t}(F)(x)=\sum_{i\in [N]-\{k\}}~\left[  L_t(F_{x^{-i}})(x^i)+V_t(x^i)~\int~(F_{x^{-i}}(z)-F(x))~m(x)(dz) \right]\\
  \\
 \hskip5cm \displaystyle + ~L_t(F_{x^{-k}})(x^k)+\frac{1}{N}~\sum_{i\in [N]-\{k\}}~V_t(x^i)~(F_{x^{-i}}(x^k)-F(x))
\end{array}
\end{equation}
Given $\Ia$, the genealogical tree $ \YY_{t}$ associated with the process $\widehat{\goodchi}_t$  is the $\in D_t(S)^N$-valued process starting at $\YY_{0}=\xi_0$ with generator
$  \GG_{\Ia,t}$ given  for any  function $F\in \Da(\widehat{\La})$  and $x=(x_s)_{s\leq t}\in D_{t}(S)^N$ by the formula
\begin{equation}\label{ref-GG-k-mn}
\begin{array}{l}
 \displaystyle \GG_{\Ia,t}(F)(x)=\sum_{i\in [N]-\{\Ia\}}~\left[  \widehat{L}_t(F_{x^{-i}})(x^i)+\widehat{V}_t(x^i)~\int~(F_{x^{-i}}(z)-F(x))~m(x)(dz) \right]\\
  \\
 \hskip3cm \displaystyle + \widehat{L}_t(F_{x^{-\Ia}})(x^{\Ia})+\frac{1}{N}~\sum_{i\in [N]-\{\Ia\}}~\widehat{V}_t(x^i)~(F_{x^{-i}}(x^\Ia)-F(x))
\end{array}
\end{equation}

As in~(\ref{ref-GG-GG-k}),  $\GG_{k,t}$ can be deduced from the operator $ \Ga_{k,t}$ by replacing in (\ref{ref-zeta-GG})
the generator and the potential function
 $(L_t,V_t)$ and the state $S$ by $(\widehat{L}_t,\widehat{V}_t)$ and by the path space $\widehat{S}$.

Observe that the historical process  $ \widehat{\zeta}^{\Ia}_t=(\zeta^{k}_s)_{ s\leq t}$ of the $\Ia$-th particle in $\zeta_t:=(\zeta_t^i)_{i\in [N]}$ evolves as the historical process $\widehat{X}_t$ and doesn't depends on the 
remaining $(N-1)$ interacting particles.

This yields the following proposition.
 \begin{prop}
  The $\Ia$-th ancestral line $\YY^{\Ia}_t$ of the genealogical tree $\YY_t:=(\YY^i_t)_{i\in [N]}\in D_t(S)^N$ coincides with the historical process  $ \widehat{\zeta}^{\Ia}_t:=(\zeta^{\Ia}_s)_{ s\leq t}$ of the $\Ia$-th particle $\zeta^{\Ia}_t$; that is, we have that
  \begin{equation}\label{ref-YY-w-zeta}
  \YY^{\Ia}_t= \widehat{\zeta}^{\Ia}_t\stackrel{law}{=}\widehat{X}_t\in D_t(S)
  \end{equation}
  \end{prop}
  
Observe that the jump rate of the system $\widehat{\goodchi}_t$ is given by the stochastic intensity
  \begin{equation}\label{def-lambda-Ia-k}
\widehat{\lambda}_{\Ia,t}(\widehat{\goodchi}_t)= (N-1)~m(\widehat{\goodchi}_t^{-\Ia})(\widehat{V}_t)
 \end{equation}

  A synthetic description of the evolution of $\widehat{\goodchi}$ given $\Ia=k$ and the  jumps 
 $$
x\leadsto  x^{\sigma_{j_1,k}}\vee x^{c_{j_1,k}}_t\qquad
 x\leadsto x\vee x^{c_{j_2,j_1}}_t\quad \mbox{\rm and}   \quad x\leadsto   x^{\sigma_{j_2,k}}\vee x^{c_{j_2,k}}_t
 $$ 
on three jump times $(T_{k,1},T_{k,2},T_{k,3})$
  is given below in figure~\ref{figure-xi-3-goodchi}.
  \begin{figure}[H]
\centering 
\begin{subfigure}[c]{0.11\textwidth}
\centering    
\resizebox{\linewidth}{!}{   
\xymatrix@C=4em@R=1.5em{
\bullet&\ar@<-.25ex>@2{-}[l]_{k}\bullet\\
\bullet&\ar@<-.25ex>@{-}[l]_{j_1}\bullet\\
\bullet&\ar@<-.25ex>@{-}[l]_{j_2}\bullet
}}
\end{subfigure} $\leadsto$
\centering 
\begin{subfigure}[c]{0.19\textwidth}
\centering    
\resizebox{\linewidth}{!}{   
\xymatrix@C=4em@R=1.5em{
&\ar@<-.25ex>@{-}[l]_{j_1}\circ\bullet\ar@<.4ex>@{.>}[d] &\ar@<-.25ex>@2{-}[l]_{k}\bullet\\
\bullet&\ar@<-.25ex>@2{-}[l]_{k}\bullet&\ar@<-.25ex>@{-}[l]_{j_1}\bullet \\
\bullet&\ar@<-.25ex>@{-}[l]_{j_2}\bullet&\ar@<-.25ex>@{-}[l]_{j_2}\bullet
}}
\end{subfigure}
$\leadsto$
\centering 
\begin{subfigure}[c]{0.26\textwidth}
\centering    
\resizebox{\linewidth}{!}{   
\xymatrix@C=4em@R=1.5em{
\bullet&\ar@<-.25ex>@{-}[l]_{j_1}\circ\bullet\ar@<.4ex>@{.>}[d]  &\ar@<-.25ex>@2{-}[l]_{k}\bullet&\ar@<-.25ex>@2{-}[l]_{k} \bullet\\
\bullet&\ar@<-.25ex>@2{-}[l]_{k}\bullet&\ar@<-.25ex>@{-}[l]_{j_1}\bullet  &\ar@<-.25ex>@{-}[l]_{j_1}\bullet\\
\bullet&\ar@<-.25ex>@{-}[l]_{j_2} \bullet&\ar@<-.25ex>@{-}[l]_{j_2}\circ\bullet \ar@<-.4ex>@{.>}[u]&\ar@<-.25ex>@{-}[l]_{j_2}\bullet
}}
\end{subfigure}
$\leadsto$
\centering 
\begin{subfigure}[c]{0.32\textwidth}
\centering    
\resizebox{\linewidth}{!}{   
\xymatrix@C=4em@R=1.5em{
\bullet&\ar@<-.25ex>@{-}[l]_{j_2}\bullet &\ar@<-.25ex>@{-}[l]_{j_2}\circ\bullet\ar@<.7ex>@{.>}[d]&\ar@<-.25ex>@{-}[l]_{j_2}\circ\bullet \ar@<.7ex>@{.>}[dd] &\ar@<-.25ex>@2{-}[l]_{k}\bullet \\
\bullet&\ar@<-.25ex>@2{-}[l]_{k}\bullet&\ar@<-.25ex>@{-}[l]_{j_1}\bullet  &\ar@<-.25ex>@{-}[l]_{j_1}\bullet&\ar@<-.25ex>@{-}[l]_{j_1}\bullet \\
\bullet&\ar@<-.25ex>@{-}[l]_{j_1} \circ\bullet\ar@<-.4ex>@{.>}[u]&\ar@<-.25ex>@2{-}[l]_{k}~~\bullet &\ar@<-.25ex>@2{-}[l]_{k}\bullet&\ar@<-.25ex>@{-}[l]_{j_2}\bullet 
}}
\end{subfigure}
 \caption{  }\label{figure-xi-3-goodchi}
\end{figure}
In the above display the $k$-th ancestral line $\YY^k_t$ is represented by the doubled lines $(\xymatrix@1{&\ar@2{-}[l]})$.

A graphical description of the corresponding process $\widehat{\chi}_t$   is given below in figure~\ref{figure-good-chi}.
 \begin{figure}[H]
\centering  $$\widehat{\goodchi}_t=\left(
\begin{subfigure}[c]{0.5\textwidth}
\centering    
\resizebox{\linewidth}{!}{   
\xymatrix@C=4em@R=1.5em{
\bullet&\ar@<-.25ex>@{-}[l]_{j_2}\bullet &\ar@<-.25ex>@{-}[l]_{j_2}\circ\bullet\ar@<.7ex>@{.>}[d]&\ar@<-.25ex>@{-}[l]_{j_2}\circ\bullet \ar@<.7ex>@{.>}[dd] &\ar@<-.25ex>@2{-}[l]_{k}\bullet \\
\bullet&\ar@<-.25ex>@2{-}[l]_{k}\bullet&\ar@<-.25ex>@{-}[l]_{j_1}\bullet  &\ar@<-.25ex>@{-}[l]_{j_1}\bullet&\ar@<-.25ex>@{-}[l]_{j_1}\bullet \\
\bullet&\ar@<-.25ex>@{-}[l]_{j_1} \circ\bullet\ar@<-.4ex>@{.>}[u]&\ar@<-.25ex>@2{-}[l]_{k}~~\bullet &\ar@<-.25ex>@2{-}[l]_{k}\bullet&\ar@<-.25ex>@{-}[l]_{j_2}\bullet \\
0 & T_{k,1}\ar@<-.25ex>@{-}[l]& T_{k,2}  \ar@<-.25ex>@{-}[l]&  \ar@<-.25ex>@{-}[l] T_{k,3} &\ar@<-.25ex>@{-}[l]t 
}}
\end{subfigure} 
\right)$$
 \caption{Graphical description of the process $\widehat{\goodchi}_t$}\label{figure-good-chi}
\end{figure}

The genealogical tree $\YY_t$ associated with the historical process $\widehat{\goodchi}_t$ in figure~\ref{figure-good-chi}
is given below in figure~\ref{fig-coalescents-goodchi}.
.
  \begin{figure}[H]
\centering
\begin{subfigure}[c]{0.5\textwidth}
\centering    
\resizebox{\linewidth}{!}{   
\xymatrix@C=4em@R=1.5em{
\bullet&\ar@<-.25ex>@{-}[l]_{j_1}\circ\ar@<-.25ex>@{.>}[d]&\bullet\ar@<-.2ex>@{-}[ld]_{j_1}\ar@{-}[rd]^{j_2}&\bullet\ar@<-.25ex>@{-}[l]_{j_1}&\bullet\ar@<-.25ex>@{-}[l]_{j_1}\\
&\bullet\ar@<.2ex>@2{-}[ld]_{k}&&~\circ\ar@<-.25ex>@{.>}[dd]&\\
\bullet&&\bullet\ar@<.25ex>@2{-}[lu]_{k}& & \ar@<-.2ex>@{-}[ld]_{j_2} \\
\bullet&\ar@<-.25ex>@{-}[l]_{j_2}\bullet&\ar@<-.25ex>@{-}[l]_{j_2}\circ\ar@<-.25ex>@{.>}[uuu]&
\ar@<.25ex>@2{-}[lu]^{k}\bullet  &   \ar@<-.25ex>@2{-}[l]^{k} \bullet&     \\
0 & T_{k,1}\ar@<-.25ex>@{-}[l]& T_{k,2}  \ar@<-.25ex>@{-}[l]&  \ar@<-.25ex>@{-}[l] T_{k,3} &\ar@<-.25ex>@{-}[l]t&  
}}
 \end{subfigure} 
 \caption{Genealogical tree associated with the process $\widehat{\goodchi}_t$ in figure~\ref{figure-good-chi}}\label{fig-coalescents-goodchi}
\end{figure}
As in figure~\ref{figure-xi}, the plain lines $(\stackrel{j}{-\!\!\!-\!\!\!-})$ represent the path associated with N independent copies $X^j_{u,v}$ of the stochastic semigroup $X_{u,v}$ of the process $X_t$ with generator $L_t$. In this notation, we have the formula
\begin{equation}\label{YY-k-X-k}
\YY^k_t=\widehat{\zeta}^{k}_t=\left(X^k_{0,s}(\xi^k_0)\right)_{ s\leq t}\in D_t(S)
\end{equation}
  
 Observe that the generator $\GG_{\Ia,t}$ of the genealogical tree $\YY_t$ defined in (\ref{ref-GG-k-mn}) can be rewritten for any $F\in\Da(\widehat{\La})$ as follows
\begin{equation}\label{ref-GG-k-replacement}
\begin{array}{l}
   \displaystyle \GG_{\Ia,t}(F)(x) \\
   \\
     \displaystyle
   =\widehat{L}_t(F_{x^{-\Ia}})(x^\Ia)+\sum_{i\in [N]-\{\Ia\}}~\left[ \widehat{L}_{x^{\Ia},t}(F_{x^{-i}})(x^i)+\widehat{V}^-_t(x^i)~\int~(F_{x^{-i}}(z)-F(x))~m(x^{-\Ia})(dz) \right]
\end{array}
\end{equation}
In the above display, $  \widehat{L}_{z,t}$   stands for the collection of operators indexed by $z\in D_t(S)$ and defined for any function $f\in \Da(\widehat{L})$ and any $y\in D_t(S)$ by
  $$
  \widehat{L}_{z,t}(f)(y):=  \widehat{L}_{t}(f)(y)+\frac{2}{N}~\widehat{V}_t(y)~(f(z)-f(y))
  $$
  In the same vein, replacing 
 $(\widehat{L}_t,\widehat{V}_t)$ by  $(L_t,V_t)$ in (\ref{ref-GG-k-replacement}), for any function $F\in \Da(\La) $ and any $x=(x^i)_{1\leq i\leq N}\in S^N$
 we obtain the formula
 $$
\begin{array}{l}
   \displaystyle \Ga_{\Ia,t}(F)(x) \\
   \\
     \displaystyle
   =L_t(F_{x^{-\Ia}})(x^{\Ia})+\sum_{i\in [N]-\{\Ia\}}~\left[ L_{x^{\Ia},t}(F_{x^{-i}})(x^i)+V^-_t(x^i)~\int~(F_{x^{-i}}(z)-F(x))~m(x^{-\Ia})(dz) \right]
\end{array}
$$
the collection of operators $  L_{z,t}(f)(y)$ indexed by $z\in S$ and defined for any function $f\in \Da(L)$ and any $y\in S$ by
$$
  L_{z,t}(f)(y):=  L_{t}(f)(y)+\frac{2}{N}~V_t(y)~(f(z)-f(y))
\quad \mbox{\rm and}\quad   V^-_t:=\left(1-\frac{1}{N}\right)~V_t
  $$
  
  \begin{defi}
  For any final time horizon $T$ and given the ancestral line  $ \YY^{\Ia}_T=\widehat{\zeta}^{\,\Ia}_T\in D_T(S)$, 
 we denote by $\QQ^-_t=\widehat{\eta}^{\,-}_{t}$ and $\eta^-_{t}$  the conditional Feynman-Kac measures  defined as in (\ref{red-QQ-nonlinear}) and (\ref{ref-FK})  by replacing $(\widehat{L}_t,\widehat{V}_t)$ by $(\widehat{L}^-_{
t},\widehat{V}^-_t)$, and respectively   $(L_t,V_t)$ by $(L_{t}^-,V^-_t)$, with the quenched generators
\begin{equation}\label{def-L-V}
\widehat{L}^-_{t}:=
\widehat{L}_{\widehat{\zeta}^{\,\Ia}_t,t}\quad\mbox{and}\quad L^-_{t}:=L_{\zeta^{\,\Ia}_t,t}
\end{equation}
  \end{defi}
Using (\ref{ref-YY-w-zeta}) we readily check the following theorem.
     \begin{theo}\label{theo-frozen-FK}
     
 For any time horizon $T$ the $\Ia$-ancestral line  $ \YY^{\Ia}_T\in D_T(S)$ has the same law as the historical process $\widehat{X}_t$.
 
 In addition, given the ancestral line  $ \YY^{\Ia}_T$, the 
 processes defined by
$$
t\in [0,T]~\mapsto~ \YY^{-}_t=\left(\YY_t^{i}\right)_{i\in [N]-\{\Ia\}}\in D_t(S)^{N-1}\quad \mbox{and}\quad
  \zeta^{-}_t=\left(\zeta_t^{i}\right)_{i\in [N]-\{\Ia\}}\in S^{N-1}
$$
 coincides with the conditional $(N-1)$-mean field particle interpretation of the conditional Feynman-Kac measures $\QQ^-_t$ and $\eta^-_{t}$. 
 \end{theo}
  
   \subsection{Coupled historical and genealogical trees}

\begin{defi}
Let $\widehat{\Xi}_{s,t}$ be the transition semigroup of the flow $\left(\widehat{\xi}_{s,t}\left(x\right),\XX_{s,t}(y)\right)$ defined for any $F\in\Ba\left(D_t(S)^{(2N)}\right)$ and $(x,y)\in D_s(S)^{(2N)}$ by the formula
\begin{equation}\label{def-w-Pi}
\widehat{\Xi}_{s,t}(F)(x,y):=\EE\left(F\left(\widehat{\xi}_{s,t}\left(x\right),\XX_{s,t}(y)\right)\right)\end{equation}
\end{defi}

Observe that the processes $\widehat{\xi}_{s,t}\left(x\right)$ and $\XX_{s,t}(x)$ have the same terminal states $\xi_{s,t}^i\left(x^i_s\right)$.
  
  Using (\ref{ref-wG-m}) and (\ref{ref-GG-GG-k}) we check that  the 
generator $\widehat{\Ha}_{t}$ of the coupled process $(\widehat{\xi}_t,\XX_t)$   is defined  by the formula
  \begin{equation}\label{ref-wH-m}
 \widehat{\Ha}_{t}(F)(x,y)=\widehat{\La}_t^{(2)}(F)(x,y)+\widehat{\Ja}^{(2)}_{t}(F)(x,y)-\widehat{\lambda}_t(x)~F(x,y)
\end{equation}
with the  function $\widehat{\lambda}_t$ defined in (\ref{def-lambda-intensity}) and the jump-intensity  integral operator $\widehat{J}^{(2)}_{t}$ given by
  \begin{equation}\label{def-wJ-xy}
\widehat{\Ja}^{(2)}_{t}(F)(x,y):=\sum_{\iota\in [ N]^2_0}~\widehat{\lambda}_t^{\iota}(x)~F(\CC_{e,c_{\iota}}(x),\CC_{c_{\iota}}(y))\quad \mbox{\rm with}\quad\widehat{\lambda}_t^{i,j}(x):= \frac{1}{N}~\widehat{V}_t(x^i)
\end{equation}
In the above display,  $c_{i,j}$ and $\CC_{a,b}$ stands for the coalescent maps and operators  defined in (\ref{def-sigma-c}) and (\ref{def-varsigma-ab}).
\begin{defi}
Let $\widehat{\Xi}^{\,\Ia}_{s,t}$ be the conditional transition semigroup of the flow $\left(\widehat{\goodchi}_{s,t}\left(x\right),\YY_{s,t}(y)\right)$ defined for any  $F\in \Ba(D_t(S)^{(2N)})$ and $(x,y)\in D_t(S)^{(2N)}$ by the formula
\begin{equation}\label{def-w-Pi-Ia}
\widehat{\Xi}^{\,\Ia}_{s,t}(F)(x,y):=\EE\left(F\left(\widehat{\goodchi}_{s,t}\left(x\right),\YY_{s,t}(y)\right)~|~\Ia\right)\end{equation}
\end{defi}

Using (\ref{ref-wG-m-k}) and (\ref{ref-GG-k-mn}) we check that  the 
generator $\widehat{\Ha}_{\Ia,t}$ of the conditional coupled process $(\widehat{\goodchi}_t,\YY_t)$ given $\Ia$  is defined for any
$F\in \Da(\widehat{\La}^{(2)})$ and $(x,y)\in D_t(S)^{(2N)}$ by the formula
  \begin{equation}\label{ref-wH-m-k}
   \displaystyle \widehat{\Ha}_{\Ia,t}(F)(x,y)=\widehat{\La}_t^{(2)}(F)(x,y)+\widehat{\Ja}^{(2)}_{\Ia,t}(F)(x,y)-\widehat{\lambda}_{\Ia,t}(x)~F(x,y)
\end{equation}
  with the rate function $\widehat{\lambda}_{k,t}$ defined in (\ref{def-lambda-Ia-k}) and the jump-intensity  integral operator $\widehat{\Ja}^{(2)}_{k,t}$ given by
$$
\begin{array}{l}
\displaystyle\widehat{\Ja}^{(2)}_{\Ia,t}(F)(x,y)=\sum_{i\in [N]-\{\Ia\}} \sum_{j\in [N]-\{i\}}\widehat{\lambda}_{\Ia,t}^{i,j}(x)~F(\CC_{e,c_{i,j}}(x),\CC_{c_{i,j}}(y)) \\
\\
\displaystyle\hskip3cm+\sum_{j\in [N]-\{\Ia\}}~\widehat{\lambda}_{\Ia,t}^{\Ia,j}(x)~F(\CC_{\sigma_{\Ia,j},c_{j,\Ia}}(x),\CC_{c_{j,\Ia}}(y))\quad \mbox{\rm with}\quad\widehat{\lambda}_{k,t}^{i,j}:=1_{i\not=k}~\widehat{\lambda}_t^{i,j}+1_{i=k}~\widehat{\lambda}_t^{j,k}
\end{array}
$$
In the above display,  $c_{i,j}$,  $ \sigma_{i,j}$ and $\CC_{a,b}$ stands for the coalescent maps and operators 
defined in (\ref{def-sigma-c}) and (\ref{def-varsigma-ab}). 
Following word-for-word the proof of lemma~\ref{lem-connect-generators}, we check the following lemma.
\begin{lem}\label{lem-connect-generators-2}
The jump intensities $(\widehat{\Ja}^{(2)}_{t},\widehat{\Ja}^{(2)}_{k,t})$ and the generators $(\widehat{\Ha}_{t},\widehat{\Ha}_{k,t})$   defined in (\ref{ref-wH-m}) and  (\ref{ref-wH-m-k}) 
are connected  for any
$F\in \Da(\widehat{\La}^{(2)})$ and $(x,y)\in D_t(S)^{(2N)}$  by the formulae
\begin{eqnarray}
\widehat{\Ja}^{(2)}_{t}(F)(x,y)&=&\widehat{\Ja}^{(2)}_{k,t}(F)(x,y)+\frac{1}{N}~\sum_{i\in [N]-\{k\}}~\epsilon^{\prime}_{k,i}(F)(x,y)
\label{ref-anti-sym}\\
\widehat{\Ha}_{t}(F)(x,y)-m(x)(\widehat{V}_t)~F(x,y)&=&\widehat{\Ha}_{k,t}(F)(x,y)-\widehat{V}_t(x^k)~F(x,y)+\frac{1}{N}~\sum_{i\in [N]-\{k\}}~\epsilon^{\prime}_{k,i}(F)(x,y)\nonumber
\end{eqnarray}
with the bounded integral operator $\epsilon^{\prime}_{i,k}$ defined by the anti-symmetric functions
\begin{equation}\label{def-epsilon-prime}
\epsilon^{\prime}_{k,i}(F)(x,y):=\widehat{V}_t(x^k)~F(x\vee x_t^{c_{k,i}},y^{c_{k,i}})-\widehat{V}_t(x^i)~F(x^{\sigma_{i,k}}\vee x^{c_{i,k}}_t,y^{c_{i,k}})=-\epsilon_{k,i}(F)(x^{\sigma_{i,k}},y^{\sigma_{i,k}})
\end{equation}
\end{lem}

Using (\ref{ref-c-sigma-prop})
we check that
$$
\begin{array}{l}
\displaystyle\widehat{\Ja}^{(2)}_{\Ia,t}(F)(x,y)=\sum_{i\in [N]-\{\Ia\}} \sum_{j\in [N]-\{i\}}\widehat{\lambda}_{\Ia,t}^{i,j}(x)~F(\CC_{e,c_{i,j}}(x),\CC_{c_{i,j}}(y)) \\
\\
\displaystyle\hskip3cm+\sum_{j\in [N]-\{\Ia\}}~\widehat{\lambda}_{\Ia,t}^{\Ia,j}(x)~F(\CC_{e,c_{\Ia,j}}\left( \CC_{\sigma_{\Ia,j}}(x)\right),\CC_{c_{\Ia,j}}\left(\CC_{\sigma_{\Ia,j}}(y)\right)
\end{array}
$$
Rewritten in a more compact form in terms of the collection of coalescent operators $\CC^{k}_{e,a}$ and $  \CC^{k}_{a}$ indexed by $k\in [N]$ and $a\in\mathfrak{C}$ introduced in (\ref{def-CC-k-1}), we obtain for any
$F\in \Da(\widehat{\La}^{(2)})$ and $(x,y)\in D_t(S)^{(2N)}$ the formula
\begin{equation}\label{ref-Ja-2-k}
\begin{array}{l}
\displaystyle\widehat{\Ja}^{(2)}_{\Ia,t}(F)(x,y)=\sum_{\iota\in [N]_0^2}\widehat{\lambda}_{\Ia,t}^{\iota}(x)~F(\CC^{\Ia}_{e,c_{\iota}}(x),\CC^{\Ia}_{c_{\iota}}(y)) 
\end{array}
\end{equation}

\subsection{Dyson-Phillips expansions}

 The main objective of this section is to express the stochastic flows $(\widehat{\xi}_{s,t},\XX_{s,t})$   in terms of the stochastic evolution flows discussed in (\ref{def-Xa-wXa}) and the jump intensity (\ref{def-wJ-xy}).

 For any given $s\geq 0$ and $x\in D_s(S)$  we denote by $T_n^{s,x}$  an increasing sequence of jump times with stochastic  intensity 
$\widehat{\lambda}_t(\widehat{\xi}_{s,t}(x))$. In the further development of this section,  $[s,t]_n$ stands for the Weyl chamber introduced in (\ref{weyl-chamber-def}). 
  
Observe that
\begin{equation}\label{Za-Ea-obs}
\Ea_{s,t}\left(\widehat{\Xa}_{s,t}(x)\right):=\exp{\left(- \int_s^t\widehat{\lambda}_{r}\left(\widehat{\Xa}_{s,r}\left(x\right)\right)~dr\right)}=\Za_{s,t}\left(\widehat{\Xa}_{s,t}\left(x\right)\right)^{N-1}
\end{equation}
with the exponential map $\Za_{s,t}$ introduced in (\ref{def-Za-xi-t}).
  \begin{defi}
For any $n\geq 0$ We denote by $\widehat{\Ua}^{(n)}_{s,t}$ and  $\widehat{\UU}^{(n)}_{s,t}$  the integral semigroups defined for any $F\in \Ba(D_t(S)^{(2N)})$ and $(x,y)\in D_s(S)^{(2N)}$  by the formulae
\begin{eqnarray}
\widehat{\Ua}^{(n)}_{s,t}(F)(x,y)&:=&\EE\left(\widehat{\UU}^{(n)}_{s,t}(F)(x,y)\right)\nonumber\\
\widehat{\UU}^{(n)}_{s,t}(F)(x,y)&:=&F\left(\widehat{\Xa}_{s,t}\left(x\right),\widehat{\Xa}_{s,t}(y)\right)~ \Za_{s,t}\left(\widehat{\Xa}_{s,t}\left(x\right)\right)^{n}\label{def-Q-2}
\end{eqnarray}
\end{defi}

Using (\ref{Za-Ea-obs}) we have
\begin{eqnarray*}
\widehat{\Ua}^{(N-1)}_{s,t}(F)(x,y)&=&\EE\left(F\left(\widehat{\xi}_{s,t}\left(x\right),\XX_{s,t}(y)\right)~1_{T^{s,x}_1>t}\right)\\
\widehat{\UU}^{(N-1)}_{s,t}(F)(x,y)&=&\EE\left(F\left(\widehat{\xi}_{s,t}\left(x\right),\XX_{s,t}(y)\right)~1_{T^{s,x}_1>t}~|~\widehat{\Xa}_{s,t}\right)
\end{eqnarray*}
In addition, the semigroups  $\left(\widehat{\Xi}_{s,t},\widehat{\Ua}^{(N-1)}_{s,t}\right)$ are connected by the Gelfand-Pettis weak-sense integration formulae
$$
\widehat{\Xi}_{s,t}=\widehat{\Ua}^{(N-1)}_{s,t}+\int_s^t \widehat{\Ua}^{(N-1)}_{s,r}\;\widehat{\Ja}_{r}^{(2)}~\widehat{\Xi}_{r,t}~~dr
$$
with the jump intensity integral operator $\widehat{\Ja}_{r}^{(2)}$ defined in (\ref{ref-Ja-2-k}). Iterating the above implicit formula, we obtain the following proposition.

\begin{prop}
For any $s\leq t$ and any $F\in \Ba(D_t(S)^{(2N)})$ and $(x,y)\in D_s(S)^{(2N)}$ we have the Dyson-Phillips expansion
 \begin{equation}\label{ref-t-decom}
\widehat{\Xi}_{s,t}(F)(x,y)=\sum_{n\geq 0}~\EE\left(F(\widehat{\xi}_{s,t}(x),\XX_{s,t}(y))~ 1_{T_n^{s,x}<t<T^{s,x}_{n+1}}\right)
 \end{equation}
with the iterated integrals
$$
 \begin{array}{l}
\displaystyle\EE\left(F(\widehat{\xi}_{s,t}(x),\XX_{s,t}(y))~ 1_{T_n^{s,x}<t<T^{s,x}_{n+1}}\right)\\
\\
\displaystyle= \int_{[s,t]_n} \left(\left(\widehat{\Ua}^{(N-1)}_{r_0,r_1}\widehat{\Ja}^{(2)}_{r_1}\right)\ldots \left(\widehat{\Ua}^{(N-1)}_{r_{n-1},r_{n}}\widehat{\Ja}^{(2)}_{r_n}\right)\left(\widehat{\Ua}^{(N-1)}_{r_n,t}(F)\right)\right)(x)~dr
\end{array}
$$

\end{prop}

 We end this section with a probabilistic interpretation of (\ref{ref-t-decom}) in terms of the   collection of maps $\widehat{\Ta}_{s,t}^{\,a} $ and $\widehat{\Xa}_{s,t}^{\,a}$, indexed by $a\in \mathfrak{C}$  and $s\leq t$ defined in (\ref{def-wTa}).
 
 By (\ref{ref-wH-m}) and (\ref{def-wJ-xy}), given the stochastic flows $\widehat{\Xa}_{s,t}$, the jump times $T_n^{s,x}=r_n$ and the randomly selected coalescent maps $a_n=c_{\iota_n}$ for some index $\iota_n\in [N]_0^2$ we have the stochastic evolution semigroup formulae
 $$
  \widehat{\xi}_{r_{n-1},r_n}= \widehat{\Ta}_{r_{n-1},r_n}^{\,c_{\iota_n}} \quad \mbox{\rm and}\quad
 \XX_{r_{n-1},r_n}= \widehat{\Xa}_{r_{n-1},r_n}^{\,c_{\iota_n}} $$
 In addition, whenever $r_n\leq t<r_{n+1}$ we have
 \begin{equation}\label{def-wT-wX}
   \widehat{\xi}_{s,t}= 
   \widehat{\Xa}_{r_{n},t}\circ \widehat{\Ta}^{\,\iota,(n)}_{r,r_{n}} \quad \mbox{\rm and}\quad
      \XX_{s,t}=\widehat{\Xa}_{r_{n},t}\circ  \widehat{\Xa}^{\,\iota,(n)}_{r,r_{n}} \end{equation}
 with the convention $r_0=s$, and the compositions
 \begin{equation}\label{def-wT-wX-2}
 \widehat{\Ta}^{\,\iota,(n)}_{r,r_{n}}:=
 \widehat{\Ta}_{r_{n-1},r_n}^{\,c_{\iota_n}} \circ\ldots\circ \widehat{\Ta}_{r_{0},r_1}^{\,c_{\iota_1}} 
 \quad \mbox{\rm and}\quad 
  \widehat{\Xa}^{\,\iota,(n)}_{r,r_{n}}:=\widehat{\Xa}_{r_{n-1},r_n}^{\,c_{\iota_n}} \circ\ldots\circ \widehat{\Xa}_{r_{0},r_1}^{\,c_{\iota_1}}
\end{equation}
The conditional probability $\mathfrak{p}_n(\widehat{\Xa}_{s,t},d(\iota,r))$ on $\left([N]_0^{(2n)}\times [s,t]_n\right)$
 of the $n$ first jump times and the selected coalescent indices is given by the formula 
  \begin{equation}\label{def-p-iota-r}
 \begin{array}{l}
 \displaystyle\mathfrak{p}^{s,x}_n(\widehat{\Xa}_{s,t},d(\iota,r))\\
 \\
 :=\PP\left((T^{s,x}_{1},\ldots,T^{s,x}_{n})\in d(r_1,\ldots,r_n),~~(a_1,\ldots,a_n)=(c_{\iota_1},\ldots, c_{\iota_n}),~~1_{T^{s,x}_{n+1}>t}~|~\widehat{\Xa}_{s,t}\right)\\
 \\
 =\displaystyle
  \Ea_{r_0,r_{1}}\left( \widehat{\Xa}_{r_0,r_{1}}(x)\right)~\widehat{\lambda}_{r_1}^{\iota_1}\left(\widehat{\Xa}_{r_0,r_{1}}(x)\right) ~\Ea_{r_1,r_{2}}\left(    \widehat{\Xa}_{r_{1},r_2}\left(\widehat{\Xa}_{r_{0},r_1}^{\,c_{\iota_1}} (x)\right)\right)~\widehat{\lambda}_{r_2}^{\iota_2}\left(\widehat{\Xa}_{r_{1},r_2}\left(\widehat{\Xa}_{r_{0},r_1}^{\,c_{\iota_1}} (x)\right)\right)\times\ldots  \times\\
  \\
  \displaystyle\times~\Ea_{r_{n-1},r_{n}}\left( \widehat{\Xa}_{r_{n-1},r_n}\left(\widehat{\Xa}^{\,\iota,(n-1)}_{r,r_{n-1}}(x)\right)\right)~\widehat{\lambda}_{r_n}^{\iota_n}\left( \widehat{\Xa}_{r_{n-1},r_n}\left(\widehat{\Xa}^{\,\iota,(n-1)}_{r,r_{n-1}}(x)\right)\right)~\Ea_{r_n,t}\left( \widehat{\Xa}_{r_{n},t} \left( \widehat{\Xa}^{\,\iota,(n)}_{r,r_{n}}(x)\right)\right)~dr
 \end{array}
\end{equation}

This construction yields a probabilistic interpretation of the weak-sense Dyson-Phillips expansion stated in
(\ref{ref-t-decom}).
\begin{theo}\label{ref-prop-sum-times}
For any $s\leq t$, any $n\geq 1$ and any $F\in \Ba(D_t(S)^{(2N)})$ and $(x,y)\in D_s(S)^{(2N)}$ we have the conditional almost sure Dyson-Phillips expansion formulae
$$
  \begin{array}{l}
\displaystyle  \EE\left(F\left(\widehat{\xi}_{s,t}(x),\XX_{s,t}(y)\right)~ 1_{T^{s,x}_n<t<T^{s,x}_{n+1}}~|~\widehat{\Xa}_{s,t}\right)\\
\\
\displaystyle= \int_{[s,t]_n} \left(\left(\widehat{\UU}^{(N-1)}_{r_0,r_1}\widehat{\Ja}^{(2)}_{r_1}\right)\ldots \left(\widehat{\UU}^{(N-1)}_{r_{n-1},r_{n}}\widehat{\Ja}^{(2)}_{r_n}\right)\left(\widehat{\UU}^{(N-1)}_{r_n,t}(F)\right)\right)(x)~dr\\
\\
\displaystyle =\sum_{\iota\in [N]_0^{(2n)}}~\int_{[s,t]_n}~F\left(\widehat{\Ta}^{\iota,n}_{r,t}(x),\widehat{\Xa}^{\iota,n}_{r,t}(y)
\right)~\mathfrak{p}^{s,x}_n(\widehat{\Xa}_{s,t},d(\iota,r))
\end{array}
$$
\end{theo}

\begin{cor}\label{cor-ref-app-proof}
For any $s\leq t$, and any $F\in \Ba(D_t(S)^{(2N)})$, $(x,y)\in D_s(S)^{(2N)}$ and $n\geq 1$ we have the almost sure formula 
$$
  \begin{array}{l}
  \EE\left(F(\widehat{\xi}_{s,t}(x),\XX_{s,t}(y))~\Za_{s,t}\left(\XX_{s,t}\left(x\right)\right)~ 1_{T^{s,x}_n<t<T^{s,x}_{n+1}}~|~\widehat{\Xa}_{s,t}\right)\\
\\
\displaystyle=\int_{[s,t]_n} \left(\left(\widehat{\UU}^{(N)}_{r_0,r_1}\widehat{\Ja}^{(2)}_{r_1}\right)\ldots \left(\widehat{\UU}^{(N)}_{r_{n-1},r_{n}}\widehat{\Ja}^{(2)}_{r_n}\right)\left(\widehat{\UU}^{(N)}_{r_n,t}(F)\right)\right)(x)~dr
\end{array}
$$
with the semigroup $\widehat{\UU}^{(N)}_{s,t}$ defined in (\ref{def-Q-2}), and the convention $r_0=s$.
\end{cor}
\proof
By theorem~\ref{ref-prop-sum-times} we have
$$
  \begin{array}{l}
  \EE\left(F\left(\widehat{\xi}_{s,t}(x),\XX_{s,t}(y)\right)~\Za_{s,t}\left(\XX_{s,t}\left(x\right)\right)~ 1_{T^{s,x}_n<t<T^{s,x}_{n+1}}~|~\widehat{\Xa}_{s,t}\right)\\
\\
\displaystyle=\sum_{\iota\in [N]_0^{(2n)}}~\int_{[s,t]_n}~F\left(\widehat{\Ta}^{\iota,n}_{r,t}(x),\widehat{\Xa}^{\iota,n}_{r,t}(y)
\right)~\mathfrak{q}^{s,x}_n(\widehat{\Xa}_{s,t},d(\iota,r))
\end{array}
$$
with the measure
$$
\mathfrak{q}^{s,x}_n(\widehat{\Xa}_{s,t},d(\iota,r))=\Za_{s,t}\left(\widehat{\Xa}^{\iota,n}_{r,t}(x)\right)~\mathfrak{p}^{s,x}_n(\widehat{\Xa}_{s,t},d(\iota,r))
$$
Also observe that
$$
 \begin{array}{l}
 \displaystyle
\Za_{s,t}\left(\widehat{\Xa}^{\iota,n}_{r,t}(x)\right)=
  \Za_{s,r_{1}}\left( \widehat{\Xa}_{s,r_{1}}(x)\right)~ ~\Za_{r_1,r_{2}}\left(    \widehat{\Xa}_{r_{1},r_2}\left(\widehat{\Xa}_{r_{0},r_1}^{\,c_{\iota_1}} (x)\right)\right)~ \\
  \\
\hskip3cm  \displaystyle\times\ldots  \times~\Za_{r_{n-1},r_{n}}\left( \widehat{\Xa}_{r_{n-1},r_n}\left(\widehat{\Xa}^{\,\iota,(n-1)}_{r,r_{n-1}}(x)\right)\right)~\Za_{r_n,t}\left( \widehat{\Xa}^{\iota,n}_{r,t}(x)\right)
 \end{array}$$
In addition, for any $y\in D_{r_n}(S)^N$ and $n\geq 1$ we have
$$
\Ea_{r_{n-1},r_{n}}\left( y\right)\times \Za_{r_{n-1},r_{n}}\left(y\right)=\Za_{r_{n-1},r_{n}}\left( y\right)^N
$$
The end of the proof of the corollary is now easily completed.
\cqfd

\section{Perturbation analysis}\label{sec-perturbation}
\subsection{Semigroup estimates}
We consider a collection of generators $L_t^{\epsilon}$ and potential functions $V_t^{\delta}$ of the form
$$
L_t^{\epsilon}=L_t+\epsilon~\overline{L}_t\quad \mbox{\rm and}\quad V_t^{\delta}=V_t+\delta\overline{V}_t\quad \mbox{\rm with}\quad \epsilon,\vert \delta\vert\in [0,1]
$$
In the above display, $\overline{V}_t$ stands for some uniformly bounded function and
$\overline{L}_t$ a bounded generator of an auxiliary jump type Markov process of the form
 $$
 \overline{L}_t(f)(x)=\lambda(x)~\int~(f(y)-f(x))~K_t(x,dy)
 $$
 for some jump rate function function $ \lambda(x)$ and some Markov transitions $K_t(x,dy)$ such that 
 $$
 \lambda_1\leq \lambda(x)\leq \lambda_2\quad\mbox{\rm and}\quad
 \varpi_1~\kappa_t(dy)\leq  K_t(x,dy)\leq  \varpi_2~\kappa_t(dy)
 $$
 In the above display, $\lambda_i, \varpi_i$ stands for some positive parameters and $\kappa_t$ some probability measures. 

 We let $P^{\epsilon}_{s,t}$ be the transition semigroup of the process with generator $L_t^{\epsilon}$. In this notation, we have the following technical lemma.
 
 \begin{lem}\label{lem-H0-perturbation}
 Assume that $P_{s,t}$ satisfies  $(H_0)$ for some parameters $h$ and $\rho(h)>0$ and some probability measures $\mu_{t,h}$. In this situation, for any $\epsilon \in [0,1]$ and $t\geq 0$  there exists some probability measures  $\mu^{\epsilon}_{t,h}$ such that
 \begin{equation}\label{ref-termin}
 \rho_{\epsilon}(h)~\mu^{\epsilon}_{t,h}(dy)\leq P^{\epsilon}_{t,t+h}(x,dy)\leq  \rho_{\epsilon}(h)^{-1}~\mu^{\epsilon}_{t,h}(dy)
\end{equation}
with the parameters
  \begin{eqnarray*}
 \rho_{\epsilon}(h)&:=&\rho(h)\left(e^{-\epsilon\lambda_2h}+(1-e^{-\epsilon\lambda_2h})\varpi_2\right)\min{\left( (\lambda_1/\lambda_2)(\varpi_1/\varpi_2) , e^{-\epsilon(\lambda_2-\lambda_1)h}\right)}\\
 &\geq &\rho(h)~\min{\left( (\lambda_1/\lambda_2)(\varpi_1/\varpi_2) , e^{-(\lambda_2-\lambda_1)h}\right)}
  \end{eqnarray*}
 \end{lem}
The proof of the above lemma is provided in the appendix on page~\pageref{lem-H0-perturbation-proof}.

 We  consider the Feynman-Kac semigroup $Q^{\delta,\epsilon}_{s,t}$ defined as $Q_{s,t}$ by replacing $V_t$ by $V_t^{\delta}$ and
$X_t$ by a Markov process 
with generator $L_t^{\epsilon}$. 

Also let $\phi^{(\delta,\epsilon)}_{s,t}$ be defined as $\phi_{s,t}$ by replacing $Q_{s,t}$ by $Q^{\delta,\epsilon}_{s,t}$, 
and set 
$$
L^{\delta,\epsilon}_{t}=\epsilon \overline{L}_t-\delta~\overline{\Va}_t\quad\mbox{\rm and}\quad L^{\delta,\epsilon}_{t,\eta}=\epsilon \overline{L}_t-\delta(\overline{\Va}_t-\eta(\overline{\Va}_t))\quad \mbox{\rm with}\quad \overline{\Va}_t(f):=\overline{V}_t~f
$$

 \begin{theo}\label{theo-epsilon-delta}
 For any $\vert\epsilon\vert,\vert \delta\vert\in [0,1]$ and any $s\leq t$ we have the perturbation semigroup Gelfand-Pettis weak-sense integration formulae
   \begin{equation}\label{1st-pert}
Q_{s,t}^{\delta,\epsilon}-Q_{s,t}=\int_s^t~Q_{s,u}^{\delta,\epsilon}~L^{\delta,\epsilon}_u~Q_{u,t}~du=\int_s^t~Q_{s,u}~L^{\delta,\epsilon}_u~Q_{u,t}^{\delta,\epsilon}~du \end{equation}
 In addition, for any $\eta\in \Pa(S)$ we have
  $$
 \phi^{\delta,\epsilon}_{s,t}(\eta)-\phi_{s,t}(\eta)=\int_s^t\phi^{\delta,\epsilon}_{s,u}(\eta)~
L^{\delta,\epsilon}_{u,\phi^{\delta,\epsilon}_{s,u}(\eta)}~
 \partial_{\phi^{\delta,\epsilon}_{s,u}(\eta)}\phi_{u,t}~du =\int_s^t\phi_{s,u}(\eta)
~L^{\delta,\epsilon}_{u,\phi_{s,u}(\eta)}~
 \partial_{\phi_{s,u}(\eta)}\phi^{\delta,\epsilon}_{u,t}~du$$
 \end{theo}
 \proof
We check (\ref{1st-pert}) using the fact that
$$
\partial_u (Q_{s,u}^{\delta,\epsilon}Q_{u,t})=Q_{s,u}^{\delta,\epsilon}~(L^{\epsilon}_u-L_u-\delta~\overline{\Va}_u)~Q_{u,t}
=\epsilon~Q_{s,u}^{\delta,\epsilon}~\overline{L}_u~Q_{u,t}-\delta~Q_{s,u}^{\delta,\epsilon}~\overline{\Va}_u~Q_{u,t}
$$
and
$$
\partial_u (Q_{s,u}Q_{u,t}^{\delta,\epsilon})
=-\epsilon~Q_{s,u}~\overline{L}_u~Q_{u,t}^{\delta,\epsilon}+\delta~Q_{s,u}~\overline{\Va}_u~Q_{u,t}^{\delta,\epsilon}
$$
 The perturbation analysis of the normalized semigroups $\phi^{\delta,\epsilon}_{s,t}$ is slightly more involved.

Let $\Lambda^{\delta,\epsilon}_{t}$   be  defined as $\Lambda_{t}$ by replacing $(L_t,V_t)$ by $(L^{\epsilon}_t,V_t^{\delta})$. Notice that
 $$
h^{-1}\left[\phi^{\delta,\epsilon}_{t,t+h}(\eta)-\eta\right]=\Lambda^{\delta,\epsilon}_{t}(\eta)+\mbox{\rm O}(h)
 $$
 For any given $s\leq t$, we consider the interpolating maps $u\in [s,t]\mapsto \Delta^{\delta,\epsilon}_{s,u,t}$ defined by
 $$
 \Delta^{\delta,\epsilon}_{s,u,t}:=\phi_{u,t}\circ\phi^{\delta,\epsilon}_{s,u}
 $$
 On the other hand, for any $s\leq u\leq u+h\leq t$  we have the decomposition
 $$
 \begin{array}{l}
 \Delta^{\delta,\epsilon}_{s,u+h,t}(\eta)- \Delta^{\delta,\epsilon}_{s,u,t}(\eta)\\
 \\
 \displaystyle=\phi_{u+h,t}\left(\phi^{\delta,\epsilon}_{s,u+h}(\eta)\right)-\phi_{u,t}\left(\phi^{\delta,\epsilon}_{s,u+h}(\eta)\right)
 +
 \phi_{u,t}\left(\phi^{\delta,\epsilon}_{s,u+h}(\eta)\right)-\phi_{u,t}\left(\phi^{\delta,\epsilon}_{s,u}(\eta)\right)\\
 \\
\displaystyle =-\Lambda_{u}(\phi^{\delta,\epsilon}_{s,u}(\eta))\left(\partial_{\phi^{\delta,\epsilon}_{s,u}(\eta)}\phi_{u,t}~\right)h\\
 \\
\hskip3cm \displaystyle+
 \phi_{u,t}\left(\phi^{\delta,\epsilon}_{s,u}(\eta)+\left[\phi^{\delta,\epsilon}_{s,u+h}(\eta)-\phi^{\delta,\epsilon}_{s,u}(\eta)\right]\right)-\phi_{u,t}\left(\phi^{\delta,\epsilon}_{s,u}(\eta)\right)+\mbox{\rm O}(h^2)
 \end{array}$$
 This implies that
 $$
 \begin{array}{l}
 h^{-1}\left[ \Delta^{\delta,\epsilon}_{s,u+h,t}(\eta)- \Delta^{\delta,\epsilon}_{s,u,t}(\eta)\right]\\
 \\
 \displaystyle=-\Lambda_{u}(\phi^{\delta,\epsilon}_{s,u}(\eta))\partial_{\phi^{\delta,\epsilon}_{s,u}(\eta)}\phi_{u,t}+h^{-1}\left[\phi^{\delta,\epsilon}_{s,u+h}(\eta)-\phi^{\delta,\epsilon}_{s,u}(\eta)\right]\partial_{\phi^{\delta,\epsilon}_{s,u}(\eta)}\phi_{u,t}+\mbox{\rm O}(h)
 \end{array}$$
 We conclude that
 $$
 \partial_u \Delta^{\delta,\epsilon}_{s,u,t}(\eta)=\left[\Lambda^{\delta,\epsilon}_{u}(\phi^{\delta,\epsilon}_{s,u}(\eta))-\Lambda_{u}(
 \phi^{\delta,\epsilon}_{s,u}(\eta))\right]\partial_{\phi^{\delta,\epsilon}_{s,u}(\eta)}\phi_{u,t}
 $$
 On the other hand, we have
 $$
\left[\Lambda^{\delta,\epsilon}_{t}(\eta)-\Lambda_{t}(\eta)\right](f)=\epsilon~\eta(\overline{L}_t(f))-\delta~\eta(f~
(\overline{V}_t-\eta(\overline{V}_t)))
 $$
By symmetry arguments, this ends the proof of the theorem.\cqfd
\begin{cor}\label{cor-ref-f}
 For any $s\leq t$ and any $\eta\in \Pa(S)$ we have the estimates
 \begin{eqnarray*}
 (H_1)&\Longrightarrow&\Vert  \phi^{\delta,\epsilon}_{s,t}(\eta)-\phi_{s,t}(\eta)\Vert_{\tiny\rm tv}\leq c~(\epsilon+\delta)\\
  (H_2)&\Longrightarrow&\Vert  \phi^{\delta,\epsilon}_{s,t}(\eta)-\phi_{s,t}(\eta)\Vert_{\tiny\rm tv}\leq c~(\epsilon+\delta)~(t-s)
 \end{eqnarray*}
 for some finite constant $c$ whose value doesn't depend on the parameters $(s,t,\eta)$, nor on $(\epsilon,\delta)$.

\end{cor}
  
\subsection{Particle stochastic flows}\label{sec-particle-flows}
Given a random measure $\mu$ on some measurable state space $(E,\Ea)$ and some $\sigma$-field $\Fa\subset \Ea$ we write $\EE(\mu~|~\Fa)$ the first conditional moment measure given by
$$
\EE(\mu~|~\Fa)~:~f\in \Ba(S)\mapsto \EE(\mu~|~\Fa)(f):=\EE(\mu(f)~|~\Fa)
$$

For any  $t\geq 0$, we let $\Delta m(\xi_t)$ be the random jump occupation measure 
$$
\Delta m(\xi_t):=m(\xi_t)-m(\xi_{t-})=\Delta M_t=M_{t}-M_{t-}
$$
with the martingale random field $M_t$ defined in (\ref{def-Mt}).
In this notation, we have
\begin{equation}\label{ref-fn}
\displaystyle N^{n-1}~\partial_t\EE\left[(\Delta m(\xi_t))^{\otimes n}(f_t^{(1)}\otimes\ldots \otimes f^{(n)}_t)~|~\Fa_{t-}\right]=m(\xi_{t-})\Gamma^{(n)}_{L^d_{m(\xi_{t-})}}(f_t^{(1)},\ldots, f^{(n)}_t)
\end{equation}
with the operators $\Gamma^{(n)}_{L^d_{m(\xi_{t-})}}$ defined in (\ref{n-Gamma}). When $n=2$ the above formula resumes to
$$
\begin{array}{l}
\displaystyle\partial_t\EE\left[\Delta m(\xi_t)(f_t)~\Delta m(\xi_t)(g_t)~|~\Fa_{t-}\right]
\displaystyle=\frac{1}{N}~m(\xi_{t-})\Gamma_{L^d_{m(\xi_{t-})}}(f_t,g_t)\\
\\
\displaystyle=\partial_t\langle M^d_t(f),M^d(g)\rangle_t
=\partial_t\langle \Ma^d_t(F),\Ma^d(G)\rangle_t
\end{array}
$$
 with the functions  $(F,G)$ defined in (\ref{link-fF})

\begin{defi}

For any $t\geq s$ and $n\geq 1$, we consider the integral random operators
\begin{eqnarray*}
\Delta^{n} \phi_{s,t}(m(\xi_{s}))&:=&~N^{n-1}~\frac{1}{n!}~(\Delta m(\xi_{s}))^{\otimes n}~
\overline{\partial}^n_{m(\xi_{s-})+\Delta m(\xi_{s}),m(\xi_{s-})}\phi_{s,t}
\end{eqnarray*}
and their first variational measure
$$
\Upsilon^{n}_{m(\xi_{s-})} \phi_{s,t}:=\partial_s\EE\left[\Delta^{n} \phi_{s,t}(m(\xi_{s}))~|~\Fa_{s-}\right]
$$
\end{defi}

Choosing $n=1$ we have
\begin{eqnarray*}
\Delta\phi_{s,t}(m(\xi_{s}))&:=&\Delta^1\phi_{s,t}(m(\xi_{s}))=\phi_{s,t}(m(\xi_{s}))-\phi_{s,t}(m(\xi_{s-}))
\end{eqnarray*}
Arguing as in the proof of (\ref{ref-r}) and using (\ref{ref-fn}), for any collection of functions $f^{(n)}\in \mbox{\rm Osc}(S)$ we have
the estimate
\begin{equation}\label{ref-fn-2}
N^{n-1}~\partial_s\,\EE\left[\Delta\phi_{s,t}(m(\xi_{s}))^{\otimes n}(f_t^{(1)}\otimes\ldots \otimes f^{(n)}_t)~|~\Fa_{s-}\right]\leq e^{nq}~\Vert \lambda+V\Vert
\end{equation}

\begin{prop}
For any $t\geq s$ and $n\geq 1$, we have
\begin{equation}\label{Delta-Phi-1}
\Delta^{n} \phi_{s,t}(m(\xi_{s}))=\frac{N^{n-1}}{n!}~(\Delta m(\xi_{s}))^{\otimes n}{\partial}^n_{m(\xi_{s-})}\phi_{s,t}+\frac{1}{N}~\Delta^{n+1} \phi_{s,t}(m(\xi_{s}))
\end{equation}
In addition, for any $f\in \Ba(S)$ we have
\begin{equation}\label{Delta-Phi-2}
\begin{array}{l}
\displaystyle
\displaystyle\Upsilon^{n}_{m(\xi_{s-})} \phi_{s,t}(f)\\
\\
=
\displaystyle(-1)^{n-1}~m(\xi_{s-})\Gamma^{(n)}_{L^d_{m(\xi_{s-})}}\left(Q_{s,t}^{m(\xi_{s-1})}(1),\ldots,Q_{s,t}^{m(\xi_{s-1})}(1),\partial_{m(\xi_{s-})}\phi_{s,t}(f)\right)+\frac{1}{N}~\Upsilon^{n+1}_{m(\xi_{s-})} \phi_{s,t}(f)
\end{array}
\end{equation}
\end{prop}
\proof
We have
\begin{eqnarray*}
\Delta^{n+1} \phi_{s,t}(m(\xi_{s}))&=&N^{n}~\left[\Delta\phi_{s,t}(m(\xi_{s}))-\sum_{1\leq k\leq n}\frac{1}{k!}~(\Delta m(\xi_{s}))^{\otimes k}{\partial}^k_{m(\xi_{s-})}\phi_{s,t}\right]\\
&=&N~\Delta^{n} \phi_{s,t}(m(\xi_{s}))-\frac{N^n}{n!}~(\Delta m(\xi_{s}))^{\otimes n}{\partial}^n_{m(\xi_{s-})}\phi_{s,t}\quad\Longleftrightarrow \quad(\ref{Delta-Phi-1})
\end{eqnarray*}
This implies that
$$
\begin{array}{l}
\displaystyle
\partial_s\EE\left[\Delta^{n} \phi_{s,t}(m(\xi_{s}))~|~\Fa_{s-}\right]:=\Upsilon^{n}_{m(\xi_{s-})} \phi_{s,t}~
\\
\\
\displaystyle=\frac{N^{n-1}}{n!}~\partial_s\EE\left[(\Delta m(\xi_{s}))^{\otimes n}{\partial}^n_{m(\xi_{s-})}\phi_{s,t}~|~\Fa_{s-}\right]+\frac{1}{N}~\partial_s\EE\left[\Delta^{n+1} \phi_{s,t}(m(\xi_{s}))~|~\Fa_{s-}\right]
\end{array}$$
This ends the proof of the proposition.
\cqfd

\begin{lem}\label{formula-remainder-jumps-lem}
For any $n\geq 1$ and $s\leq t$ we have the almost sure uniform estimates
\begin{equation}\label{formula-remainder-jumps}
(H_2)\Longrightarrow  \Vert \Upsilon^{n}_{m(\xi_{s-})} \phi_{s,t}\Vert_{\tiny\rm tv}\leq 2^{n-1}e^{(n+1)q}~\Vert \lambda+V\Vert
\end{equation}
\end{lem}
The detailed proof of the above estimate is provided in the appendix, on page~\pageref{formula-remainder-jumps-proof}.

In the further development of this section, 
for any  given time horizon $t$ and any $f\in \Ba(S)$ we let
$$s\in [0,t]\mapsto\Ma^d_{s}(\phi_{\point,t}(m(\point))(f))$$ be the martingale
$s\in [0,t]\mapsto \Ma^d_{s}(F)$
associated with the function
$$
(s,x)\in [0,t]\times S^N\mapsto F(s,x)=\phi_{s,t}(m(x))(f)
$$
We also denote by $$s\in [0,t]\mapsto M^c_s\left(\partial_{m(\xi_{\point})}\phi_{\point,t}(f)\right)~, \quad\mbox{\rm resp.}\quad M^c_s\left(Q^{m(\xi_{\point})}_{\point,t}(1)\right)$$  the martingale $M^c_{s}(f)$
associated with the  $\Fa$-predictable bounded function
$$
(s,x)\in [0,t]\times S\mapsto f_s(x)=\partial_{m(\xi_{s-})}\phi_{s,t}(f)(x)~,\quad\mbox{\rm resp.}\quad  f_s(x)=Q^{m(\xi_{s-})}_{s,t}(1)(x)
$$
We are now in position to state and to prove the main result of this section.
\begin{theo}\label{theo-key-decom}
For any time horizon $t\geq 0$ and any $f\in\Ba(S)$ the interpolating function $$s\in [0,t]\mapsto \phi_{s,t}(m(\xi_s))(f)\in\RR$$ satisfies the stochastic differential equation
  \begin{eqnarray*}
d\phi_{s,t}(m(\xi_s))(f)
&=&\displaystyle\frac{1}{\sqrt{N}}~ dM^c_s\left(\partial_{m(\xi_{\point})}\phi_{\point,t}(f)\right)+ d\Ma^d_{s}\left(\phi_{\point,t}(m(\point))(f)\right)\\
&&\displaystyle+\frac{1}{N}~\Upsilon^{2}_{m(\xi_{s-})} \phi_{s,t}(f)~ds- \frac{1}{N}~m(\xi_s)\Gamma_{L^c_s}\left(Q^{m(\xi_{s})}_{s,t}(1),\partial_{m(\xi_s)}\phi_{s,t}(f)\right)~ds
 \end{eqnarray*}
\end{theo}
\proof
Observe that
 $$
 dm(\xi_s)=\Lambda_{s}(m(\xi_s))~ds+ \frac{1}{\sqrt{N}}~dM^c_s+\Delta m(\xi_s)-\underbrace{\EE(\Delta m(\xi_s)~|~\Fa_{s-})}_{\ll ds} $$
  Using It\^o formula and the backward formula (\ref{forward-backward-ref}) we have
  $$
 \begin{array}{l}
d~\phi_{s,t}(m(\xi_s))(f)
\displaystyle=-\Lambda_{s}(m(\xi_s))\left(\partial_{m(\xi_s)}\phi_{s,t}(f)\right)~ds+\left[\phi_{s,t}(m(\xi_{s-})+dm(\xi_s))-\phi_{s,t}(m(\xi_{s-}))\right](f)\\
\\
\hskip.5cm\displaystyle= \frac{1}{\sqrt{N}}~dM^c_s\left(\partial_{m(\xi_{\point})}\phi_{\point,t}(f)\right)+d\Ma^d_{s}(\phi_{\point,t}(m(\point))(f))\\
\\
\displaystyle\hskip1cm+\frac{1}{2N}~(dM_s^c\otimes dM^c_s)~\partial^2_{m(\xi_s)}\phi_{s,t}(f)+\partial_s\EE\left[\Delta \phi_{s,t}(m(\xi_{s}))(f)-\Delta m(\xi_{s})\partial_{m(\xi_{s-})}\phi_{s,t}(f)~|~\Fa_{s-}\right] ~ds\end{array}
 $$
This ends the proof of the theorem.
\cqfd

The above theorem can be interpreted as an Aleeksev-Gr\"obner formula for interpolating semigroups on the space of measures~\cite{dm-ssd-19}.

Next corollary is a direct consequence of the recursion (\ref{Delta-Phi-2}).

\begin{cor}
For any $t\geq 0$ and any $f\in\Ba(S)$ we have the almost sure formula
  \begin{equation}\label{cor-ups-3}
   \begin{array}{l}
   \phi_{s,t}(m(\xi_s))(f)-\phi_{0,t}(m(\xi_0))(f)\\
   \\
=\displaystyle\frac{1}{\sqrt{N}}~M^c_s\left(\partial_{m(\xi_{\point})}\phi_{\point,t}(f)\right)+
\Ma^d_{s}(\phi_{\point,t}(m(\point))(f))\\
\\
\hskip1cm\displaystyle- \frac{1}{N}~\int_0^s~m(\xi_u)\Gamma_{L_{u,m(\xi_u)}}\left(Q^{m(\xi_{u})}_{u,t}(1),\partial_{m(\xi_u)}\phi_{u,t}(f)\right)~du+\frac{1}{N^2}~\int_0^s~\Upsilon^{3}_{m(\xi_{u-})} \phi_{u,t}(f) du
 \end{array}
 \end{equation}
\end{cor}
Choosing $s=t$ and taking the expectation in (\ref{cor-ups-3}) we obtain the following result.

 \begin{cor}\label{cor-EE-m}
 For any $t\geq 0$ and $f\in \Da(S)$ we have the formula
$$
 \begin{array}{l}
\EE(m(\xi_t)(f))-\EE(\phi_{0,t}(m(\xi_0))(f))\\
\\
\displaystyle=- \frac{1}{N}~\int_0^t~\EE\left[m(\xi_s)\Gamma_{L_{s,m(\xi_s)}}
\left(Q^{m(\xi_{s})}_{s,t}(1),\partial_{m(\xi_{s})}\phi_{s,t}(f)
\right)\right]~ds+\frac{1}{N^2}\int_0^t~\EE\left[\Upsilon^{3}_{m(\xi_{s-})} \phi_{s,t}(f)\right] ~ds
 \end{array}
 $$
 \end{cor}

\subsection{Some non-asymptotic estimates}\label{sec-non-asympt}
 
\begin{theo}\label{cor-ref-zeta-m}
 For any time horizon $t\geq 0$ and any function $f\in \mbox{\rm Osc}(S)$  we have
 \begin{eqnarray}
 (H_1)&\Longrightarrow&\vert\EE(m(\xi_t)(f))-\eta_t(f)\vert\leq c/N\nonumber\\
  (H_2)&\Longrightarrow&\vert\EE(m(\xi_t)(f))-\eta_t(f)\vert\leq c~t/N\label{ref-bias-path-H1}
 \end{eqnarray}
 In addition,  for any function $F\in \mbox{\rm Osc}(D_t(S))$   we have
 \begin{equation}\label{ref-bias-path-space}
   (H_1)\Longrightarrow\vert\EE(m(\XX_t)(F))-\QQ_t(F)\vert\leq c~t/N
 \end{equation}
 for some finite constant $c$ whose value doesn't depend on the parameters $(t,N)$.

\end{theo}
The assertion (\ref{ref-bias-path-space}) is a direct consequence of (\ref{ref-H-wH-2}) and (\ref{ref-bias-path-H1}).
The proof of the estimates (\ref{ref-bias-path-H1}) is mainly based on the decomposition presented in corollary~\ref{cor-EE-m}.
The estimates rely on elementary but rather technical  carr\'e du champ inequalities, and semigroup techniques. Thus, the detail of the proof  is housed in the appendix, on page~\pageref{cor-ref-zeta-m-proof}.
 
 The first estimate stated in the above corollary extends the bias estimate obtained in~\cite{mathias} to time varying Feynman-Kac models. The central difference between homogeneous and time varying models lies on the fact that we cannot use $h$-process techniques. The latter allows to interpret the 
 Feynman-Kac semigroups in terms of more conventional  Markov semigroups.
 
 We end this section with a some more or less direct consequences of the above estimates in the analysis of the measures discussed in theorem~\ref{theo-frozen-FK}.
 
 By corollary~\ref{cor-ref-f},  for any $N>1$ we have
$$
 (H_1)\Longrightarrow\Vert  \eta^-_t-\eta_t\Vert_{\tiny\rm tv}\leq c/N\quad\mbox{\rm and}\quad
  (H_2)\Longrightarrow\Vert  \eta^-_t-\eta_t\Vert_{\tiny\rm tv}\leq c~t/N
$$
with the Feynman-Kac measures $\eta^-_t$ defined in (\ref{def-L-V}).

Arguing as in the proof of (\ref{ref-bias-path-space}) we also have the almost sure estimate
$$
 (H_1)\Longrightarrow\Vert  \QQ^-_t-\QQ_t\Vert_{\tiny\rm tv}\leq c~t/N
$$
with the Feynman-Kac measures $\QQ^-_t$ defined in (\ref{def-L-V}).

By (\ref{ref-termin}),
when $(H_0)$ is satisfied, the Feynman-Kac model defined in terms of  the pair
$(L_{s}^-,V^-_s)$ given in (\ref{def-L-V}) satisfies the stability property
$(H_1)$. Thus, using theorem~\ref{cor-ref-zeta-m} we readily deduce the following estimates.
\begin{cor}\label{cor-ref-gibbs-2}
Assume that $ (H_1)$ is met. In this situation,
for any $f\in \mbox{\rm Osc}(S)$ and $F\in \mbox{\rm Osc}(D_t(S))$ we have almost sure and uniform estimates given by
$$
 \vert\EE\left(m(\zeta^{-}_t)(f)~|~\widehat{\zeta}^{\,\Ia}_t\right)-\eta_t(f)\vert\leq c /N
\quad \mbox{and}\quad
 \vert\EE\left(m(\YY^{-}_t)(F)~|~\YY^{\,\Ia}_t\right)-\QQ_t(F)\vert\leq c~t/N
$$ 
\end{cor}
 
 The above results give some information on the bias of the occupation measures.
 We end this section with some propagation of chaos estimate. Using (\ref{cor-ups-3}), for any functions $f_i\in \mbox{\rm Osc}(S)$ we have
 $$
    \begin{array}{l}
\EE\left(m(\xi_t)(f_1)~m(\xi_t)(f_2)\right)-\EE\left(\phi_{0,t}(m(\xi_0))(f_1)~\phi_{0,t}(m(\xi_0))(f_2)\right)\\\
   \\
=-\displaystyle\frac{1}{N}\sum_{(k,l)\in \{(1,2),(2,1)\}}\int_0^t\EE\left[ \phi_{s,t}(m(\xi_s))(f_k)~m(\xi_u)\Gamma_{L_{u,m(\xi_u)}}\left(Q^{m(\xi_{u})}_{u,t}(1),\partial_{m(\xi_u)}\phi_{u,t}(f_l)\right)~\right]~du\\
\\
\displaystyle\hskip1cm+\frac{1}{N}\int_0^t\EE\left[ m(\xi_u)\Gamma_{L_{u,m(\xi_u)}}\left(\partial_{m(\xi_u)}\phi_{u,t}(f_1),\partial_{m(\xi_u)}\phi_{u,t}(f_2)\right)~\right]\\
\\
\displaystyle\hskip2cm+\int_0^t~ \partial_s\,\EE\left[ \Delta\phi_{s,t}(m(\xi_{s}))(f_1)~ \Delta\phi_{s,t}(m(\xi_{s}))(f_2)\right]~ds\\
\\
\displaystyle\hskip3cm+\frac{1}{N^2}~\sum_{(k,l)\in \{(1,2),(2,1)\}}~\int_0^t~ \EE\left[\phi_{s,t}(m(\xi_s))(f_k)~\Upsilon^{3}_{m(\xi_{s-})} \phi_{s,t}(f_l) \right]~ds
 \end{array}
 $$
By (\ref{ref-fn-2}) and using the same lines of arguments as in the proof of theorem~\ref{cor-ref-zeta-m}  we check the following estimates.
\begin{cor}\label{cor-pchaos}
For any time horizon $t\geq 0$, any $f,g\in \mbox{\rm Osc}(S)$ and $i\not=j$  we have
 \begin{eqnarray*}
 (H_1)&\Longrightarrow&\vert\EE\left(f(\xi^i_t)~g(\xi^j_t)\right)-\eta_t(f)~\eta_t(g)\vert\leq c/N\\
  (H_2)&\Longrightarrow&\vert\EE\left(f(\xi^i_t)~g(\xi^j_t)\right)-\eta_t(f)~\eta_t(g)\vert\leq c~t/N
 \end{eqnarray*}
In addition, when $(H_1)$ is met, for any $F,G\in \mbox{\rm Osc}(D_t(S))$ we have
 $$
\vert\EE\left(F(\XX^i_t)~G(\XX^j_t)\right)-\QQ_t(F)~\QQ_t(G)\vert\leq c~t/N
$$
as well as the almost sure estimates
$$ \vert\EE\left(F(\YY^i_t)~G(\YY^j_t)~|~\YY^{\,\Ia}_t\right)-\QQ_t(F)~\QQ_t(G)\vert\leq c~t/N
 $$
 \end{cor}
 We can extend the above  arguments to any finite block of particles.

 \section{A duality formula}\label{duality-sec}
  \subsection{Many-body Feynman-Kac semigroups}\label{sec-many-body-FK}

 For any given $s\geq 0$ and $x\in D_s(S)$  we denote by $T_{k,n}^{s,x}$  an increasing sequence of jump times with intensity  $\widehat{\lambda}_{k,t}(\widehat{\goodchi}_{s,t}(x))$, on the time horizon $t\in [s,\infty[$. We use the convention $T_{k,0}^{s,x}=0$ for $n=0$.

\begin{defi}
Let $\widehat{\Pi}_{s,t}$ be the many-body Feynman-Kac semigroup defined for any $F\in \Ba(D_t(S)^{(2N)})$ and $(x,y)\in D_t(S)^{(2N)}$ by the formula
\begin{equation}\label{def-w-Pi}
\widehat{\Pi}_{s,t}(F)(x,y):=\EE\left(F\left(\widehat{\xi}_{s,t}\left(x\right),\XX_{s,t}(y)\right)~\Za_{s,t}\left(\XX_{s,t}\left(x\right)\right)\right)\end{equation}
\end{defi}

We are now in position to state and prove the main result of this section.

\begin{theo}\label{theo-duality-in}
For any $s\leq t$, and any $F\in \Ba(D_t(S)^{(2N)})$, any $(x,y)\in D_s(S)^{(2N)}$ and $n\geq 1$ we have
the almost sure formula
\begin{equation}\label{theo-a-e-ref}
\begin{array}{l}
\displaystyle  \EE\left(F(\widehat{\xi}_{s,t}(x),\XX_{s,t}(y))~\Za_{s,t}\left(\XX_{s,t}\left(x\right)\right)~ 1_{T^{s,x}_n<t<T^{s,x}_{n+1}}\right)\\
\\
\displaystyle=  \EE\left(F(\widehat{\goodchi}_{s,t}(x),\YY_{s,t}(y))~Z_{s,t}\left(\YY^{\Ia}_{s,t}\left(x\right)\right)~ 1_{T^{s,x}_{\Ia,n}<t<T^{s,x}_{\Ia,n+1}}~|~\Ia\right)
\end{array}
\end{equation}
with the exponential map ${Z}_{s,t}$ on $D_t(S)$ defined in
(\ref{def-Zt}).  In addition, we have the conditional transition semigroup formula
$$
\widehat{\Pi}_{s,t}(F)(x,y)=\EE\left(F(\widehat{\goodchi}_{s,t}(x),\YY_{s,t}(y))~Z_{s,t}\left(\YY^{\Ia}_{s,t}(x)\right)~|~ \Ia\right)
$$

\end{theo}

\proof
By lemma~\ref{lem-connect-generators-2} we have
$$
\begin{array}{l}
\left(\widehat{\Ua}^{(N)}_{s,t}\widehat{\Ja}^{(2)}_{t}(F)\right)(x,y)\\
\\
\displaystyle=\left(\widehat{\Ua}^{(N)}_{s,t}\widehat{\Ja}^{(2)}_{k,t}(F)\right)(x,y)+
\frac{1}{N}~\sum_{i\in [N]-\{k\}}~\EE\left(\epsilon^{\prime}_{k,i}(F)\left(\widehat{\Xa}_{s,t}\left(x\right),\widehat{\Xa}_{s,t}(y)\right)~ \Za_{s,t}\left(\widehat{\Xa}_{s,t}\left(x\right)\right)^{N} \right)
\end{array}
$$
with the bounded integral operators $\epsilon^{\prime}_{i,k}$ defined in (\ref{def-epsilon-prime}).
Replacing in the r.h.s. expectation  $\widehat{\Xa}_{s,t}$ by $\widehat{\Xa}_{s,t}^{\,\sigma_{i,j}}$ and using the fact that
$$
\Za_{s,t}\left(\widehat{\Xa}^{\,\sigma_{i,k}}_{s,t}\left(x\right)\right)=\Za_{s,t}\left(\widehat{\Xa}_{s,t}\left(x\right)\right)
\quad\mbox{\rm
we conclude that}\quad 
\widehat{\Ua}^{(N)}_{s,t}\widehat{\Ja}^{(2)}_{k,t}=\widehat{\Ua}^{(N)}_{s,t}\widehat{\Ja}^{(2)}_{t}.
$$
This implies that
for any $n\geq 1$ and $k\in [N]$ we have 
\begin{equation}\label{inter-prop}
  \begin{array}{l}
  \EE\left(F(\widehat{\xi}_{s,t}(x),\XX_{s,t}(y))~\Za_{s,t}\left(\widehat{\xi}_{s,t}\left(x\right)\right)~ 1_{T^{s,x}_n<t<T^{s,x}_{n+1}}\right)\\
  \\
  \displaystyle=\int_{[s,t]_n} \left(\left(\widehat{\Ua}^{(N)}_{r_0,r_1}\widehat{\Ja}^{(2)}_{k,r_1}\right)\ldots \left(\widehat{\Ua}^{(N)}_{r_{n-1},r_{n}}\widehat{\Ja}^{(2)}_{k,r_n}\right)\left(\widehat{\Ua}^{(N)}_{r_n,t}(F)\right)\right)(x)~dr
  \end{array}
\end{equation}
In the above display, $\widehat{\Ua}^{(N)}_{s,t}$ stands for
 the Feynman-Kac semigroup  defined in (\ref{def-Q-2})  and  $[s,t]_n$  the Weyl chamber introduced in (\ref{weyl-chamber-def}). The end of proof of (\ref{theo-a-e-ref}) is rather technical but it follows the same arguments as the ones used in the proof of theorem~\ref{ref-prop-sum-times}, thus it is provided in the appendix, on page~\pageref{theo-a-e-ref-proof}.

The second assertion is obtained by summing over the number of jumps in the interval $[s,t]$.
Alternatively, in terms of generators, using (\ref{ref-anti-sym}) and recalling that
$$
\widehat{V}_t(\XX_{s,t}^{\,k}\left(x\right))=V_t(\xi_{s,t}^{\,k}\left(x_s\right))=\widehat{V}_t(\widehat{\xi}_{s,t}^{\,k}\left(x\right))
$$
we check that
$$
 \begin{array}{l}
   \displaystyle
\partial_t\widehat{\Pi}_{s,t}(F)(x,y)\\
\\
   \displaystyle=\EE\left(\left(\widehat{\Ha}_{k,t}(F)\left(\widehat{\xi}_{s,t}\left(x\right),\XX_{s,t}(y)\right)-\widehat{V}_t(\widehat{\xi}_{s,t}^{\,k}\left(x\right))~F\left(\widehat{\xi}_{s,t}\left(x\right),\XX_{s,t}(y)\right)\right)~\Za_{s,t}\left(\widehat{\xi}_{s,t}\left(x\right)\right)\right)
\end{array}
$$
for any $k\in[N]$. This implies that 
$$
\partial_t\widehat{\Pi}_{s,t}(F)=\widehat{\Pi}_{s,t}\left(\widehat{\Ha}_{k,t}^V(F)\right)
\quad\mbox{\rm 
with}\quad
\widehat{\Ha}_{k,t}^V(F)(x,y):=\widehat{\Ha}_{k,t}(F)\left(x,y\right)-\widehat{V}_t(x^k)~F\left(x,y\right)
$$
Consider the semigroup
\begin{eqnarray*}
\widehat{\Pi}^{(\Ia)}_{s,t}(F)(x,y)&:=&\EE\left(F(\widehat{\goodchi}_{s,t}(x),\YY_{s,t}(y))~Z_{s,t}\left(\YY^{\Ia}_{s,t}(x)\right)~|~ \Ia\right)
\end{eqnarray*}
Arguing as above and recalling that $$
\widehat{V}_t(\widehat{\goodchi}^k_{s,t}(x))=V_t\left(\zeta^k_{s,t}\left(x_s\right)\right)=\widehat{V}_t\left(\YY^k_{s,t}(x)\right)
\quad\mbox{\rm we check that}
\quad
\partial_t\,\widehat{\Pi}^{(\Ia)}_{s,t}(F)=\widehat{\Pi}^{(\Ia)}_{s,t}\left(\widehat{\Ha}_{\Ia,t}^V(F)\right).
$$
This ends the proof of the theorem.\cqfd

Next corollary is a direct consequence of the above theorem and it extends the duality formula presented in~\cite{DMKP:16} to continuous time Feynman-Kac models.

\begin{cor}\label{cor-duality}
For any $t\geq 0$ and $F\in \Ba(D_t(S)^{(2N)})$ we have the almost sure formula
$$
\widehat{\Pi}_t(F):=\EE\left(F\left(\widehat{\xi}_{t},\XX_{t}\right)~\overline{\Za}_{t}\left(\XX_{t}\right)\right)=
\EE\left(F(\widehat{\goodchi}_{t},\YY_{t})~\overline{Z}_{t}\left(\YY^{\Ia}_t\right)\,|~\Ia~\right)
$$
with the function $\overline{\Za}_{t}$ defined as $\Za_{t}$
by replacing $V_t$ by the normalized potential $\overline{V}_t$ defined in (\ref{def-over-V}).
This yields for any $F\in \Ba(D_t(S)^{N+1})$ the duality formula
\begin{equation}\label{def-Pi-in}
\Pi_t(F):=\EE\left(F\left(\XX_{t},\XX^{\Ia}_t\right)~\overline{\Za}_{t}\left(\XX_{t}\right)\right)=
\EE\left(F(\YY_{t},\YY^{\Ia})~\overline{Z}_{t}\left(\YY^{\Ia}_t\right)\right)
\end{equation}
\end{cor}

 \subsection{Particle Gibbs samplers}\label{pg-section}
 For any given time horizon $t\geq 0$, the probability measure introduced in (\ref{def-Pi-in}) takes the form
 $$
 \Pi_t(d(z_1,z_2))\in \Pa(E_1\times E_2)\quad \mbox{\rm with}\quad E_1=D_t(S)^N\quad \mbox{\rm and}\quad
E_2:= D_t(S)
 $$
The first marginal of the measure $\Pi_t$ is given for any  $F\in \Ba(D_t(S)^{N})$  by the formula
$$
\pi_t(F):=\EE\left(F\left(\XX_{t}\right)~\overline{\Za}_{t}\left(\XX_{t}\right)\right)
$$

On the other hand, by corollary~\ref{cor-duality}, for any $k\in [N]$ and $F\in \Ba(\Da_t(S))$ we have
 $$
 \EE\left(F\left(\XX^k_t\right)~\overline{\Za}_{t}\left(\XX_{t}\right)\right)=
\EE\left(F(\YY^{k}_t)~\overline{Z}_{t}\left(\YY^{k}_t\right)\,|~\Ia=k~\right)=\QQ_t(F)
 $$
 In addition, for any $k,l\in [N]$ we have
 $$
 \QQ_t(F)=\EE\left(F\left(\XX^k_t\right)~\overline{\Za}_{t}\left(\XX_{t}\right)\right)=
\EE\left(F(\YY^{k}_t)~\overline{Z}_{t}\left(\YY^{l}_t\right)\,|~\Ia=l~\right)
 $$
 Thus, under $\widehat{\Pi}_t$ all the ancestral lines in $\YY_t$ are distributed according to $\QQ_t$.
 In addition, the second marginal of the measure $\Pi_t$ is given for any  $F\in \Ba(D_t(S))$ by the formula
$$
\QQ_t(F)=\EE\left(F(\YY^{\Ia})~\overline{Z}_{t}\left(\YY^{\Ia}_t\right)\right)=\EE\left(F\left(\XX^{\Ia}_t\right)~\overline{\Za}_{t}\left(\XX_{t}\right)\right)=\EE\left(m\left(\XX_t\right)(F)~\overline{\Za}_{t}\left(\XX_{t}\right)\right)
$$
This yields the duality formula stated in (\ref{intro-equation-duality-ref}) and in corollary~\ref{cor-duality}.\label{intro-equation-duality-ref-proof}

 The transition of the Gibbs-sampler with target measure $ \Pi_t$ on $E:=(E_1\times E_2)$ is defined by
 \begin{equation}\label{def-GGibbs}
 \SS_t((z_1,z_2),d(\overline{z}_1,\overline{z}_2)):=
 \BB_t(z_2,d\overline{z}_1)~\AC_t(\overline{z}_1,d\overline{z}_2)~
 \end{equation}
 This transition is summarized in the following synthetic diagram
 $$
 \left(\begin{array}{l}
 z_1\\
 z_2
 \end{array}\right)\longrightarrow \left(\begin{array}{l}
 \overline{z}_1~ \sim~(\YY_t~|~\YY_t^{\Ia}=z_2)\\
 z_2
 \end{array}\right)\longrightarrow \left(\begin{array}{l}
  \overline{z}_1\\
  \overline{z}_2~\sim~m(\overline{z}_1)
 \end{array}\right)
 $$
 By construction, we have the duality property
 \begin{equation}\label{def-GGibbs-rev}
    \Pi_t(d(z_1,z_2))~ \SS_t((z_1,z_2),d(\overline{z}_1,\overline{z}_2))=    \Pi_t(d(\overline{z}_1,\overline{z}_2))~  \SS^{-}_t((\overline{z}_1,\overline{z}_2),d(z_1,z_2))
 \end{equation}
  with the backward transition
 $$
 \SS^{-}_t((\overline{z}_1,\overline{z}_2),d(z_1,z_2))
=\AC_t(\overline{z}_1,dz_2)~\BB_t(z_2,dz_1)
 $$

 Integrating (\ref{def-GGibbs-rev}) w.r.t. $\overline{z}_1$ we also have the reversibility property
  $$
 \QQ_t(dz_2)~\KK_t(z_2,d\overline{z}_2)~= \QQ_t(d\overline{z}_2)~ \KK_t(\overline{z}_2,dz_2)
 $$
 with the Markov transition $ \KK_t=\BB_t\circ \AC_t$ from $D_t(S)$ into itself defined for any $F\in\Ba(D_t(S))$ and $z\in D_t(S)$ by
$$
\KK_t(F)(z):=\int~ \KK_t(z,d\overline{z})~F(\overline{z})=\EE\left(m(\YY_t)(F)~|~\YY^{\Ia}_t=z\right)
$$
We further assume that the Markov transitions of $X_t$ satisfy condition $(H_0)$ and thus $(H_1)$. In this situation, using corollary~\ref{cor-ref-gibbs-2}, for any time horizon $t\geq 0$,
any function $F\in \mbox{\rm Osc}(D_t(S))$ and and $n\geq 1$
we check that
$$
\Vert\KK_t(F)-\QQ_t(F)\Vert\leq c~(t\vee 1)/N~$$
  for some finite constant $c$ whose value doesn't depend on the parameters $(F,t,N)$.
This implies that
$$
\mbox{\rm osc}(\KK_t(F))\leq c~(t\vee 1)/N~\mbox{\rm osc}(F)\quad
\mbox{\rm and}\quad   \beta_{\mbox{\tiny\rm dob}}(\KK_t)\leq c~(t\vee 1)/N 
$$
For any $\mu\in \Pa(D_t(S))$, we conclude that
$$
\Vert\mu \KK_t^n-\QQ_t\Vert_{\tiny\rm tv}\leq (c~(t\vee 1)/N)^n\times \Vert
\mu-\QQ_t\Vert_{\tiny\rm tv}.
$$
The last assertion is a direct consequence of the contraction inequalities (\ref{ref-contraction}).
   
   \subsubsection*{Acknowledgements}

The authors are  supported by the ANR Quamprocs on quantitative analysis of metastable processes.
P. Del Moral is also supported in part from the Chair Stress Test, RISK Management and Financial Steering, led by the French Ecole polytechnique and its Foundation and sponsored by BNP Paribas. 

We also thank the anonymous reviewers for their excellent suggestions for improving the paper.
Their detailed comments greatly improved the 
presentation of the article.

   \section*{Appendix}
   
   \subsection*{Proof of (\ref{MM-eta})}\label{MM-eta-proof}
   Observe that
    $$
     \begin{array}{l}
  \displaystyle\MM_t(F)=\eta_0Q_{0,t}(f)+  \sum_{n\geq 1}\int_{0\leq s_1<\ldots<s_n\leq t}~ \left[\prod_{0\leq p< n}  \eta_{s_p}Q_{s_p,s_{p+1}}(V_{s_{p+1}}) \right] ~(\eta_{s_{n}}Q_{s_{n},t})(f)
~ds_1\ldots ds_n
   \end{array}
  $$
  Recalling that
  $$
  \eta_{s}Q_{s,t}(f)=\frac{  \gamma_{s}Q_{s,t}(f)}{\gamma_s(1)}=\frac{  \gamma_{t}(f)}{\gamma_s(1)}=
 \frac{\gamma_{t}(1)  }{\gamma_s(1)}\times  \eta_{t}(f)
  $$
  we find that
     \begin{eqnarray*}
  \displaystyle \MM_t(F)&=&\gamma_t(f)+ \sum_{n\geq 1}\int_{0\leq s_1<\ldots<s_n\leq t}~ \left[\prod_{0\leq p< n}  \eta_{s_{p+1}}(V_{s_{p+1}}) \times \frac{\gamma_{s_{p+1}}(1)}{\gamma_{s_p}(1)}\right] ~\frac{\gamma_t(f)}{\gamma_{s_n}(1)}
~ds_1\ldots ds_n\\
&=&\gamma_t(f)~\left(1+\sum_{n\geq 1}\int_{0\leq s_1<\ldots<s_n\leq t}~ \left[\prod_{0\leq p< n}  \eta_{s_{p+1}}(V_{s_{p+1}})\right] ds_1\ldots ds_n\right)
   \end{eqnarray*}
   This ends the proof of (\ref{MM-eta}).
   \cqfd
  \subsection*{Proof of lemma~\ref{lem-H0-perturbation}}\label{lem-H0-perturbation-proof}

 Let $X_{s,t}(x)$, with $t\geq s$, be the stochastic flow associated with the generator $L_t$ starting at $X_{s,s}(x)=x$ at time $t=s$.
  In this notation, we  have the perturbation formula
 $$
 \begin{array}{l}
 P^{\epsilon}_{s,t}(f)(x)\\
 \\
 \displaystyle=\EE\left[ f(X_{s,t}(x))~
 e^{-\epsilon\int_s^t\lambda(X_{s,u}(x))du}\right]+\int_s^t
 \EE\left[\epsilon\lambda(X_{s,u}(x))~ e^{-\epsilon\int_s^u\lambda(X_{s,v}(x))dv} K_u( P^{\epsilon}_{u,t}(f))(X_{s,u}(x))\right]
~du
 \end{array}$$
 For non negative functions $f$ and any $t\geq 0$ and $h>0$ we have
  \begin{eqnarray*}
 P^{\epsilon}_{t,t+h}(f)&\leq& 
 e^{-\epsilon\lambda_1
  h} P_{t,t+h}(f)+\epsilon\lambda_2\varpi_2\int_t^{t+h}e^{-\epsilon\lambda_1 (u-t)}~\kappa_u P^{\epsilon}_{u,t+h}(f)~du\\
  &\leq &e^{\epsilon(\lambda_2-\lambda_1)
  h}\left[e^{-\epsilon\lambda_2h} P_{t,t+h}(f)+\epsilon\lambda_2\varpi_2\int_t^{t+h}e^{-\epsilon\lambda_2 (u-t)}~\kappa_u P^{\epsilon}_{u,t+h}(f)~du\right]
   \end{eqnarray*}
In the same vein, we have
  \begin{eqnarray*}
  P^{\epsilon}_{t,t+h}(f)\geq (\lambda_1/\lambda_2)(\varpi_1/\varpi_2)\left[
 e^{-\epsilon\lambda_2
  h} P_{t,t+h}(f)+\epsilon\lambda_2\varpi_2\int_t^{t+h}e^{-\epsilon\lambda_2 (u-t)}~\kappa_u P^{\epsilon}_{u,t+h}(f)~du\right]
   \end{eqnarray*}
This shows that
$$
\rho_{\epsilon}(h)~\leq \frac{d(\delta_xP^{\epsilon}_{t,t+h})}{d\mu^{\epsilon}_{t,h}}(y)\leq \rho_{\epsilon}(h)^{-1}
$$
 with the probability measure
 $$\mu_{t,h}^{\epsilon}~\propto~ e^{-\epsilon\lambda_2h} \mu_{t,h}+\epsilon\lambda_2\varpi_2\int_t^{t+h}e^{-\epsilon\lambda_2 (u-t)}~\kappa_u P^{\epsilon}_{u,t+h}(f)~du
 $$
 This ends the proof of the lemma.\cqfd

\subsection*{Proof of lemma~\ref{formula-remainder-jumps-lem}}\label{formula-remainder-jumps-proof}

For any functions $f_i\in \mbox{\rm Osc}(S)$ and any $l\leq k$ we have
$$
\begin{array}{l}
\displaystyle\left\vert \EE\left[\prod_{1\leq l\leq k}\Delta m(\xi_t)(f_l)~|~\Fa_{t-}\right]\right\vert
 \displaystyle\leq \frac{1}{N^k}~\left[\sum_{i\in [N]} \left(V_t(\xi_t^i)+\lambda_t(\xi^i_t)\right)\right]\leq \frac{1}{N^{k-1}}~\Vert \lambda+V\Vert dt\\
\\
\displaystyle\Longrightarrow N^k~\vert\Upsilon^{k+1}_{m(\xi_{s-})} \phi_{s,t}(f)\vert=\frac{N^k}{(k+1)!}~\left\vert \partial_s\EE\left[ (\Delta m(\xi_{s}))^{\otimes (k+1)}~\overline{\partial}^{(k+1)}_{m(\xi_{s-})}\phi_{s,t}(f)~|~\Fa_{s-}\right]\right\vert\\
\\
\hskip1cm\displaystyle=N^k~
\left\vert \partial_s\EE\left[ \frac{1}{m(\xi_s) Q^{m(\xi_{s-})}_{s,t}(1)}~\left(
\Delta m(\xi)\left(Q^{m(\xi_{s-})}_{s,t}(1)\right)\right)^k~\Delta m(\xi)\partial_{m(\xi_{s-})}\phi_{s,t}(f)~|~\Fa_{s-}\right]\right\vert\\
\\
\hskip1cm\displaystyle \leq ~e^{(2+k)q}~\Vert \lambda+V\Vert 
\end{array}
$$
\subsection*{Proof of theorem~\ref{cor-ref-zeta-m}}\label{cor-ref-zeta-m-proof}
We use  (\ref{ref-pr-final}) to check that
$$
\begin{array}{l}
 \displaystyle~m(\xi_s)\Gamma_{L_{s,m(\xi_s)}}\left(  Q^{m(\xi_{s})}_{s,t}(1),\partial_{m(\xi_{s})}\phi_{s,t}(f)\right) \\
 \\
\displaystyle = (\eta_s Q_{s,t}^{m(\xi_{s})}(1))^{2}~m(\xi_s)\Gamma_{L_{s,m(\xi_s)}}\left(Q^{\eta_s}_{s,t}(1), \partial_{\eta_s}\phi_{s,t}(f)\right)\\
\\
\displaystyle\hskip3cm+ (\eta_s Q_{s,t}^{m(\xi_{s})}(1))^{2}~~[\phi_{s,t}(\eta_s)-\phi_{s,t}(m(\xi_s))](f)~m(\xi_s)\Gamma_{L_{s,m(\xi_s)}}\left(Q^{\eta_s}_{s,t}(1)\right)
\end{array}$$
 Using (\ref{ref-osc-partial-phi}) we also have the estimate
\begin{equation}\label{ref-Gamma-estimates}
\begin{array}{l}
 \displaystyle~\left\vert m(\xi_s)\Gamma_{L_{s,m(\xi_s)}}\left(  Q^{m(\xi_{s})}_{s,t}(1),\partial_{m(\xi_{s})}\phi_{s,t}(f)\right) \right\vert\leq e^{3q}~~\mbox{\rm osc}(\overline{Q}_{s,t}(f))~m(\xi_s)\Gamma_{L_{s,m(\xi_s)}}\left(Q^{\eta_s}_{s,t}(1)\right)\\
\\
\displaystyle\hskip5cm+e^{2q}~\sqrt{m(\xi_s)\Gamma_{L_{s,m(\xi_s)}}\left(Q^{\eta_s}_{s,t}(1)\right)}~\sqrt{m(\xi_s)\Gamma_{L_{s,m(\xi_s)}}\left(\partial_{\eta_s}\phi_{s,t}(f)\right)}
\end{array}\end{equation}

On the other hand, we have
$$
\begin{array}{l}
\partial_{\eta_s}\phi_{s,t}(f)=Q^{\eta_s}_{s,t}\left[f-\eta_t(f)\right]
\quad\mbox{\rm
and}
\quad
Q^{\eta_s}_{s,t}(f)(x)=\EE\left(f(X_t)~e^{-\int_s^t \overline{V}_u(X_u)~du}~|~X_s=x\right)\\
\\
\displaystyle\Longrightarrow
\partial_sQ^{\eta_s}_{s,t}(f)=-L_s(Q^{\eta_s}_{s,t}(f))+\overline{V}_s~Q^{\eta_s}_{s,t}(f)\\
\\
\displaystyle\Longrightarrow
\partial_s(Q^{\eta_s}_{s,t}(f)Q^{\eta_s}_{s,t}(g))=- Q^{\eta_s}_{s,t}(f)~L_s(Q^{\eta_s}_{s,t}(g))-Q^{\eta_s}_{s,t}(g)~L_s(Q^{\eta_s}_{s,t}(f))+2 \overline{V}_s~Q^{\eta_s}_{s,t}(f)Q^{\eta_s}_{s,t}(g)
\end{array}
$$
We also have
$$
L_{t,\mu}(f)=L_t(f)+V_t~[\mu(f)-f]\Longleftrightarrow
L_t(f)-L_{t,\mu}(f)=V_t~[f-\mu(f)]
$$
This yields the formula
\begin{eqnarray*}
\eta\Gamma_{L_{t,\eta}}(f,g)-\eta\Gamma_{L_t}(f,g)&=&\int~\eta(dx)~\eta(dy)~V_t(y)~[f(y)-f(x)][g(y)-g(x)]\\
&=&\eta(V_t (fg))+\eta(V_t)~ \eta(fg)-\eta(fV_t)~\eta(g)-\eta(gV_t)~\eta(f)
\end{eqnarray*}
For any given time horizon $t$ and $s\in [0,t]$ we have
$$
\begin{array}{l}
dm(\xi_s)(Q^{\eta_s}_{s,t}(f)Q^{\eta_s}_{s,t}(g))-\frac{1}{\sqrt{N}}~dM_s(Q^{\eta_{\point}}_{\point,t}(f)Q^{\eta_{\point}}_{\point,t}(g))\\
\\
\displaystyle=m(\xi_s)\left[L_{s,m(\xi_s)}(Q^{\eta_s}_{s,t}(f)Q^{\eta_s}_{s,t}(g))- Q^{\eta_s}_{s,t}(f)~L_s(Q^{\eta_s}_{s,t}(g))\right.\\
\\
\hskip5cm\left.-Q^{\eta_s}_{s,t}(g)~L_s(Q^{\eta_s}_{s,t}(f))+2 \overline{V}_s~Q^{\eta_s}_{s,t}(f)Q^{\eta_s}_{s,t}(g)\right]~ds
\end{array}$$
This yields
$$
\begin{array}{l}
dm(\xi_s)(Q^{\eta_s}_{s,t}(f)Q^{\eta_s}_{s,t}(g))-\frac{1}{\sqrt{N}}~dM_s(Q^{\eta_{\point}}_{\point,t}(f)Q^{\eta_{\point}}_{\point,t}(g))\\
\\
\displaystyle=m(\xi_s)\left[\Gamma_{s,L_{s,m(\xi_s)}}(Q^{\eta_s}_{s,t}(f),Q^{\eta_s}_{s,t}(g))+ V_s~Q^{\eta_s}_{s,t}(f)~~(m(\xi_s)Q^{\eta_s}_{s,t}(g)-Q^{\eta_s}_{s,t}(g))
\right.\\
\\
\left.\hskip3cm+V_t~Q^{\eta_s}_{s,t}(g)~~(m(\xi_s)Q^{\eta_s}_{s,t}(f)-Q^{\eta_s}_{s,t}(f))+2 \overline{V}_s~Q^{\eta_s}_{s,t}(f)Q^{\eta_s}_{s,t}(g)\right]~ds
\end{array}$$
from which we check that
$$
\begin{array}{l}
\displaystyle m(\xi_t)(fg)-m(\xi_0)(Q^{\eta_0}_{0,t}(f)Q^{\eta_0}_{0,t}(g))-\frac{1}{\sqrt{N}}~M_t(Q^{\eta_{\point}}_{\point,t}(f)Q^{\eta_{\point}}_{\point,t}(g))\\
\\
\displaystyle=\int_0^t~m(\xi_s)\Gamma_{s,L_{s,m(\xi_s)}}(Q^{\eta_s}_{s,t}(f),Q^{\eta_s}_{s,t}(g))~ds\\
\\
\displaystyle+\int_0^t~m(\xi_s)\left[V_s~Q^{\eta_s}_{s,t}(f)~~(m(\xi_s)Q^{\eta_s}_{s,t}(g)-Q^{\eta_s}_{s,t}(g))
\right.\\
\\
\left.\hskip2cm+V_s~Q^{\eta_s}_{s,t}(g)~~(m(\xi_s)Q^{\eta_s}_{s,t}(f)-Q^{\eta_s}_{s,t}(f))+2 (V_s-\eta_s(V_s))~Q^{\eta_s}_{s,t}(f)Q^{\eta_s}_{s,t}(g)\right]~ds
\end{array}$$
After some simplifications we check that
$$
\begin{array}{l}
\displaystyle \int_0^t~m(\xi_s)\Gamma_{s,L_{s,m(\xi_s)}}(Q^{\eta_s}_{s,t}(f),Q^{\eta_s}_{s,t}(g))~ds\\
\\
\displaystyle=m(\xi_t)(fg)-m(\xi_0)(Q^{\eta_0}_{0,t}(f)Q^{\eta_0}_{0,t}(g))-\frac{1}{\sqrt{N}}~M_t(Q^{\eta_{\point}}_{\point,t}(f)Q^{\eta_{\point}}_{\point,t}(g))\\
\\
\displaystyle+\int_0^t~\left[2\eta_s(V_s)~m(\xi_s)(Q^{\eta_s}_{s,t}(f)Q^{\eta_s}_{s,t}(g))\right.\\
\\
\left.\displaystyle\hskip2cm-m(\xi_s)(V_s~Q^{\eta_s}_{s,t}(f))~~m(\xi_s)Q^{\eta_s}_{s,t}(g)
-m(\xi_s)(V_s~Q^{\eta_s}_{s,t}(g))~~m(\xi_s)Q^{\eta_s}_{s,t}(f)\right]~ds
\end{array}$$
Choosing $f=g=1$ and taking the expectations we find that
$$
\begin{array}{l}
\displaystyle \int_0^t~\EE\left[m(\xi_s)\Gamma_{s,L_{s,m(\xi_s)}}(Q^{\eta_s}_{s,t}(1))\right]~ds\\
\\
\displaystyle=1-\eta_0(Q^{\eta_0}_{0,t}(1)^2)+2\int_0^t~\EE\left[\eta_s(V_s)~m(\xi_s)(Q^{\eta_s}_{s,t}(1)^2)-
m(\xi_s)(V_s~Q^{\eta_s}_{s,t}(1))~~m(\xi_s)Q^{\eta_s}_{s,t}(1)\right]~ds\\
\\
\leq 1+2e^{2q}~\Vert V\Vert~t
\end{array}$$
Choosing $f=g=h-\eta_t(h)$, with $h\in\mbox{\rm Osc}(S)$ and taking the expectations we find that
$$
\begin{array}{l}
\displaystyle \int_0^t~\EE\left[m(\xi_s)\Gamma_{s,L_{s,m(\xi_s)}}(\partial_{\eta_s}\phi_{s,t}(h))\right]~ds\\
\\
\displaystyle=\EE\left[m(\xi_t)([h-\eta_t(h)]^2)\right]-\eta_0([\partial_{\eta_0}\phi_{0,t}(h)]^2)\\
\\
\displaystyle+2\int_0^t~\EE\left[\eta_s(V_s)~m(\xi_s)([\partial_{\eta_s}\phi_{s,t}(h)]^2)-m(\xi_s)(V_s~\partial_{\eta_s}\phi_{s,t}(h))~~m(\xi_s)(\partial_{\eta_s}\phi_{s,t}(h))
\right]~ds\\
\\
\leq 1+4e^{2q}~\Vert V\Vert~t
\end{array}$$

For any $f\in \mbox{\rm Osc}(S)$ combining (\ref{ref-Gamma-estimates}) with Cauchy-Schwartz inequality we find that that
$$
\begin{array}{l}
 \displaystyle~\int_0^t~\vert m(\xi_s)\Gamma_{L_{s,m(\xi_s)}}\left(  Q^{m(\xi_{s})}_{s,t}(1),\partial_{m(\xi_{s})}\phi_{s,t}(f)\right)\vert~ds\leq 2e^{3q}~\left[1+4e^{2q}~\Vert V\Vert\right]~t
\end{array}
$$
Combining the above estimate with (\ref{formula-remainder-jumps}) and corollary~\ref{cor-EE-m} whenever $(H_2)$ is met we conclude that
$$
\vert \EE(m(\xi_t)(f))-\EE(\phi_{0,t}(m(\xi_0))(f))\vert\leq 2e^{3q}\left(1+4e^{2q}~\Vert V\Vert+\frac{1}{N^2}~2e^{q}~\Vert \lambda+V\Vert\right)~t/N
$$
We further assume that $(H_1)$ is satisfied. In this case, using (\ref{ref-Gamma-estimates}) we also have
$$
\begin{array}{l}
 \displaystyle~\vert m(\xi_s)\Gamma_{L_{s,m(\xi_s)}}\left(  Q^{m(\xi_{s})}_{s,t}(1),\partial_{m(\xi_{s})}\phi_{s,t}(f)\right) \vert\leq e^{3q}~\alpha~e^{-\beta(t-s)}~m(\xi_s)\Gamma_{L_{s,m(\xi_s)}}\left(Q^{\eta_s}_{s,t}(1)\right)\\
\\
\displaystyle\hskip5cm+~e^{2q}~\sqrt{m(\xi_s)\Gamma_{L_{s,m(\xi_s)}}\left(Q^{\eta_s}_{s,t}(1)\right)}~\sqrt{m(\xi_s)\Gamma_{L_{s,m(\xi_s)}}\left(\partial_{\eta_s}\phi_{s,t}(f)\right)}\end{array}
$$

For any $\widetilde{\beta}\in \RR$ we set
\begin{eqnarray*}
\widetilde{Q}^{\eta_s}_{s,t}(f)(x)&:=&e^{\widetilde{\beta}(t-s)}~
Q^{\eta_s}_{s,t}(f)(x)\\
&=&\EE\left(f(X_t)~e^{-\int_s^t \widetilde{V}_u(X_u)~du}~|~X_s=x\right)\quad \mbox{\rm with}\quad
\widetilde{V}_t(x)=\overline{V}_t(x)-\widetilde{\beta}
\end{eqnarray*}
Arguing as above, we have
$$
\partial_s(\widetilde{Q}^{\eta_s}_{s,t}(f)\widetilde{Q}^{\eta_s}_{s,t}(g))=- \widetilde{Q}^{\eta_s}_{s,t}(f)~L_s(\widetilde{Q}^{\eta_s}_{s,t}(g))-\widetilde{Q}^{\eta_s}_{s,t}(g)~L_s(\widetilde{Q}^{\eta_s}_{s,t}(f))+2 ~\widetilde{V}_s~\widetilde{Q}^{\eta_s}_{s,t}(f)~\widetilde{Q}^{\eta_s}_{s,t}(g)
$$
and
$$
\begin{array}{l}
\displaystyle dm(\xi_s)(\widetilde{Q}^{\eta_s}_{s,t}(f)\widetilde{Q}^{\eta_s}_{s,t}(g))-\frac{1}{\sqrt{N}}~dM_s(\widetilde{Q}^{\eta_{\point}}_{\point,t}(f)\widetilde{Q}^{\eta_{\point}}_{\point,t}(g))\\
\\
\displaystyle=m(\xi_s)\left[\Gamma_{s,L_{s,m(\xi_s)}}(\widetilde{Q}^{\eta_s}_{s,t}(f),\widetilde{Q}^{\eta_s}_{s,t}(g))+ V_s~\widetilde{Q}^{\eta_s}_{s,t}(f)~~(m(\xi_s)\widetilde{Q}^{\eta_s}_{s,t}(g)-\widetilde{Q}^{\eta_s}_{s,t}(g))
\right.\\
\\
\left.\hskip5cm+V_t~\widetilde{Q}^{\eta_s}_{s,t}(g)~~(m(\xi_s)\widetilde{Q}^{\eta_s}_{s,t}(f)-\widetilde{Q}^{\eta_s}_{s,t}(f))+2 \widetilde{V}_s~\widetilde{Q}^{\eta_s}_{s,t}(f)~\widetilde{Q}^{\eta_s}_{s,t}(g)\right]~ds
\end{array}$$

This implies that
$$
\begin{array}{l}
\displaystyle\int_0^t~e^{2\widetilde{\beta}(t-s)}~m(\xi_s)\Gamma_{s,L_{s,m(\xi_s)}}(Q^{\eta_s}_{s,t}(f),Q^{\eta_s}_{s,t}(g))~ds\\
\\
\displaystyle=m(\xi_t)(fg)-e^{\widetilde{\beta}t}~m(\xi_0)(Q^{\eta_0}_{0,t}(f)Q^{\eta_0}_{0,t}(g))-\frac{1}{\sqrt{N}}~M_t(\widetilde{Q}^{\eta_{\point}}_{\point,t}(f)\widetilde{Q}^{\eta_{\point}}_{\point,t}(g))\\
\\
\displaystyle+\int_0^t~e^{2\widetilde{\beta}(t-s)}~\left[2 (\eta_s(V_s)+\widetilde{\beta})~
m(\xi_s)(Q^{\eta_s}_{s,t}(f)Q^{\eta_s}_{s,t}(g))-m(\xi_s)(V_s~Q^{\eta_s}_{s,t}(f))~~m(\xi_s)Q^{\eta_s}_{s,t}(g)
\right.\\
\\
\left.\hskip3cm-m(\xi_s)(V_s~Q^{\eta_s}_{s,t}(g))~~m(\xi_s)Q^{\eta_s}_{s,t}(f)\right]~ds
\end{array}$$
Choosing $f=g=1$ and $\widetilde{\beta}<0$ we have
$$
\begin{array}{l}
\displaystyle\int_0^t~e^{2\widetilde{\beta}(t-s)}~\EE\left[m(\xi_s)\Gamma_{s,L_{s,m(\xi_s)}}(Q^{\eta_s}_{s,t}(1))\right]~ds\\
\\
\displaystyle=1-e^{\widetilde{\beta}t}~\eta_0(Q^{\eta_0}_{0,t}(1)^2)\\
\\
\displaystyle+2~\int_0^t~e^{2\widetilde{\beta}(t-s)}~\EE\left[(\eta_s(V_s)+\widetilde{\beta})~
m(\xi_s)(Q^{\eta_s}_{s,t}(1)^2)-m(\xi_s)(V_s~Q^{\eta_s}_{s,t}(1))~~m(\xi_s)Q^{\eta_s}_{s,t}(1)
\right]~ds\\
\\
\displaystyle\leq 1+e^{2q}~\left(1+2\vert\widetilde{\beta}\vert^{-1}\Vert V\Vert \right)=1+e^{2q}~\left(1+4\beta^{-1}\Vert V\Vert \right)
\quad \mbox{\rm when}\quad \widetilde{\beta}=-\beta/2
\end{array}$$
Choosing $f=g=[h-\eta_t(h)]$, with $h\in \mbox{\rm Osc}(S)$ and $0<\widetilde{\beta}<\beta$ we have
$$
\begin{array}{l}
\displaystyle\int_0^t~e^{2\widetilde{\beta}(t-s)}~\EE\left[m(\xi_s)\Gamma_{s,L_{s,m(\xi_s)}}(\partial_{\eta_s}\phi_{s,t}(h)\right]~ds\\
\\
\displaystyle\leq \EE\left[m(\xi_t)([h-\eta_t(h)]^2)\right]\\
\\
\displaystyle+2\int_0^t~e^{2\widetilde{\beta}(t-s)}~\left[(\eta_s(V_s)+\widetilde{\beta})~
m(\xi_s)\left(\left[\partial_{\eta_s}\phi_{s,t}(h)\right]^2\right)-m(\xi_s)(V_s~\partial_{\eta_s}\phi_{s,t}(h))~~m(\xi_s)\partial_{\eta_s}\phi_{s,t}(h)\right]~ds\\
\\
\displaystyle\leq 1+2r^2(2\Vert V\Vert+\widetilde{\beta})\int_0^t~e^{-2(\beta-\widetilde{\beta})(t-s)}~ds\\
\\
\displaystyle\leq 1+r^2(2\Vert V\Vert+\widetilde{\beta})~(\beta-\widetilde{\beta})^{-1}=1+r^2(4\Vert V\Vert \beta^{-1}+1)
\quad \mbox{\rm when}\quad \widetilde{\beta}=\beta/2
\end{array}$$
We end the proof of the theorem using the fact that
$$
\begin{array}{l}
 \displaystyle~\vert m(\xi_s)\Gamma_{L_{s,m(\xi_s)}}\left(  Q^{m(\xi_{s})}_{s,t}(1),\partial_{m(\xi_{s})}\phi_{s,t}(f)\right) \vert\\
\\
\leq e^{3q} (1+\alpha)~e^{-\beta (t-s)/2}
m(\xi_s)\Gamma_{L_{s,m(\xi_s)}}\left(Q^{\eta_s}_{s,t}(1)\right)+e^{2q}~e^{\beta (t-s)/2}
m(\xi_s)\Gamma_{L_{s,m(\xi_s)}}\left(\partial_{\eta_s}\phi_{s,t}(f)\right)
\end{array}
$$
In the last assertion we have used the fact that the estimate $\sqrt{ab}\leq c a+b/c$, for all $a,b,c>0$.
 This ends the proof of the theorem.
 \cqfd
\subsection*{Proof of (\ref{theo-a-e-ref})}\label{theo-a-e-ref-proof}

The proof follows the same lines of arguments as the ones used in the proof of corollary~\ref{cor-ref-app-proof}.

 By (\ref{ref-wH-m-k}) and (\ref{ref-Ja-2-k}), given $\Ia=k$, the historical flow $\widehat{\Xa}_{s,t}$, the jump times $T_{k,n}^{s,x}=r_n$ as well as the randomly selected coalescent maps $a_n=c_{\iota_n}$ for some $\iota_n\in [N]_0^2$,  for any $r_n\leq t<r_{n+1}$ we have
 \begin{equation}
   \widehat{\goodchi}_{s,t}=   \widehat{\Xa}_{r_{n},t}\circ \widehat{\Ta}^{\,\iota,(k,n)}_{r,r_{n}}\quad\mbox{\rm and}\quad
      \YY_{s,t}=   \widehat{\Xa}_{r_{n},t}\circ \widehat{\Xa}^{\,\iota,(k,n)}_{r,r_{n}}
 \label{def-wT-wX-k}
 \end{equation}
with the semigroups $( \widehat{\Ta}^{\,\iota,(k,n)}_{r,r_{n}}, \widehat{\Xa}^{\,\iota,(k,n)}_{r,r_{n}})$
 defined as $( \widehat{\Ta}^{\,\iota,(n)}_{r,r_{n}}, \widehat{\Xa}^{\,\iota,(n)}_{r,r_{n}})$ by replacing in (\ref{def-wT-wX-2}) the maps $( \widehat{\Ta}_{s,t}^{\,a}, \widehat{\Xa}_{s,t}^{\,a} )$ by the maps  $( \widehat{\Ta}_{s,t}^{\,k,a}, \widehat{\Ta}_{s,t}^{\,k,a} )$ introduced in (\ref{def-wTa}). 
This yields the almost sure formula
$$
  \begin{array}{l}
\displaystyle\int_{[s,t]_n} \left(\left(\widehat{\UU}^{(N)}_{r_0,r_1}\widehat{\Ja}^{(2)}_{\Ia,r_1}\right)\ldots \left(\widehat{\UU}^{(N)}_{r_{n-1},r_{n}}\widehat{\Ja}^{(2)}_{\Ia,r_n}\right)\left(\widehat{\UU}^{(N)}_{r_n,t}(F)\right)\right)(x)~dr\\
\\
\displaystyle=\sum_{\iota\in [N]_0^{(2n)}}~\int_{[s,t]_n}~F\left(  \left(
\widehat{\Xa}_{r_{n},t}\circ \widehat{\Ta}^{\,\iota,(\Ia,n)}_{r,r_{n}}\right)(x),\left(\widehat{\Xa}_{r_{n},t}\circ \widehat{\Xa}^{\,\iota,(\Ia,n)}_{r,r_{n}}\right)(y)
\right)~\mathfrak{q}^{s,x}_{\Ia,n}(\widehat{\Xa}_{s,t},d(\iota,r))
\end{array}
$$
with the random measures
$$
 \begin{array}{l}
 \displaystyle\mathfrak{q}^{s,x}_{\Ia,n}(\widehat{\Xa}_{s,t},d(\iota,r))\\
 \\
: =\displaystyle
  \Za_{r_0,r_{1}}\left( \widehat{\Xa}_{r_0,r_{1}}(x)\right)^N~\widehat{\lambda}_{\Ia,r_1}^{i_1,j_1}\left(\widehat{\Xa}_{r_0,r_{1}}(x)\right) ~\\
 \\
 \displaystyle \times~\Za_{r_1,r_{2}}\left(    \widehat{\Xa}_{r_{1},r_2}\left(\widehat{\Xa}_{r_{0},r_1}^{\,\Ia,c_{\iota_1}} (x)\right)\right)^N~\widehat{\lambda}_{\Ia,r_2}^{i_2,j_2}\left(\widehat{\Xa}_{r_{1},r_2}\left(\widehat{\Xa}_{r_{0},r_1}^{\,\Ia,c_{\iota_1}} (x)\right)\right)\\
  \\
  \displaystyle\times\ldots  \times~\Za_{r_{n-1},r_{n}}\left( \widehat{\Xa}_{r_{n-1},r_n}\left(\widehat{\Xa}^{\,\iota,(\Ia,n-1)}_{r,r_{n-1}}(x)\right)\right)^N~\widehat{\lambda}_{\Ia,r_n}^{i_n,j_n}\left( \widehat{\Xa}_{r_{n-1},r_n}\left(\widehat{\Xa}^{\,\iota,(\Ia,n-1)}_{r,r_{n-1}}(x)\right)\right)~\\
  \\
  \displaystyle\times~\Za_{r_n,t}\left(\widehat{\Xa}_{r_{n},t}\left(\widehat{\Xa}^{\,\iota,(n)}_{\Ia,(r,r_{n})}(x)\right)\right)^N~dr
 \end{array}
$$
On the other hand, using (\ref{def-lambda-Ia-k}) we have
$$
\begin{array}{l}
\Za_{r_{n-1},r_{n}}\left( \widehat{\Xa}_{r_{n-1},r_n}\left(\widehat{\Xa}^{\,\iota,(\Ia,n-1)}_{r,r_{n-1}}(x)\right)\right)^N\\
\\
=Z_{r_{n-1},r_{n}}\left(\widehat{\Xa}^k_{s,r_n}(x)\right)\times 
\Ea^{k}_{r_{n-1},r_{n}}\left( \widehat{\Xa}_{r_{n-1},r_n}\left(\widehat{\Xa}^{\,\iota,(\Ia,n-1)}_{r,r_{n-1}}(x)\right)\right)
\end{array}
$$
with the stochastic exponential functional
$$
\Ea^k_{s,t}\left(\widehat{\Xa}_{s,t}(x)\right):=\exp{\left(- \int_s^t\widehat{\lambda}_{k,r}\left(\widehat{\Xa}_{s,r}\left(x\right)\right)~dr\right)}\quad \mbox{\rm with}\quad \widehat{\lambda}_{\Ia,t}(x)= (N-1)~m(x^{-\Ia})(\widehat{V}_t)
$$
The last assertion comes from the fact that
$$
\widehat{\Xa}^k_{r_{n-1},r_n}\left(\widehat{\Xa}^{\,\iota,(\Ia,n-1)}_{r,r_{n-1}}(x)\right)=\widehat{\Xa}^k_{s,r_n}(x)\quad \mbox{\rm by (\ref{ref-ancestral-k-proof})}
$$
We conclude that
$$
\mathfrak{q}^{s,x}_{\Ia,n}(\widehat{\Xa}_{s,t},d(\iota,r))=Z_{s,t}\left(\widehat{\Xa}^k_{s,t}(x)\right)~\times~\mathfrak{p}^{s,x}_{\Ia,n}(\widehat{\Xa}_{s,t},d(\iota,r))
$$
with the conditional probability measures of the jump times and coalescent maps
$$
 \begin{array}{l}
 \displaystyle\mathfrak{p}^{s,x}_{\Ia,n}(\widehat{\Xa}_{s,t},d(\iota,r))\\
 \\
 :=\PP\left((T^{s,x}_{\Ia,1},\ldots,T^{s,x}_{\Ia,n})\in d(r_1,\ldots,r_n),~~(a_1,\ldots,a_n)=(c_{\iota_1},\ldots, c_{\iota_n}),~~1_{T^{s,x}_{\Ia,n+1}>t}~|~\widehat{\Xa}_{s,t},~\Ia\right)
 \end{array}
$$
We end the proof of  (\ref{theo-a-e-ref}) taking the expectations and using (\ref{inter-prop}).
\cqfd


\begin{thebibliography}{10}

\bibitem{aldous}
D. Aldous, U. Vazirani. Go with the winners algorithms. In Foundations of Computer Science. Proceedings., 35th Annual Symposium on (pp. 492--501). IEEE (1994).

\bibitem{adh-2010}
C. Andrieu, A. Doucet,  R. Holenstein. Particle Markov chain Monte Carlo for efficient
numerical simulation. In L'Ecuyer, P. and Owen, A. B., editors, Monte Carlo and Quasi-Monte
Carlo Methods 2008, pp. 45--60. Spinger-Verlag Berlin Heidelberg (2009).

\bibitem{adm-18}
M. Arnaudon, P.  Del Moral. A variational approach to nonlinear and interacting diffusions. Stochastic Analysis and Applications, vol. 37, no. 5, pp. 717--748 (2019).

\bibitem{adm-19}
Marc Arnaudon, Pierre Del Moral.
A second order analysis of McKean-Vlasov semigroups. to appear in Annals of Appl. Probab. (2020).

\bibitem{aronson}
D. G. Aronson, Bounds for the fundamental solution of a parabolic equation, Bulletin of American Math. Soc., vol. 73, pp. 890--896 (1967).

\bibitem{bally-Cont}
V. Bally, L. Caramellino, R. Cont. Stochastic integratio by parts and functional It\^o calculus.
Advanced courses in Mathematics, Barcelona, C.R.M. Birkh\"auser (2016).
 
 \bibitem{bene-1}
 D. Benedetto, E. Caglioti, M. Pulvirenti. A kinetic equation for granular media. RAIRO Mod\`el. Math. Anal. Num\'er. vol. 31, no. 5, pp. 615--641 (1997).
 
  \bibitem{bene-2}
   D. Benedetto, E. Caglioti, E., Carrillo, M. Pulvirenti. A non-Maxwellian steady distribution for one-dimensional granular media. J. Statist. Phys.vol.  91, pp. 979--990 (1998).
   
      \bibitem{gbg-13}
F. Bolley, I. Gentil, A. Guillin. Uniform convergence to equilibrium for granular media.
Archive for Rational Mechanics and Analysis, vol. 208, no. 2, pp. 429--445 (2013).
   
   \bibitem{beskos}
  A. Beskos, G.O. Roberts.  Exact simulation of diffusions. The Annals of Applied Probability, vol. 15, no. 4, pp. 2422--2444 (2005).

   \bibitem{beskos-2}
A. Beskos, O. Papaspiliopoulos, G.O. Roberts. Retrospective exact simulation of diffusion sample paths with applications. Bernoulli, vol. 12, no. 6, 
pp. 1077--1098. (2006). 



   \bibitem{beskos-3}
A. Beskos, O. Papaspiliopoulos, G.O. Roberts. A factorisation of diffusion measure and finite sample path constructions. 
Methodol. Comput. Appl. Probab., vol. 10, pp. 85--104  (2008). 



   \bibitem{beskos-4}
A. Beskos, O. Papaspiliopoulos, G.O. Roberts. P. Fearnhead.  Exact and computationally efficient likelihood-based estimation for discretely observed diffusion processes. J. R. Stat. Soc. Ser. B Stat. Methodol., vol. 68, pp. 333--382 (2006). 

\bibitem{caffarel}
M. Caffarel, R. Assaraf.
A pedagogical introduction to quantum Monte Carlo.
In Mathematical models and methods for ab initio Quantum Chemistry in Lecture Notes in Chemistry, eds. M. Defranceschi and C.Le Bris, Springer p.45 (2000).

\bibitem{casella}
B. Casella, G.O. Roberts. Exact simulation of jump-diffusion processes with Monte Carlo applications. Methodology and Computing in Applied Probability, vol. 13, no. 3, pp. 449--473 (2011).

\bibitem{cances-tony}
E. Canc\`es,  B. Jourdain, T. Leli\`evre. Quantum Monte Carlo simulations of fermions. A mathematical analysis of the fixed-node approximation, Mathematical Models and Methods in Applied Sciences, 16(9), 1403-1440, (2006).

\bibitem{cappe-moulines}
O. Capp\'e, E. Moulines, and T. Ryd\`en. Inference in Hidden Markov Models,
Springer-Verlag (2005).


\bibitem{cmcv-2006}
J. A. Carrillo, R. J. McCann, and C. Villani. Contractions in the 2-Wasserstein length space and thermalization of granular media. Arch. Rational Mech. Anal., vol. 179, pp. 217--263, (2006).

\bibitem{cattiaux}
P. Cattiaux, A. Guillin, and F. Malrieu. Probabilistic approach for granular media equations in the non uniformly convex case. Prob. Theor. Rel. Fields, vol. 140, no. 1-2, pp. 19--40 (2008).

 \bibitem{cordero}
D. Cordero-Erausquin, W. Gangbo, and C. Houdr\'e. Inequalities for generalized entropy and optimal transportation. In Recent Advances in the Theory and Applications of Mass Transport, Contemp. Math. 353. A. M. S., Providence (2004).

\bibitem{cont}
R. Cont, D.A. Fourni\'e. Functional It\^o calculus and stochastic integral representation of martingales. The Annals of Probability, vol.  41, no. 1, pp. 109--133 (2013).

\bibitem{cox}
A.M. Cox, S.C. Harris, E.L. Horton, A. E. Kyprianou. Multi-species neutron transport equation. Journal of Statistical Physics, pp. 1--31,  (2019).

\bibitem{dalayan}
A. S. Dalalyan. Theoretical guarantees for approximate sampling from smooth and log-concave densities. In: Journal of the Royal Statistical Society: Series B (Statistical Methodology), vol. 79, no. 3  pp. 651--676 (2017).

\bibitem{d-2004}
P. {Del Moral}.
\newblock \href{http://www.math.u-bordeaux1.fr/~delmoral/gips.html}{\tt Feynman-{K}ac formulae}.
\newblock  Genealogical and interacting particle systems with applications.
\newblock  Probability and its Applications (New York). (573p.) Springer-Verlag, New
  York (2004).
\bibitem{d-2013}
P. {Del Moral}.
\newblock \href{http://www.math.u-bordeaux1.fr/~pdelmora/Intro+Refs-Mean-Field-Simulation.pdf}{\tt Mean field simulation for Monte Carlo integration.}
\newblock \href{http://www.crcpress.com/product/isbn/9781466504059}{\tt Chapman \& Hall.  Monographs on Statistics \& Applied Probability}  (2013).

\bibitem{d-1996}
P. {Del Moral}. Nonlinear filtering: interacting particle solution.
 {\em Markov Process. \& Related Fields}, vol. 2, no. 4, pp. 555--579 (1996).
 
 \bibitem{dd-2004}
P. {Del Moral}, A. Doucet.Particle motions in absorbing medium with hard and soft obstacles. Stochastic Analysis and Applications, vol. 22, no. 5, pp. 1175--1207  (2004). .
 
 \bibitem{dm-guionnet}
P. Del Moral, A. Guionnet. On the stability of interacting processes with applications to filtering and genetic algorithms. Annales de l'Institut Henri Poincaré (B) Probability and Statistics, vol. 37, no. 2, pp. 155--194 (2001).

 \bibitem{dm-guionnet-2}
P. Del Moral, A. Guionnet. On the stability of measure valued processes with applications to filtering. Comptes Rendus de l'Acad\'emie des Sciences-Series I-Mathematics, vol. 329, no. 5, pp. 429--434 (1999).

\bibitem{dmj-18-1}
P. {Del Moral}, A. Jasra. A sharp first order analysis of Feynman-Kac particle models, Part I: Propagation of chaos. Stochastic Processes and their Applications, 128(1), pp. 332--353 (2018).

\bibitem{dmj-18-2}
P. {Del Moral}, A. Jasra. A sharp first order analysis of Feynman-Kac particle models, Part II: Particle Gibbs samplers. Stochastic Processes and their Applications, 128(1), pp. 354-371  (2018).



  
  \bibitem{dm-jacob-lee-murrau-peters}
  P. Del Moral, P. Jacob, A. Lee, L. Murray, G.W. Peters. 
Feynman-Kac particle integration with geometric interacting jumps, 
Stochastic Analysis and Applications, vol. 31, no. 5, pp. 830--871 (2013). 

\bibitem{dm-jacod-02}
P. Del Moral, J. Jacod. The Monte-Carlo method for filtering with discrete-time observations: Central limit theorems. Numerical Methods and stochastics, vol. 34, pp. 29--53 (2002).
  

\bibitem{DMKP:16}
P. Del Moral, R. Kohn and F. Patras \emph{On particle Gibbs sampler}. Ann. Inst. Henri Poincar\'e Probab. Stat. vol. 52, no. 4, pp. 1687--1733 (2016).

\bibitem{dm-ledoux-miclo}
P.~Del~Moral, M. Ledoux, L.~Miclo. On contraction properties of Markov kernels. Probab. Theory Relat. Fields, vol. 126,  pp. 395--420
(2003).

\bibitem{dm-2000-moran}
P.~Del~Moral, L.~Miclo. A Moran particle system approximation of Feynman-Kac formulae. Stochastic processes and their applications, 86(2), pp. 193--216
(2000).

\bibitem{dm-2000}
P.~Del~Moral, L.~Miclo.
\newblock \href{http://www.math.u-bordeaux1.fr/~delmoral/seminaire.ps}{\tt Branching and interacting particle systems} ap\-pro\-ximations of
  {F}eynman-{K}ac formulae with applications to non-linear filtering.
\newblock In {\em S\'eminaire de {P}robabilit\'es, {XXXIV}}, volume 1729,
  {\em Lecture Notes in Math.}, pages 1--145. Springer, Berlin (2000).
  
  \bibitem{dm-2001-g}
  P.~Del~Moral, L.~Miclo.
  Genealogies and Increasing Propagations of Chaos for Feynman-Kac and Genetic Models.   Annals of Applied Probability, vol. 11, no.  4, pp. 1166--1198 (2001). 

  \bibitem{dm-2007}
  P.~Del~Moral, L.~Miclo.
  Strong propagations of chaos in Moran's type particle interpretations of Feynman-Kac measures.   Stochastic Analysis and Applications, vol. 25, no.  3, pp. 519--575 (2007). 



\bibitem{dm-stab}
P. Del Moral, L. Miclo.
On the Stability of Non Linear Semigroup of Feynman-Kac Type. 
Ann. de la Facult\'e des Sciences de Toulouse, vol.11, no.2, pp. (2002).

\bibitem{dm-sch}
P. Del Moral, L. Miclo.  Particle approximations of Lyapunov exponents connected to Schr\"odinger operators and Feynman-Kac semigroups. ESAIM: Probability and Statistics, no. 7, pp. 171--208 (2003).



\bibitem{dm-penev}
P. Del Moral, S. Penev.  Stochastic Processes: From Applications to Theory. Chapman and Hall/CRC (2016).

\bibitem{dm-patras-ruben}
P. Del Moral, F. Patras, S. Rubenthaler. Convergence of U-statistics for interacting particle systems. Journal of Theoretical Probability, vol. 24, no. 4:1002  (2011).

\bibitem{dm-ssd-19}
P. Del Moral, S.S. Singh 
A forward-backward stochastic analysis of diffusion flows.
\href{https://hal.inria.fr/hal-02161914/document}{Hal-02161914} (2019). 

\bibitem{dm-villemonais}
P. Del Moral, D. Villemonais.  Exponential mixing properties for time inhomogeneous diffusion processes with killing. Bernoulli, vol. 24, no. 2, pp. 1010--1032 (2018).

\bibitem{delyon-cerou-18}
B. Delyon, F. C\'erou, A. Guyader, M. Rousset. A central limit theorem for Fleming-Viot particle systems with hard killing. arXiv preprint arXiv:1709.06771 (2017).

\bibitem{ddg}
A. Doucet, N. De Freitas, N. Gordon, Editors. An introduction to sequential Monte Carlo methods. In Sequential Monte Carlo methods in practice. Springer, New York, NY (2001).

\bibitem{dupire}
B. Dupire. Functional it\^o calculus.  Bloomberg Portfolio Research Paper No. 2009--04--Frontiers (2009).


\bibitem{durmus}
A. Durmus, E. Moulines. Nonasymptotic convergence analysis for the unadjusted Langevin algorithm. In: Ann. Appl. Probab., vol.  27, no. 3 , pp. 1551--1587
(2017).


\bibitem{fearnhead}
P. Fearnhead, K. Latuszynski, G. O. Roberts, G. Sermaidis. Continuous-time importance sampling: Monte Carlo methods which avoid time-discretization error. arXiv preprint arXiv:1712.06201 (2017).

\bibitem{gillespie}
D. T. Gillespie. Exact stochastic simulation of coupled chemical reactions. J. Phys. Chem., vol. 81, no. 25, pp. 2340--2361 (1977).

\bibitem{yuri}
S. Jazaerli, Y.F. Saporito. Functional It\^o calculus, path-dependence and the computation of Greeks. Stochastic Processes and their Applications, vol. 127, no.12, pp. 3997--4028  (2017).



\bibitem{kloeden}
P.E. Kloeden, E. Platen. Numerical Solution of Stochastic Differential Equations. Springer, Berlin (1992).

\bibitem{kubo}
R. Kubo. The fluctuation-dissipation theorem. Reports on progress in physics, vol. 29, no. 1, p. 255 (1966).

\bibitem{kyprianou}
A.E. Kyprianou, S. Palau.  Extinction properties of multi-type continuous-state branching processes. Stochastic Processes and their Applications, vol. 128, no.10, pp.  3466--3489 (2018).


\bibitem{lande}
R. Lande, S. Engen, B.E. Saether.  Stochastic population dynamics in ecology and conservation. Oxford University Press (2003).

\bibitem{langevin}
P. Langevin. Sur la th\'eorie du mouvement brownien, Comptes-Rendus de l'Acad/'emie des Sciences, vol. 146, pp. 530--532 (1908).


\bibitem {tony}
T. Leli\`evre, M. Rousset and G. Stoltz. \emph{Free energy
computations: A mathematical perspective}. Imperial College Press (2010).

\bibitem{tony-rs}
T. Leli\`evre, M. Rousset and G. Stoltz. Computation of free energy differences through non-equilibrium stochastic dynamics: the reaction coordinate case, Journal of
Computational Physics, 222(2), 624-643, (2007).


\bibitem{lindsten}
F. Lindsten, T. Sch\"on, M. I. Jordan.
Ancestor Sampling for Particle Gibbs.
{\em Conference in Advances in Neural Information Processing Systems}, vol. 25, pp. 2600--2608 (2012).

\bibitem{lindsten2}
F. Lindsten, T. Sch\"on.
On the use of backward simulation in the particle Gibbs sampler.
Acoustics, Speech and Signal Processing (ICASSP), 2012 IEEE International Conference on. IEEE (2012).	
\bibitem{malrieu}
F. Malrieu. Logarithmic Sobolev inequalities for some nonlinear PDE's. Stochastic Process. Appl. vol. 95,  no. 1, 109--132 (2001).

\bibitem{malrieu-2}
F.  Malrieu. Convergence to equilibrium for granular media equations and their Euler schemes. Ann. Appl. Probab., vol. 13, no. 2, pp. 540--560 (2003).



 \bibitem{mckean-1}
 H. P. McKean. A class of Markov processes associated with nonlinear parabolic equations. Proc. Nat. Acad. Sci. U.S.A., vol. 56, pp. 1907--1911 (1966).

 \bibitem{mckean-2}
 H. P. McKean. Propagation of chaos for a class of non-linear parabolic equations. In Stochastic Differential Equations, Lecture Series in Differential Equations, Session 7, Catholic Univ. (1967), pp. 41--57. Air Force Office Sci. Res., Arlington, Va., (1967).
 

\bibitem{mel}
S. M\'el\'eard. Asymptotic behaviour of some interacting particle systems, McKean-Vlasov and Boltzmann models. In: Talay, D., Tubaro, L. (Eds.), Probabilistic Models for Nonlinear Partial Di erential Equations, Montecatini Terme, 1995, Lecture Notes in Mathematics, Vol. 1627. Springer, Berlin
(1996).

\bibitem{moran-1}
P.A.P.  Moran. The Statistical Processes of Evolutionary Theory. Oxford: Clarendon Press (1962).

\bibitem{moran-2}
P.A.P.  Moran. Random processes in genetics. Mathematical Proceedings of the Cambridge Philosophical Society. 54 (1): pp. 60--71 (1958). 


\bibitem{nash} 
J. Nash, Continuity of solutions of parabolic and elliptic equations, American J. Math. vol. 80, pp 931--954 (1958).

\bibitem{otto}
F. Otto. The geometry of dissipative evolution equations: the porous medium equation. Comm. Partial Differential Equations, vol. 26, pp. 101--174 (2001).

\bibitem{otto-2}
F. Otto, V. Villani. Generalization of an inequality by Talagrand, and links with the logarithmic Sobolev inequality. J. Funct. Anal. vol. 
173, pp. 361--400 (2000).


\bibitem{revuz}
D. Revuz. Markov chains. North Holland (2001). 
 
\bibitem{roberts}
G.O. Roberts, J.S. Rosenthal. Optimal scaling for various Metropolis-Hastings algorithms. Statistical science, vol. 16, no. 4, pp. 351--367 (2001). 
 
 
\bibitem{mathias}
M. Rousset. On the control of an interacting particle approximation of Schr\" odinger ground states. SIAM J. Math. Anal., vol. 38, no. 3, pp. 824--844 (2006).

\bibitem{mathias-2}
M.Rousset and G.Stoltz. Equilibrium sampling from non equilibrium dynamics. J. Stat. Phys., vol. 123, no. 6, pp. 1251--1272  (2006).

\bibitem{rosen}
M.N. Rosenbluth, A.W. Rosenbluth. Monte Carlo calculation of the average extension of molecular chains. The Journal of Chemical Physics, 23(2), pp. 356--359
(1955).



\bibitem{yuri-2}
Y.F. Saporito. Topics on Functional It\^o Calculus and Multiscale
Stochastic Volatility Modeling. PhD. thesis Santa Barbara (2014).

\bibitem{singer}
R.A. Singer. Estimating optimal tracking filter performance for manned maneuvering targets. IEEE Transactions on Aerospace and Electronic Systems, vol. 4, pp. 473--483 (1970). 

\bibitem{tamura}
Y. Tamura. On asymptotic behaviors of the solution of a nonlinear diffusion equation. J. Fac. Sci. Univ. Tokyo Sect. IA Math., vol.31, no. 1, pp. 195--221 (1984).

\bibitem{tamura-2}
Y. Tamura. Free energy and the convergence of distributions of diffusion processes of McKean type. J. Fac. Sci. Univ. Tokyo Sect. IA Math., vol. 34, no. 2, pp. 443--484, (1987).


\bibitem{varopoulos}
N. Th. Varopoulos. Small time gaussian estimates of heat diffusion kernels. II. The theory of large deviations.
Journal oh Functional Analysis, vol. 93, pp. 1-33 (1990).

\bibitem{denis-14}
D. Villemonais. General approximation method for the distribution of Markov processes conditioned not to be killed. ESAIM Probab. Stat., vol. 18, pp.
441--467 (2014).

\end{thebibliography}
\end{document}